\documentclass[a4paper,11pt]{article}
\usepackage{geometry,graphicx, amsmath, amsthm, multirow,algorithm,algorithmic,amssymb,lscape,array,listings,textcomp,hyperref,color, subfigure, psfrag}

\begin{document}
\title{Automated parameters for troubled-cell indicators using outlier detection}

\author{Mathea J. Vuik\footnote{Email:  M.J.Vuik@tudelft.nl.  Delft Institute of Applied Mathematics, Delft University of Technology, Mekelweg 4, 2628CD Delft, The Netherlands.}
\renewcommand{\thefootnote}{\fnsymbol{footnote}}
and
 \hspace{.05in}  Jennifer~K. Ryan\footnote{Corresponding Author. Email:  Jennifer.Ryan@uea.ac.uk. Telephone: +44 (0)1603 592586.   School of Mathematics, University of East Anglia, Norwich NR4 7TJ, United Kingdom.  Supported by the Air Force Office of Scientific Research, Air Force Material Command, USAF, under grant number FA8655-09-1-3055.}
}

\maketitle

\begin{abstract}
In Vuik and Ryan (2014) we studied the use of troubled-cell indicators for discontinuity detection in nonlinear hyperbolic partial differential equations and introduced a new multiwavelet technique to detect troubled cells. We found that these methods perform well as long as a suitable, problem-dependent parameter is chosen. This parameter is used in a threshold which decides whether or not to detect an element as a troubled cell. Until now, these parameters could not be chosen automatically. The choice of the parameter has impact on the approximation: it determines the strictness of the troubled-cell indicator. An inappropriate choice of the parameter will result in detection (and limiting) of too few or too many elements. The optimal parameter is chosen such that the minimal number of troubled cells is detected and the resulting approximation is free of spurious oscillations. 

In this paper we will see that for each troubled-cell indicator the sudden increase or decrease of the indicator value with respect to the neighboring values is important for detection. Indication basically reduces to detecting the outliers of a vector (one dimension) or matrix (two dimensions). This is done using Tukey's boxplot approach to detect which coefficients in a vector are straying far beyond others (Tukey, 1977). We provide an algorithm that can be applied to various troubled-cell indication variables. Using this technique the problem-dependent parameter that the original indicator requires is no longer necessary as the parameter will be chosen automatically.
\end{abstract}



\pagestyle{myheadings}
\thispagestyle{plain}
\markboth{}{M.~J. Vuik and J.~K. Ryan, Automated parameters via outlier detection}
\section{Introduction}
In \cite{Vui14R}, we studied the use of troubled-cell indicators for discontinuity detection in nonlinear hyperbolic partial differential equations and introduced a new multiwavelet technique to detect troubled cells. We compared the troubled-cell indicator of Qiu et al. using Harten's subcell resolution \cite{Har89,Qiu05S_2}, the KXRCF shock detector \cite{Kri04XRCF} and the multiwavelet troubled-cell indicator \cite{Vui14R}. We found that these methods perform well as long as a suitable, problem-dependent parameter is chosen, which was also observed in \cite{Qiu05S_2}. This parameter is used in a threshold which decides whether or not to detect an element as a troubled cell. Until now, these parameters could not be chosen automatically such that the indicator works well in a variety of situations \cite{Zhu13CQ}. Similarly, a parameter is required for adaptive mesh refinement \cite{Ger15IMMS}. Here, the used threshold parameter does depend on the discretization and the order of accuracy. This parameter is always chosen in the same way, similar to the KXRCF method \cite{Kri04XRCF}.

The choice of the parameter has impact on the approximation: it determines the strictness of the troubled-cell indicator. An inappropriate choice of the parameter will result in detection (and limiting) of too few or too many elements. Detection of too few elements leads to spurious oscillations, since not enough elements are limited. If too many elements are detected, then the limiter is applied too often, and therefore the method is more costly and the approximation smooths out after a long time. The optimal parameter is chosen such that the minimal number of troubled cells is detected and the resulting approximation is free of non-physical spurious oscillations. In general, many tests are required to obtain this optimal parameter for each problem \cite{Qiu05S_2, Zhu13CQ}. 

In this paper we will see that for each troubled-cell indicator the sudden increase or decrease of the indicator value with respect to the neighboring values is important for detection. Indication basically reduces to detecting the outliers of a vector (one dimension) or matrix (two dimensions). This is done using Tukey's boxplot approach to detect which coefficients in a vector are straying far beyond others \cite{Tuk77}, which is commonly used in statistical analysis \cite{Fri89HI, Vel81H}. This method is designed in such a way that only a few 'false positives' are found if the data are well behaved (i.e., Gaussian \cite{Hoa83MT}). Another advantage of this method is that it is not necessary to specify the number of possible outliers in advance. This is in contrast to many standard outlier-detection techniques which require a statement of the exact or the maximum number of outliers that may be present \cite{Hoa86IT}. 
    
We provide an algorithm that can be applied to various troubled-cell indication variables. Using this technique the problem-dependent parameter that the original indicator requires is no longer necessary as the parameter will be chosen automatically.

The numerical results in this paper are computed using the discontinuous Galerkin (DG) method \cite{Coc89S, Coc89LS, Coc90HS, Coc98S} together with a third-order strong stability-preserving Runge-Kutta time-stepping scheme \cite{Got98S}. We apply either the original troubled-cell indicators (with an optimal parameter), or the outlier-detection technique in combination with the indication variable. In that way, the performance of the new technique can be easily compared to the current state-of-the-art methods. We will apply the techniques to various test problems in one and two dimensions. Here, we investigate the modified multiwavelet troubled-cell indicator \cite{Vui15R}, the KXRCF shock detector \cite{Kri04XRCF} and the minmod-based TVB indicator \cite{Coc89S} in more detail. The moment limiter is used in the detected troubled cells \cite{Kri07}, but other limiting techniques can be used.

The outline of this paper is as follows:  in \S \ref{sec:background} we present the relevant background information on discontinuous Galerkin methods, troubled-cell indicators and the moment limiter. In \S \ref{sec:outlier} we introduce our new outlier-detection algorithm. The effectivity of this new method compared with the corresponding parameter-using troubled-cell indicators is presented in \S \ref{sec:results} for standard numerical examples. The computational costs of our algorithm are addressed in \S \ref{sec:costs}. We conclude with a discussion of our method and future work in \S \ref{sec:conclude}.

\section{Background}\label{sec:background}
This section contains some background information about the discontinuous Galerkin (DG) method \cite{Coc89S, Coc89LS, Coc90HS, Coc98S},  as well as the theory behind troubled-cell indicators  \cite{Kri04XRCF, Coc89S, Qiu05S_2, Vui14R, Vui15R} and the moment limiter \cite{Kri07}. This information will be used to apply our new outlier-detection scheme (\S \ref{sec:outlier}) in numerical examples. 
\subsection{Discontinuous Galerkin method}
We briefly explain the DG method using the initial-value problem
\setlength\arraycolsep{2pt}
\begin{equation}
\begin{array}{rlll}
  u_t + f(u)_x \!\! &= 0, & x\in[-1,1],  & t > 0, \\
u(x,0)\!\! &=u^0(x), &  x \in [-1,1], &
\end{array}
\label{eqn:linadvsystem}
\end{equation}
where $u = u(x,t)$, and $f(u)$ describes the flux function. 

Discretization in space is obtained by dividing $[-1,1]$ into $2^n$ elements (used in the multiwavelet expansion, \S \ref{sec:modmul}), defined as
\[
I_j = [x_{j-\frac{1}{2}}, x_{j+\frac{1}{2}}), \ j = 0,\ldots,2^n-1. 
\]
The choice for half-open intervals follows from the paper of Archibald et al. \cite{Arc11FS}. Different choices are available in the literature, for example closed intervals (Hovhannisyan et al. \cite{Hov10MS}), or open intervals (Gerhard et al. \cite{Ger15IMMS}).

The approximation space that we use on each element is $V_h(I_j) = \{ v \in \mathbb{P}^k(I_j) \},$ where $\mathbb{P}^k$ is the space of polynomials of degree $k$. In order to take advantage of the  multiwavelet properties, the basis for $\mathbb{P}^k$ is constructed using the scaled Legendre polynomials, which are defined as
\begin{equation}
 \phi_\ell(x) = \sqrt{\ell + \frac{1}{2}} P^{(\ell)}(x),\label{eqn:Legendre}
\end{equation}
where $P^{(\ell)}$ is the Legendre polynomial of degree $\ell$, $\ell = 0,\ldots, k$. 

The weak form of the PDE in problem (\ref{eqn:linadvsystem}) is obtained by multiplying the equation by a test function $v \in V_h(I_j)$ and integrating over element $I_j$. Using integration by parts, this yields
\begin{equation}
\int_{I_j} u_tvdx = \int_{I_j}f(u)v_xdx + \hat{f}_{j-\frac{1}{2}}v_{j-\frac{1}{2}}^+ - \hat{f}_{j+\frac{1}{2}}v_{j+\frac{1}{2}}^-, \label{eqn:weakformflux}
\end{equation}
where $\hat{f}_{j\pm\frac{1}{2}}$ denote the flux values through boundaries $x_{j\pm1/2}$. These are approximated using the local Lax Friedrichs flux \cite{LeV02}:
\[
 \hat{f}_{j-\frac{1}{2}} = \frac{1}{2} \left( f(u_{j-\frac{1}{2}}^-) + f(u_{j-\frac{1}{2}}^+) - a_{j-\frac{1}{2}}(u_{j-\frac{1}{2}}^+ - u_{j-\frac{1}{2}}^-) \right),
\]
where we use that $f$ is convex, such that
\[
 a_{j-\frac{1}{2}} = \max(|f^\prime(u_{j-\frac{1}{2}}^-)|, |f^\prime(u_{j-\frac{1}{2}}^+)|).
\]
The third-order strong stability-preserving Runge-Kutta scheme \cite{Got98S} is used for time evolution. Note that this is only a choice, and that other time-stepping schemes are also possible \cite{Got01ST, Ket09MG, Shu88}.

\subsection{Troubled-cell indicators}\label{sec:troubled}
In this section, various troubled-cell indicators are described on which our new outlier-detection algorithm will be tested. In particular, we will investigate the modified multiwavelet troubled-cell indicator \cite{Vui14R, Vui15R}, the KXRCF indicator \cite{Kri04XRCF}, and the minmod-based TVB indicator \cite{Coc89S}. In earlier papers, Harten's subcell resolution idea \cite{Har89} was used for indication \cite{Qiu05S_2, Vui14R}. However, this method was unstable in several numerical experiments \cite{Zhu09Q}, and therefore, we will not investigate this method here.

The outlier-detection algorithm will require that we pass in a vector of troubled-cell indication variables and therefore we provide the vector form as well.

\subsubsection{Modified multiwavelet troubled-cell indicator}\label{sec:modmul}
In \cite{Vui14R}, a multiwavelet troubled-cell indicator was constructed, which was modified in \cite{Vui15R} by using multiwavelet coefficients instead of multiwavelet contributions. In this section, we repeat the important definitions. Here, we only investigate the domains $[-1,1]$ (one dimension) and $[-1,1]\times[-1,1]$ (two dimensions). The corresponding definitions can be easily extended to general domains in one and two dimensions \cite{Vui14R}. 

A global DG approximation of degree $k$ on $2^n$ elements in $[-1,1]$ can be written as
\[
  u_h(x) = 2^{-\frac{n}{2}}\sum_{j=0}^{2^n - 1} \sum_{\ell=0}^k u_j^{(\ell)} \phi_{\ell j}^n(x), \quad x \in [-1,1],
\]
where the scaling functions $\phi_{\ell j}^n$ are defined as
\[
\phi_{\ell j}^n(x) = 2^{n/2} \phi_\ell(2^n(x+1)-2j-1), \quad \ell = 0,\ldots, k, \ j = 0,\ldots, 2^n - 1,
\]
and $\phi_{\ell}$ are the scaled Legendre polynomials as in equation (\ref{eqn:Legendre}).

The corresponding multiwavelet decomposition of the DG approximation can be written as
\[
 u_h(x) =  \sum_{\ell=0}^{k}s_{\ell0}^0 \phi_\ell(x) + \sum_{m=0}^{n-1} \sum_{j = 0}^{2^m - 1} \sum_{\ell = 0}^k d_{\ell j}^m\psi_{\ell j}^m(x),
\]
where $s_{\ell 0}^0$ are the scaling-function coefficients belonging to $u_h$, the multiwavelets on higher levels are defined as $\psi_{\ell j}^m(x) = 2^{m/2} \psi_\ell(2^m(x+1)-2j-1)$, and $d_{\ell j}^m$ are the corresponding multiwavelet coefficients \cite{Arc11FS, Vui14R}, which are determined using the orthogonal projection of the DG approximation onto the multiwavelet basis:
\[
 d_{\ell j}^m = \langle u_h, \psi_{\ell j}^m \rangle_{L_2([-1 + 2^{-m+1}j , -1 + 2^{-m+1}(j + 1) ])}. 
\]
In practice, these coefficients are efficiently computed using the quadrature mirror filter coefficients \cite{Alp02BGV, Arc11FS}. The multiwavelets $\psi_{\ell}$ have been developed by Alpert \cite{Alp93}, and are also explained in \cite{Hov10MS}. 

As we have seen in \cite{Vui15R}, the coefficients on level $n-1$ are strongly related to the inter-element jumps in (the derivatives of) the DG approximation. The multiwavelet coefficients on level $n-1$ of the decomposition equal
\begin{subequations}
\begin{equation}
  d_{\ell j}^{n-1} = 2^{-\frac{n-1}{2}} \sum_{m=0}^k c_{m,\ell}^n \cdot \left( u_h^{(m)}(x_{2j+1/2}^+) - u_h^{(m)}(x_{2j+1/2}^-) \right),\label{eqn:dsecondlast}
\end{equation}
with
\begin{equation}
 c_{m,\ell}^n = \frac{2^{(-n+1)m}}{m!}  \cdot \int_0^1 x^m \psi_\ell(x) \, dx,\label{eqn:c}
\end{equation}\label{eqn:defd}\end{subequations}where $\ell = 0,\ldots, k$, $j = 0,\ldots, 2^{n-1}-1$, and $u_h^{(m)}$ is the $m$th derivative of $u_h$. Note that this only includes half of the element-boundary jumps. In order to also compute the rest of the jumps, we use the same renumbering technique as was proposed in \cite{Vui15R}. This gives rise to $2^n$ coefficients for level $n-1$, in the following denoted by $\tilde{d}_{\ell j}^{n-1}, \ \ell = 0,\ldots,k, j = 0,\ldots,2^n-1$.

\bigskip
In general, the DG approximation is discontinuous at element boundaries. Therefore, the multiwavelet coefficients are usually not exactly equal to zero. However, when the solution is sufficiently smooth, then the element-boundary jumps in the approximation and its derivatives will be noticeably smaller than when a discontinuity in (one of the derivatives of) the solution is present due to the cancellation property of multiwavelets \cite{Hov14MS}. The multiwavelet coefficients $\tilde{d}_{kj}^{n-1}$ are used to detect troubled cells when 
\begin{equation}
 |\tilde{d}_{kj}^{n-1}| > C \cdot \max \{ |\tilde{d}_{kj}^{n-1}|, j=0,\ldots,2^n-1 \}, \ C \in [0,1]. \label{eqn:mworiginal}
\end{equation}Since coefficient $\tilde{d}_{kj}^{n-1}$ contains information about the jump in (derivatives of) the DG approximation at $x_{j+1/2}$, elements $I_j$ and $I_{j+1}$ are limited if $\tilde{d}_{kj}^{n-1}$ satisfies inequality (\ref{eqn:mworiginal}).

Note that only the multiwavelet coefficients on level $n-1$ are used for indication. Therefore, a twoscale representation of the DG approximation would suffice, and the number of elements in the domain might be even (instead of the restriction to a power of two). However, since this would change the definitions, we have chosen to use $2^n$ elements in this work.

\bigskip
The boundary of which $|\tilde{d}_{kj}^{n-1}|$ is maximal, is assumed to be the location where the strongest shock occurs. If $C=1$, then no element will be detected, and the smaller $C$ is, the more elements will be limited. In this way, the value of $C$ is a useful tool to prescribe the strictness of the limiter. In general, it is hard to choose a sufficient value for $C$. For each problem, several tests should be done in order to obtain an optimal parameter \cite{Vui14R}.

\bigskip
In order to remove the problem-dependent parameter $C$ that occurs in this indicator, we propose to use the multiwavelet coefficients in our outlier-detection algorithm (\S \ref{sec:outlier}). The indication vector is defined as $\mathbf{D} = (\tilde{d}_{k,0}^{n-1},\ldots,\tilde{d}_{k,2^n-1}^{n-1})^\top$.

\bigskip
In two dimensions, the relations for the multiwavelet coefficients on level $\mathbf{n}-1$ follow naturally from the one-dimensional coefficients:
\begin{align*}
 d_{\boldsymbol{\ell}, \mathbf{j}}^{\alpha, \mathbf{n}-1} &= 2^{- \frac{n_y - 1}{2}}\sum_{m_y=0}^k c_{m_y,\ell_y}^{n_y} \int_{x_{2i-\frac{1}{2}}}^{x_{2i+\frac{3}{2}}} 
\hspace{-0.5mm} \left(\frac{\partial^{m_y} u_h}{\partial y^{m_y}}(x,y_{2j+\frac{1}{2}}^+) -\frac{\partial^{m_y} u_h}{\partial y^{m_y}}(x,y_{2j+\frac{1}{2}}^-)\right) \phi_{\ell_x, i}^{n_x - 1}(x) dx, \\
 d_{\boldsymbol{\ell}, \mathbf{j}}^{\beta, \mathbf{n}-1} &= 2^{- \frac{n_x - 1}{2}}\sum_{m_x=0}^k c_{m_x,\ell_x}^{n_x} \int_{y_{2j-\frac{1}{2}}}^{y_{2j+\frac{3}{2}}} 
\left(\frac{\partial^{m_x} u_h}{\partial x^{m_x}}(x_{2i+\frac{1}{2}}^+,y) -\frac{\partial^{m_x} u_h}{\partial x^{m_x}}(x_{2i+\frac{1}{2}}^-,y)\right) \phi_{\ell_y, j}^{n_y - 1}(y) dy,\\
d_{\boldsymbol{\ell}, \mathbf{j}}^{\gamma, \mathbf{n}-1}& = 2^{- \frac{n_x - 1}{2}} 2^{- \frac{n_y - 1}{2}} \sum_{m_x=0}^k \sum_{m_y=0}^k c_{m_x,\ell_x}^{n_x} c_{m_y,\ell_y}^{n_y} \cdot \left(\frac{\partial^{m_x}}{\partial x^{m_x}} \frac{\partial^{m_y}}{\partial y^{m_y}}\right.\\
&\left.  \left( u_h(x_{2i+\frac{1}{2}}^+, y_{2j+\frac{1}{2}}^+) - u_h(x_{2i+\frac{1}{2}}^+, y_{2j+\frac{1}{2}}^-) - u_h(x_{2i+\frac{1}{2}}^-, y_{2j+\frac{1}{2}}^+) + u_h(x_{2i+\frac{1}{2}}^-, y_{2j+\frac{1}{2}}^-)\right) \right),
\end{align*}
where $\boldsymbol{\ell} = (\ell_x, \ell_y)^\top$, $\mathbf{j} = (i,j)^\top$, and $\mathbf{n} = (n_x, n_y)^\top$, $i=0,\ldots,2^{n_x-1} - 1$ and $j=0,\ldots,2^{n_y-1}-1$, and $c_{m,\ell}^n$ is defined as in equation (\ref{eqn:c}). Using the renumbering technique as mentioned before \cite{Vui15R}, we include the second half of the boundaries in the multiwavelet direction, and find coefficients $\tilde{d}_{\boldsymbol{\ell}, \mathbf{j}}^{\alpha, \mathbf{n}-1}$ ($j = 0, \ldots, 2^{n_y}-1$), $\tilde{d}_{\boldsymbol{\ell}, \mathbf{j}}^{\beta, \mathbf{n}-1}$ ($i=0,\ldots,2^{n_x}-1$) and $\tilde{d}_{\boldsymbol{\ell}, \mathbf{j}}^{\gamma, \mathbf{n}-1}$ ($i=0,\ldots,2^{n_x}-1, \ j = 0, \ldots, 2^{n_y}-1$).

Note that these relations indeed confirm the observations that the $\alpha$ mode detects discontinuities in the $y$-, the $\beta$ mode in the $x$-, and the $\gamma$ mode in the $xy$-direction, as was stated in \cite{Mal98} and seen in  \cite{Vui14R}. 

In two dimensions, the approach of inequality (\ref{eqn:mworiginal}) is applied for each mode separately. In the $\alpha$ mode, we take the coefficients with index $\boldsymbol{\ell} = (0,k)^\top$ for indication. In the $\beta$ mode, the indices $\boldsymbol{\ell} = (k,0)^\top$ are used, and for $\gamma$ we take $\boldsymbol{\ell} = (k,k)^\top$.

\bigskip
Using the directions of each mode, the one-dimensional outlier-detection algorithm is applied to the $\alpha$-mode vectors for each $x$, and to the $\beta$-mode vectors for each $y$. We have found that detection on the $\gamma$ mode selects too many elements. Therefore, this mode is not used in the outlier-detection scheme. We apply outlier detection to the following vectors:
\[
  \begin{array}{rll}
  \alpha: &\mathbf{D}_i = \left( \tilde{d}_{(0,k),(i,0)}^{\alpha, \mathbf{n}-1}, \ldots,\tilde{d}_{(0,k),(i,2^{n_y} - 1)}^{\alpha, \mathbf{n}-1} \right), & i=0,\ldots,2^{n_x-1}-1,\\
  \beta: & \mathbf{D}_j= \left( \tilde{d}_{(k,0),(0,j)}^{\beta, \mathbf{n}-1}, \ldots,\tilde{d}_{(k,0),(2^{n_x}-1,j)}^{\beta, \mathbf{n}-1} \right), & j=0,\ldots,2^{n_y-1} -1.\\
  \end{array}
\]
We note that the multiwavelet indicator is equivalent to comparing the DG approximation over two different spatial meshes. An alternative strategy would be to compare the approximation at two different levels in time as in Dumbser et al. \cite{Dum14ZLD}.

\subsubsection{KXRCF indicator}\label{sec:KXRCF}
The shock-detection technique by Krivodonova et al. \cite{Kri04XRCF} uses inflow boundaries to detect troubled cells. The detector considers the jump in $u_h$ across the inflow edges of $I_j$ and examines
\[
 \mathcal{I}_j = \left|\int_{\partial I_j^-} (u_h|_{I_j} - u_h|_{I_{n_j}}) ds\right|.
\]
Here, $\partial I_j^-$ is the inflow boundary and $u_h|_{I_{n_j}}$ is the DG approximation in the neighbor of $I_j$ on the side of $\partial I_j^-$. The indicator is normalized to
\begin{equation}
  \hat{\mathcal{I}}_j = \frac{\left| \int_{\partial I_j^-} (u_h|_{I_j} - u_h|_{I_{n_j}}) ds \right|}{h^{\frac{k+1}{2}} |\partial I_j^-| ||u_h|_{I_j}||}, \quad j = 0, \ldots, 2^n-1. \label{eqn:KXRCF}
\end{equation}
Here, $h$ is the radius of the circumscribed circle in $I_j$, and the norm is based on the average in one dimension and the maximum norm in quadrature points in two dimensions. 

Near a discontinuity $\hat{\mathcal{I}}_j \rightarrow \infty$, whereas $\hat{\mathcal{I}}_j \rightarrow 0$ if $h \rightarrow 0$ or $k \rightarrow \infty$ in smooth-solution regions. In \cite{Kri04XRCF}, the threshold value is taken equal to 1, such that element $I_j$ is detected as troubled if $\hat{\mathcal{I}}_j > 1$, and in that case the limiter is applied in $I_j$. Note that this threshold parameter is chosen arbitrarily: the value 1 does not necessarily follow from the theory.

\bigskip
In order to remove this parameter, the outlier-detection mechanism was tested on the vector $\mathbf{D} = (\hat{\mathcal{I}}_0, \ldots, \hat{\mathcal{I}}_{2^n-1})^\top$. However, it turned out that the original discontinuity detector without normalization is more suitable, such that the jump across the interfaces is used in the indicator: we take $\mathbf{D} = (\mathcal{I}_0, \ldots, \mathcal{I}_{2^n-1})^\top$ for detection. 

In two dimensions, a matrix $\mathbf{D} = \{\mathcal{I}_{ij} \}$ is found. Here, the one-dimensional outlier-detection approach is applied in the $x$- and $y$-direction separately (row and column wise).

\subsubsection{Minmod-based TVB indicator}\label{sec:minmod}
In this section, the minmod-based TVB indicator of Cockburn and Shu will be explained \cite{Coc89LS, Coc89S}. For each element $I_j, j = 0,\ldots, 2^n-1$, the element-boundary approximations are split into
\[
 u_{j+\frac{1}{2}}^- = \bar{u}_j + \tilde{u}_j \quad \mbox{and} \quad u_{j-\frac{1}{2}}^+ = \bar{u}_j - \tilde{\tilde{u}}_j,
\]
where
\begin{equation}
 \tilde{u}_j = \sum_{\ell=1}^k u_j^{(\ell)} \phi_\ell(1), \quad \tilde{\tilde{u}}_j = - \sum_{\ell=1}^k u_j^{(\ell)} \phi_\ell(-1).\label{eqn:minmodparts}
\end{equation}
Element $I_j$ is detected as troubled if either $\tilde{u}_j$ or $\tilde{\tilde{u}}_j$ is modified by the functions
\begin{equation}
  \tilde{u}_j^{(\mbox{\footnotesize{mod}})} = \tilde{m}(\tilde{u}_j, \bar{u}_{j+1} - \bar{u}_j, \bar{u}_j - \bar{u}_{j-1}), \quad  
 \tilde{\tilde{u}}_j^{(\mbox{\footnotesize{mod}})} = \tilde{m}(\tilde{\tilde{u}}_j, \bar{u}_{j+1} - \bar{u}_j, \bar{u}_j - \bar{u}_{j-1}), \label{eqn:minmoddetect}
\end{equation}
where the TVB-modified minmod function is defined as
\[
 \tilde{m}(a_1, \ldots, a_q) = \left\{ \begin{array}{ll} a_1, & \mbox{if } |a_1| \leq M\Delta x^2, \\ m(a_1, \ldots, a_q), & \mbox{otherwise,}  \end{array}
\right.
\]
in contrast with the standard minmod function
\begin{equation}
 m(a_1, \ldots, a_q) = \left\{ \begin{array}{ll} s \cdot \min_{1 \leq j \leq q} |a_j|, & \mbox{if } \mbox{sign} (a_1) = \cdots = \mbox{sign} (a_q) = s, \\ 0, & \mbox{otherwise.}  \end{array}
\right. \label{eqn:minmod}
\end{equation}
Note that the parameter $M$ is difficult to tune, and hardly any difference is found when $M$ ranges from 1 to 100, \cite{Zhu13CQ}. We use the minmod-based TVB indicator for detection and then apply a chosen limiter in the detected troubled cells.

For systems of equations, characteristic field decompositions are required \cite{Coc89LS}. The corresponding eigenvector matrix is computed using Roe averages \cite{Coc89LS, Roe81}.

\bigskip
Instead of using the parameter $M$, we will apply the outlier-detection algorithm. DG coefficients $u_j^{(1)},\ldots, u_j^{(k)}$ usually differ substantially from their neighbors when $I_j$ belongs to a discontinuous region. This means that we use the vectors $\mathbf{D}_1 = ( \tilde{u}_0,\ldots,\tilde{u}_{2^n-1})^\top$ and $\mathbf{D}_2 = ( \tilde{\tilde{u}}_0,\ldots,\tilde{\tilde{u}}_{2^n-1})^\top$ in the outlier-detection technique and detect element $I_j$ as troubled if either $\tilde{u}_j$ or $\tilde{\tilde{u}}_j$ is detected as an outlier.

\bigskip
For two-dimensional systems, the procedure for $\mathbb{P}^k$ has been explained in \cite{Coc98S}. The indicator uses solution derivatives (e.g. DG coefficients) for detection. We use $\mathbb{Q}^k$, which means that more 'cross-product' coefficients exist (for example, for $k=1$: $u_{ij}^{(1,1)}$). However, using Biswas's reasoning \cite{Bis94DF}, we do not use these coefficients for detection, since they have a lesser effect on the numerical approximation than either $u_{ij}^{(1,0)}$ or $u_{ij}^{(0,1)}$. 

The minmod-based TVB indicator in two dimensions resembles the two-dimensional moment limiter \cite{Kri07}. The difference between the two approaches is that the moment limiter uses forward and backward differences of lower derivatives, whereas the minmod-based indicator uses a finite-difference approach on the element averages.

In our numerical examples, we will focus on the case $k=1$, and use the DG coefficients $u_{ij}^{(1,0)}$ and $u_{ij}^{(0,1)}$ for detection. Outlier detection will be applied to the vectors
\[
 \mathbf{D}_j = \left(u_{0,j}^{(1,0)}, \ldots, u_{2^{n_x}-1, j}^{(1,0)} \right)^\top, \ j = 0,\ldots, 2^{n_y}-1,
\]
in the $x$-direction and
\[
 \mathbf{D}_i = \left(u_{i,0}^{(0, 1)}, \ldots, u_{i, 2^{n_y}-1}^{(0, 1)} \right)^\top, \ i = 0,\ldots, 2^{n_x}-1,
\]
in the $y$-direction. In this way, it is possible to detect discontinuities in different directions as we will see in \S \ref{sec:2d}.

\subsection{Moment limiter}\label{sec:moment}
In the detected troubled cells, a limiter is applied. The limiting technique that we use in this paper is the moment limiter \cite{Kri07}. This is only a choice - other limiters are also possible.

The moment limiter reduces the DG approximation to a low order in discontinuous regions, and maintains a high order if the approximation is smooth enough. Although the limiter has its own mechanism to control which regions should be limited, we will apply troubled-cell indicators as a switch to control where the limiter is applied. This is to prevent limiting smooth extrema. 

The moment limiter limits DG coefficients, starting at the highest level $k$. For each element $I_j, j=0,\ldots,2^n-1$, the limited value of coefficient $u_j^{(k)}$ equals 
\begin{equation}
 \widetilde{u}_j^{(k)} = m \left( u_j^{(k)},\beta_k\left(u_{j+1}^{(k-1)} - u_j^{(k-1)}\right),\beta_k\left(u_j^{(k-1)}-u_{j-1}^{(k-1)}\right) \right), \label{eqn:moment}
\end{equation}
with $\beta_k = (\sqrt{k-1/2})/(\sqrt{k+1/2})$ and using the minmod function (equation (\ref{eqn:minmod})). If $\widetilde{u}_j^{(k)} = u_j^{(k)},$ then the limiting procedure is cut off for this element $I_j$. If not, then $u_j^{(k-1)}$ is limited using the same procedure, continuing until $u_j^{(1)}$ is limited, or stopping the first time $\widetilde{u}_j^{(\ell)} = u_j^{(\ell)}$ for some $\ell = k-1,\ldots,1.$  

For systems of equations the limiter is applied to the characteristic variables $\mathbf{w}_j^{(\ell)} = R^{-1}\mathbf{u}_j^{(\ell)}$. Due to this approach it is possible that negative values for density, pressure or energy are found. In that case, all higher-order coefficients are set equal to zero, and $u_j^{(1)}$ is limited using equation (\ref{eqn:moment}). If negative values are still found, then the linear coefficient is also set equal to zero. 

In two dimensions, the moment limiter uses the neighboring elements both in the $x$-, and in the $y$-direction \cite{Kri07}.

\section{Troubled-cell indication using outlier detection}\label{sec:outlier}
In this section, an outlier-detection algorithm is proposed to detect outliers in a vector. This technique will be applied to the troubled-cell indicators given in \S \ref{sec:troubled}.

In order to detect outliers we use a boxplot mechanism that is often applied in statistics \cite{Fri89HI, Vel81H}, and described by Tukey \cite{Tuk77}. Important properties of this method are that only a few 'false positives' are found if the data are well behaved (i.e., Gaussian \cite{Hoa83MT}), and that it is not necessary to specify the number of possible outliers in advance. This is in contrast to many standard outlier-detection techniques which require to state the exact or the maximum number of outliers that may be present \cite{Hoa86IT}.

\bigskip
Here, we use a vector $\mathbf{d} = (d_0,\ldots,d_N)^\top$, of which outliers (suddenly changing coefficients with respect to neighbors) should be detected. A general outline of the outlier-detection algorithm that we use is provided below. In the following we discuss the details.

\begin{algorithm}
\caption{Outlier-detection algorithm.}
\begin{algorithmic}\label{alg:outlierplain}
\STATE Send in a suitable troubled-cell indication vector $\mathbf{d}$.
\STATE Sort $\mathbf{d}$ to obtain $\mathbf{d}^s$.
\STATE Compute the quartiles of $\mathbf{d}^s$.
\STATE Construct the outer fences.
\STATE Determine the outliers.
\end{algorithmic}
\end{algorithm}

\subsection{Quartiles}\label{sec:median}
Quartiles separate the data into four equal groups \cite{Qua}. The values of $Q_1$, $Q_2$ (the median), and $Q_3$ provide useful information about the structure of $\mathbf{d}$. As a preparation, it is convenient to sort $\mathbf{d}$, such that we obtain the vector $\mathbf{d}^s$: 
\[
 \mathbf{d}^s = (d_0^s, d_1^s, \ldots, d_N^s)^\top \mbox{, where, } d_0^s \leq d_1^s \leq \ldots \leq d_N^s.
\]
The median of $\mathbf{d}$ is defined as the numerical value that separates the higher half of the vector from the lower half \cite{Med}. It equals
\[
 \mbox{med}(\mathbf{d}) =  \left\{ \begin{array}{cl} d_{N/2}^s, & \mbox{if } N \mbox{ is even,} \\ \frac{1}{2}\left(d_{(N-1)/2}^s + d_{(N+1)/2}^s\right), & \mbox{if } N \mbox{ is odd.} \end{array} \right.
\]
The median is also called the second quartile of the vector $\mathbf{d}$.

\bigskip
The first quartile is defined as the value below which 25\% of the data fall, and is denoted by $Q_1$. Similarly, the third quartile, $Q_3$, equals the value that splits off the lowest 75\% of the data from the highest 25\% \cite{Qua}. Many different definitions of the first and third quartiles are used. In this work we apply Tukey's definition (definition 6 in \cite{Fri89HI}): 
\begin{equation}
  Q_1 = (1-g)d_{j-1}^s + gd_j^s,\label{eqn:percentiles}
\end{equation}
where $[(N+4)/2]/2 = j+g$, and $[x]$ denotes the largest integer that does not exceed $x$. Note that $g=0$ or $g=1/2$. The third quartile $Q_3$ is then computed symmetrically using the upper end of the vector $\mathbf{d}^s$.

In practice, we will always use a vector with $N+1=4r$ coefficients, where $r\in \mathbb{N}$. In that case, $Q_1 = (d_{r-1}^s + d_{r}^s)/2$ and $Q_3 = (d_{3r-1}^s+d_{3r}^s)/2$.

\subsection{Fences and outlier detection}
We have already seen that the values of the quartiles provide useful information about the structure of $\mathbf{d}^s$. However, this is not enough to define outliers in the vector. Outliers are the coefficients in the vector that are straying far out beyond the others. In order to pick out certain coefficients as outliers, inner and outer fences are constructed, which were originally defined by Tukey \cite{Tuk77}. The inner fences are equal to $[Q_1 - 1.5(Q_3 - Q_1), Q_3 + 1.5(Q_3 - Q_1)]$ (coefficients outside this interval are called \emph{soft outliers}). When the data are normally distributed, only 0.7\% of the data set is seen as a soft outlier (asymptotically) \cite{Hoa86IT}. The value 1.5 is referred to as the \emph{whisker length} of the boxplot.

The outer fences of a vector are $[Q_1 - 3(Q_3 - Q_1), Q_3 + 3(Q_3 - Q_1)]$ (coefficients outside are called \emph{extreme outliers}). The coverage for this whisker length is 99.9998\%, such that  only 0.0002\% of the data in a normally distributed vector is detected as an extreme outlier (asymptotically) \cite{Hoa86IT}. The choices of the whisker lengths (1.5 and 3) were proposed by Tukey \cite{Tuk77}, and are commonly used in the literature \cite{Fri89HI, Hoa86IT, Hub08V, Hun06Y, Sch07S, Sch04OA}. We will use the extreme outliers to detect troubled cells, since then very outstanding coefficients in the vector are selected. Because the data were sorted, the outer fences and outliers can easily be determined. 

\subsection{Application to troubled-cell indication variables}\label{sec:application}
In this section, we will explain the application of outlier detection to troubled-cell indication variables in one dimension. The corresponding indication vectors for each troubled-cell indicator were given in \S \ref{sec:troubled}. In this section, we have seen that all described indicators attach a value to each element of the domain (multiwavelet coefficient, jump across inflow boundary, or approximation at boundaries). Discontinuous regions usually correspond to the locations where the indicator value suddenly increases or decreases with respect to the neighboring values. This means that indication basically reduces to detecting the outliers of a vector with troubled-cell indication values. By applying the new outlier-detection technique, the threshold to be an extreme outlier is fixed, and the indicator no longer depends on problem-dependent parameters.

\bigskip
When an approximation contains several discontinuous regions, outlier detection applied to the global vector $\mathbf{D}$ will only select the strongest discontinuities. In order to also take into account the weaker discontinuities and the local structure of the approximation, the vector $\mathbf{D}$ will be split into local vectors of fixed length. For each subvector the outlier-detection mechanism is applied. In the local approach we ignore the detected coefficients in the left half of the local region if they are not detected with respect to the left-neighboring vector, and similarly the detected coefficients in the right half of the local region are tested. In this way the spatial information can still be used.

The outlier-detection algorithm executes the steps as provided in Algorithm \ref{alg:outlier}. Below we explain this in more detail.

\begin{algorithm}
\caption{Outlier-detection algorithm using local vectors.}
\begin{algorithmic}\label{alg:outlier}
\STATE Send in a suitable troubled-cell indication vector $\mathbf{D}$.
\STATE Split this vector into local vectors, $\mathbf{d}$.
\FORALL{local vectors} \STATE Sort $\mathbf{d}$ to obtain $\mathbf{d}^s$.
\STATE Compute $Q_1$ and $Q_3$ using definition (\ref{eqn:percentiles}).
\STATE Detect $d_j^s$ in the smallest 25\% of $\mathbf{d}^s$ if $d_j^s < Q_1 - 3(Q_3 - Q_1)$, and $d_j^s$ in the biggest 25\% of $\mathbf{d}^s$ if $d_j^s > Q_3 + 3(Q_3 - Q_1)$. 
\ENDFOR
\STATE Ignore the detected outliers in the left half of the local region when they are not detected with respect to the left-neighboring vector, and similarly test the detected coefficients in the right half of the local region.
\end{algorithmic}
\end{algorithm}

Since the global vector $\mathbf{D}$ consists of $2^n$ coefficients, we propose to split $\mathbf{D}$ into $2^{n-p}$ local vectors of length $2^p$, where $p \in \{ 2,\ldots,n \}$. Each local vector is then sorted. For convenience, we denote the sorted local vector of coefficients by $\mathbf{d}^{s} = (d_0^s, d_1^s, \ldots d_{N}^s)$, where $N = 2^p-1$. By definition this vector has the following 25th and 75th percentiles (see equation (\ref{eqn:percentiles})):
\[
Q_1 = \frac{d_{2^{p-2}-1}^s + d_{2^{p-2}}^s}{2}, \quad Q_3 = \frac{d_{3\cdot 2^{p-2}-1}^s + d_{3 \cdot 2^{p-2}}^s}{2}.
\]
Next we compute outer fences. Outliers are determined by comparing the smallest vector values with $Q_1 - 3(Q_3-Q_1)$ and the biggest components with $Q_3+3(Q_3-Q_1)$. For the smallest values we start with testing whether $d_0^s < Q_1 - 3(Q_3-Q_1)$.  If $d_0^s$ is not an outlier, then there are no other outliers, since $d_j^s \geq d_0^s \geq Q_1 -3(Q_3-Q_1)$, $j=0,\ldots,N$. If $d_0^s$ is an outlier, then we test $d_1^s$, etcetera. By construction, $Q_1 - 3(Q_3-Q_1) \leq Q_1$, such that the only possibilities for low outliers are $d_0^s,\ldots, d_{2^{p-2}-2}^s$ ($2^{p-2}-1$ coefficients). This means that at most $d_0^s, d_1^s, \ldots, d_{2^{p-2}-2}^s$ should be tested.

Similarly we test $d_{N}^s$ and (possibly) $d_{N-1}^s, \ldots, d_{3\cdot 2^{p-2}+1}^s$ against $Q_3+3(Q_3-Q_1)$ (depending on the outcome). Also here, at most $2^{p-2} -1$ coefficients should be tested. 

Finally, the detected outliers in the left half of the local vector are compared with the fences of the left-neighboring vector, and the outliers in the right half are compared with the right-neighboring fences.

\medskip
Considering the number of elements in each local vector, it should be noticed that $p = 3$ results in 8 coefficients per vector, which is too few to find a boxplot which is meaningful. Using $p=4$ (16 coefficients per vector) means that at maximum six outliers can be detected per local vector. Therefore, the maximum number of possible outliers in $\mathbf{D}$ equals $2^{n-4} \cdot 6 = 3 \cdot 2^{n-3}$. If we take more coefficients per local vector, for example $p = 5$ (32 coefficients per vector), then the 'stencil' is too big to extract all local information of the approximation. Therefore, we propose to use 16 coefficients per local vector ($p=4$), which worked well in all test cases we performed.

\medskip
In two dimensions, the one-dimensional algorithm is applied in the $x$- and $y$-direction separately. The corresponding troubled-cell indication vectors were given in \S \ref{sec:troubled}.

\section{Numerical results}\label{sec:results}
In this section, the original troubled-cell indicators are compared with the new outlier-detection approaches. This is done for the modified multiwavelet troubled-cell indicator of Vuik and Ryan \cite{Vui15R}, the KXRCF indicator \cite{Kri04XRCF}, and the minmod-based TVB indicator \cite{Coc89S}. We computed the results using $k=1,2,3$. In this paper, we only present the case $k=2$.

The results for the one-dimensional test cases are presented using time-history plots of detected troubled cells, which is commonly done \cite{Vui14R, Vui15R, Zhu09Q}. 

\subsection{One-dimensional tests}
\begin{figure}[ht!]
 \centering
\subfigure[Multiwavelets, $C=0.5$]{\includegraphics[scale = 0.27]{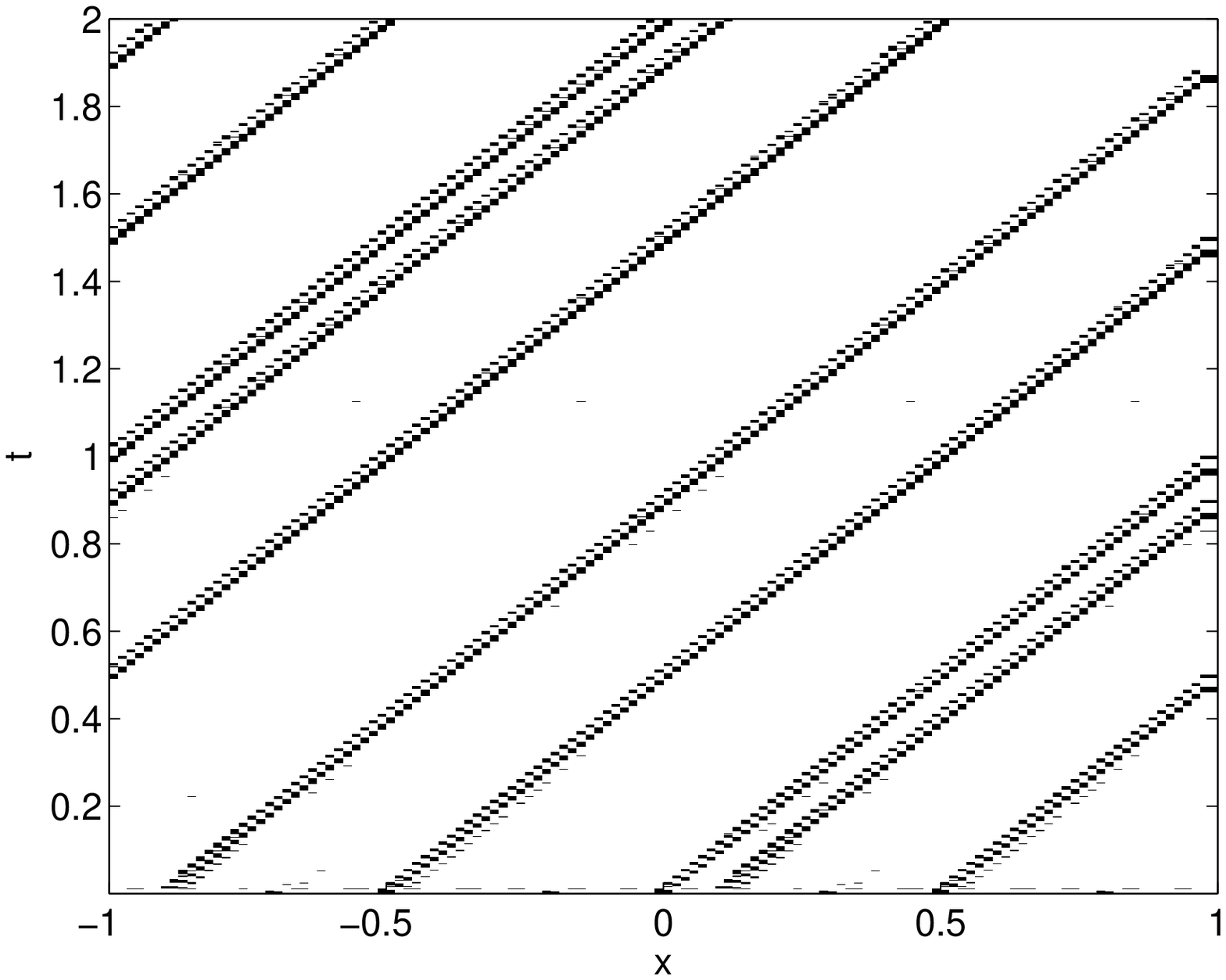}}
\subfigure[KXRCF, threshold 1]{\includegraphics[scale = 0.27]{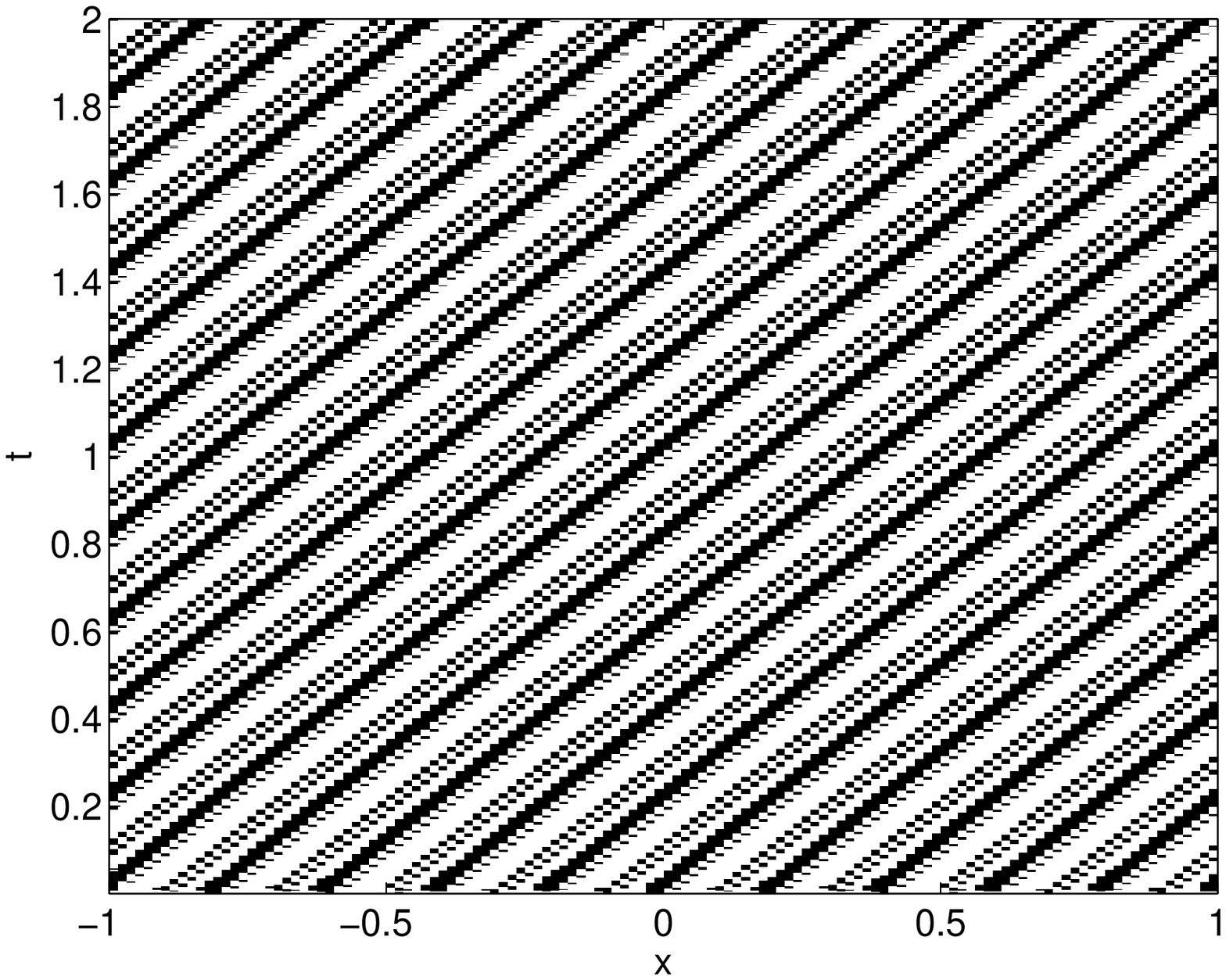}}
\subfigure[Minmod, $M=10$]{\includegraphics[scale = 0.27]{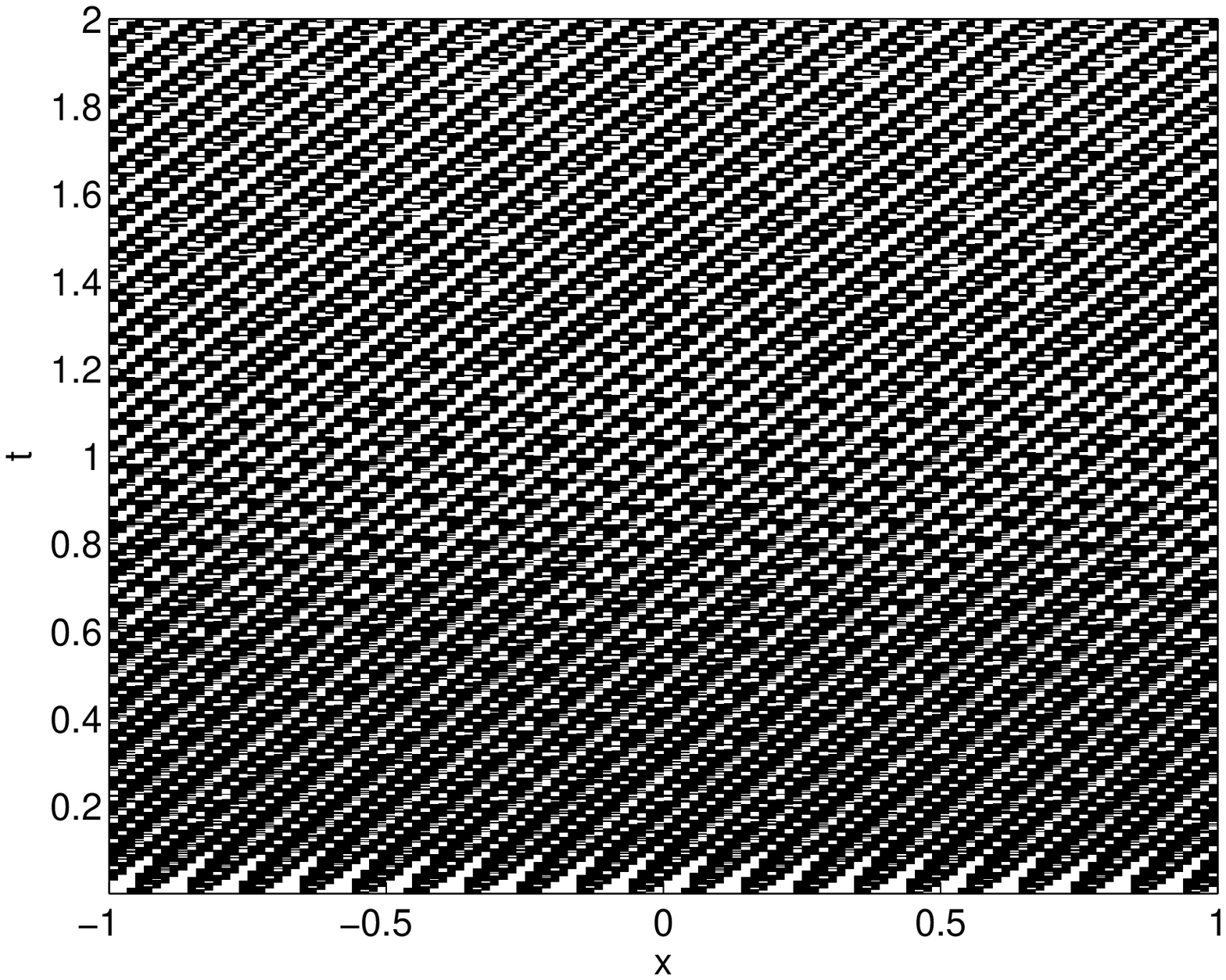}} \\
              \vspace{-0.4cm}
\caption{Detected troubled cells for $\rho(x,2) = 1 + 0.5\sin(10\pi x)$, 128 elements, $k=2$, using original troubled-cell indicators. Corresponding outlier-detection approaches do not detect any element.}\label{fig:Ccont}
\end{figure}

The test cases in one dimension include one continuous example using the Euler equations on $[-1,1]$ with initial conditions $\rho_0(x) = 1 + 0.5\sin(10\pi x)$, $u_0(x) = 1$, $p_0(x) = 1$, and periodic boundary conditions. The solution at final time $T=2$ is given by $\rho(x,2) = \rho_0(x)$. Using this example, we can validate our algorithm: since no discontinuities are present, no element should be detected. Indeed, the original troubled-cell indicators do detect certain elements (chosen parameters are reasonable, and commonly used \cite{Coc89S, Kri04XRCF, Vui14R}). This is depicted in Figure \ref{fig:Ccont}, in which the detected troubled cells using the original indicators are visualized. These so-called time-history plots show which elements are detected in space for each time step.
    
The application of the outlier-detection algorithm together with the troubled-cell indication vectors does not select any element, which is the desirable result.

The standard numerical examples for the Euler equations are also investigated.

\begin{figure}[h!]
 \centering
\subfigure[Original, $C=0.1$]{\includegraphics[scale = 0.27]{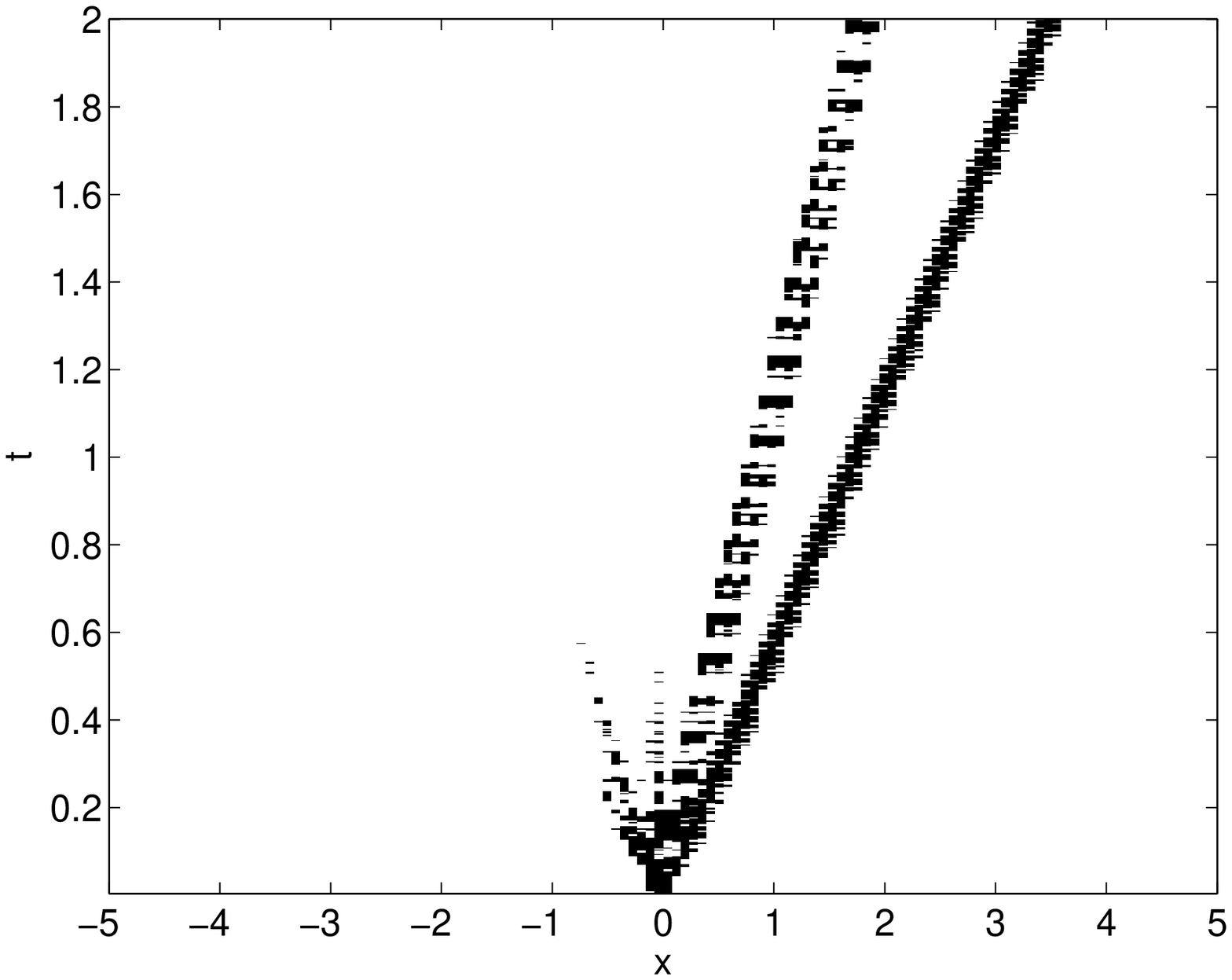}}
\subfigure[Original, KXRCF]{\includegraphics[scale = 0.27]{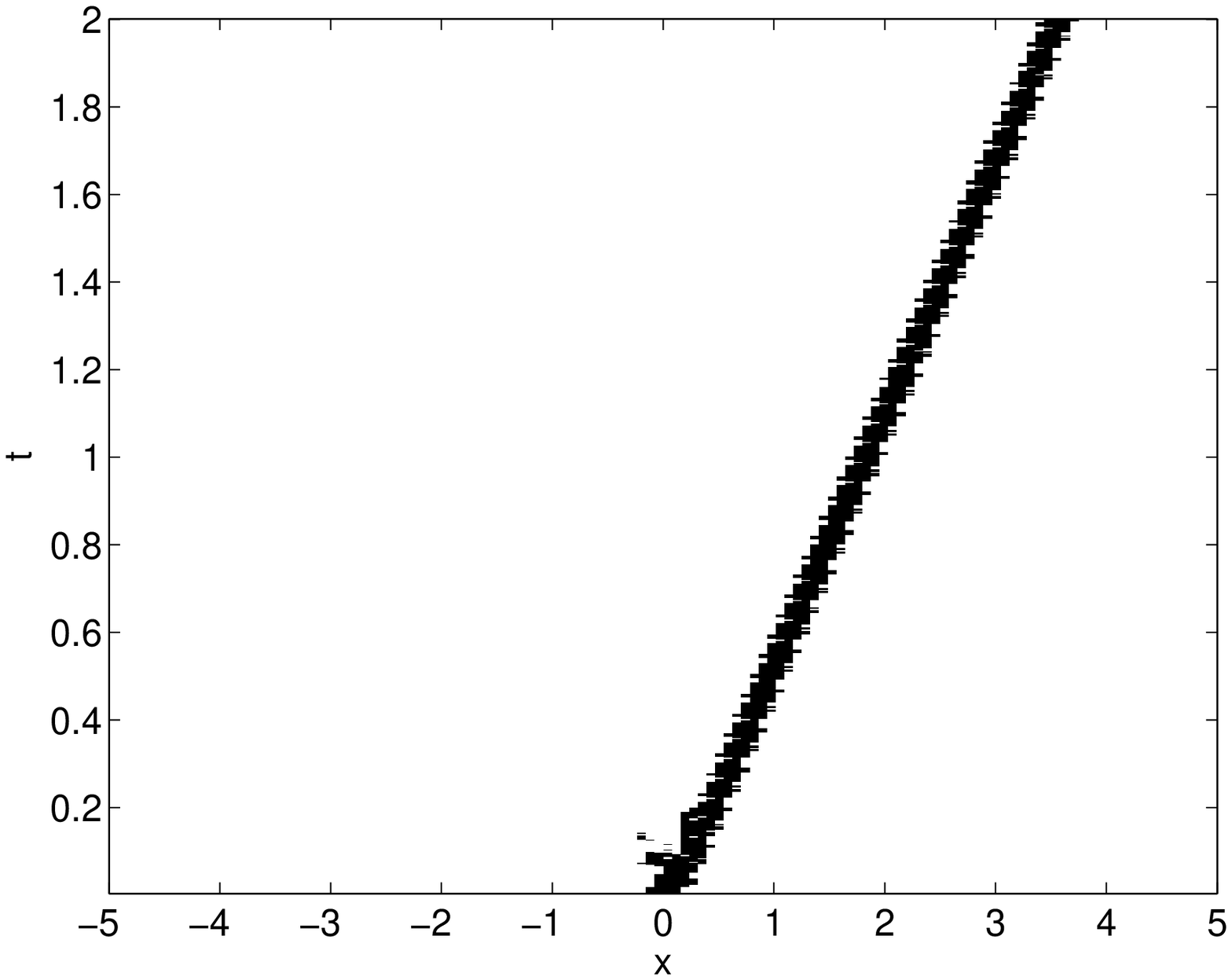}}
\subfigure[Original, $M = 10$]{\includegraphics[scale = 0.27]{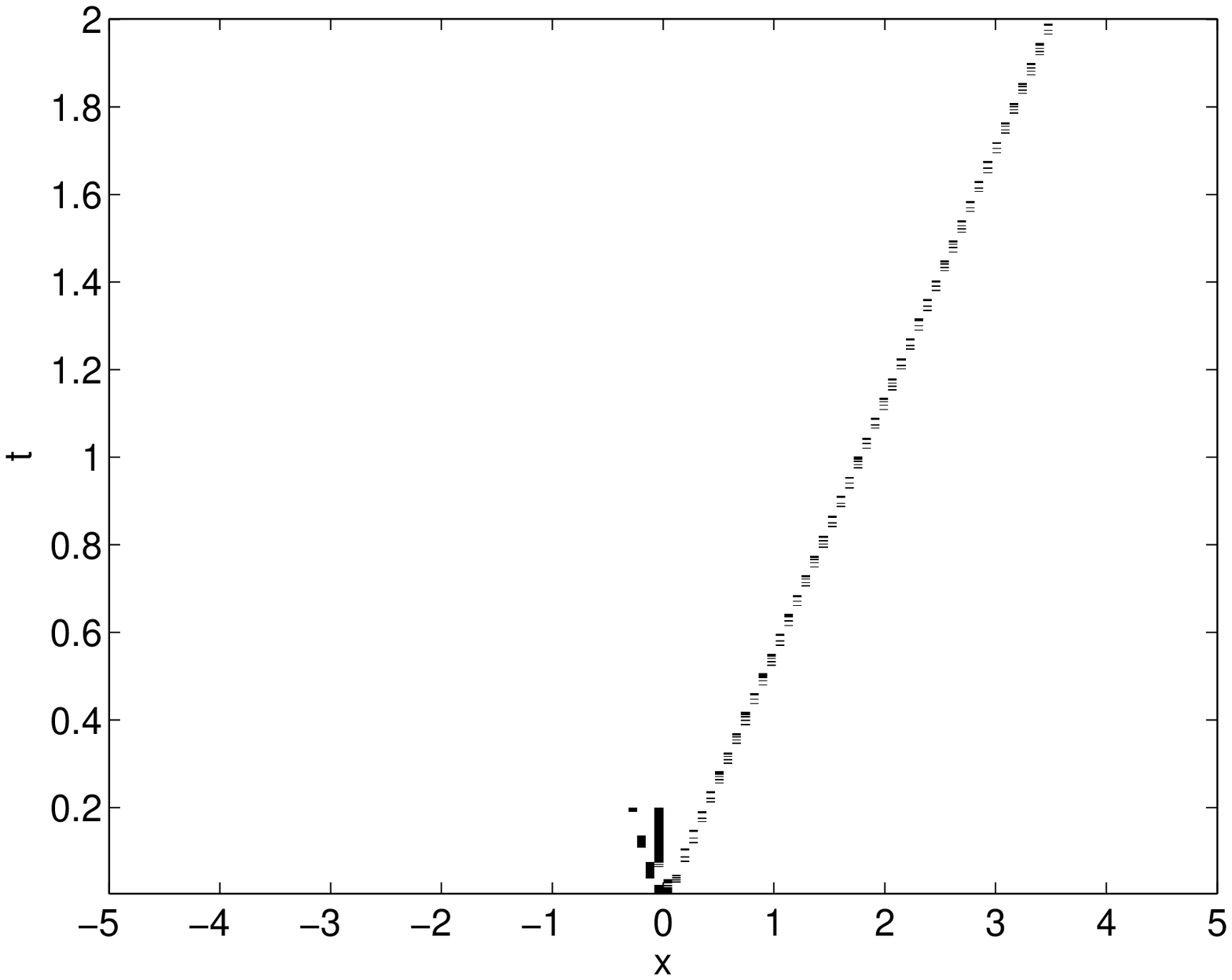}} \\
\vspace{-0.4cm}
\subfigure[Outlier, multiwavelets]{\includegraphics[scale = 0.27]{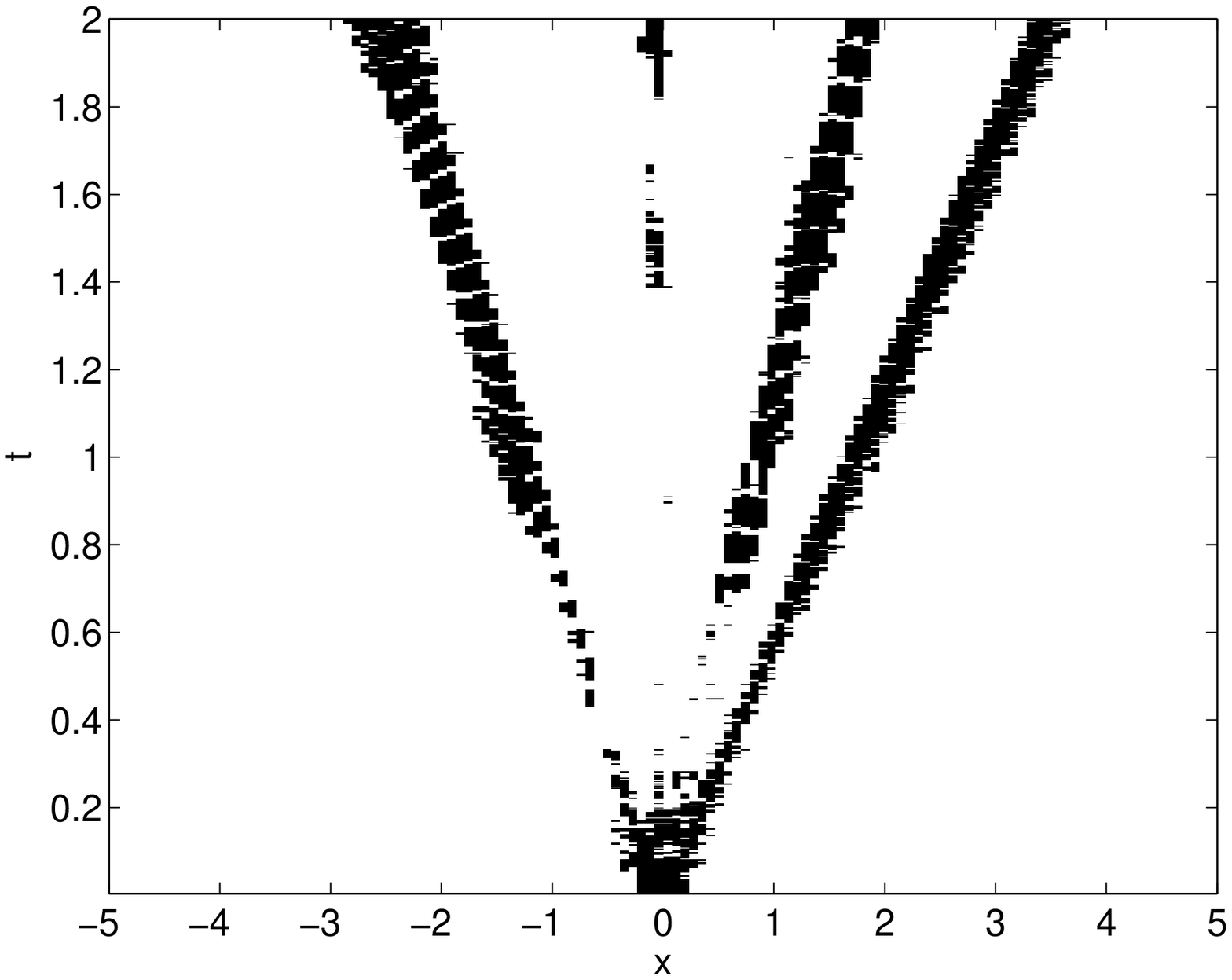}}
\subfigure[Outlier, KXRCF value]{\includegraphics[scale = 0.27]{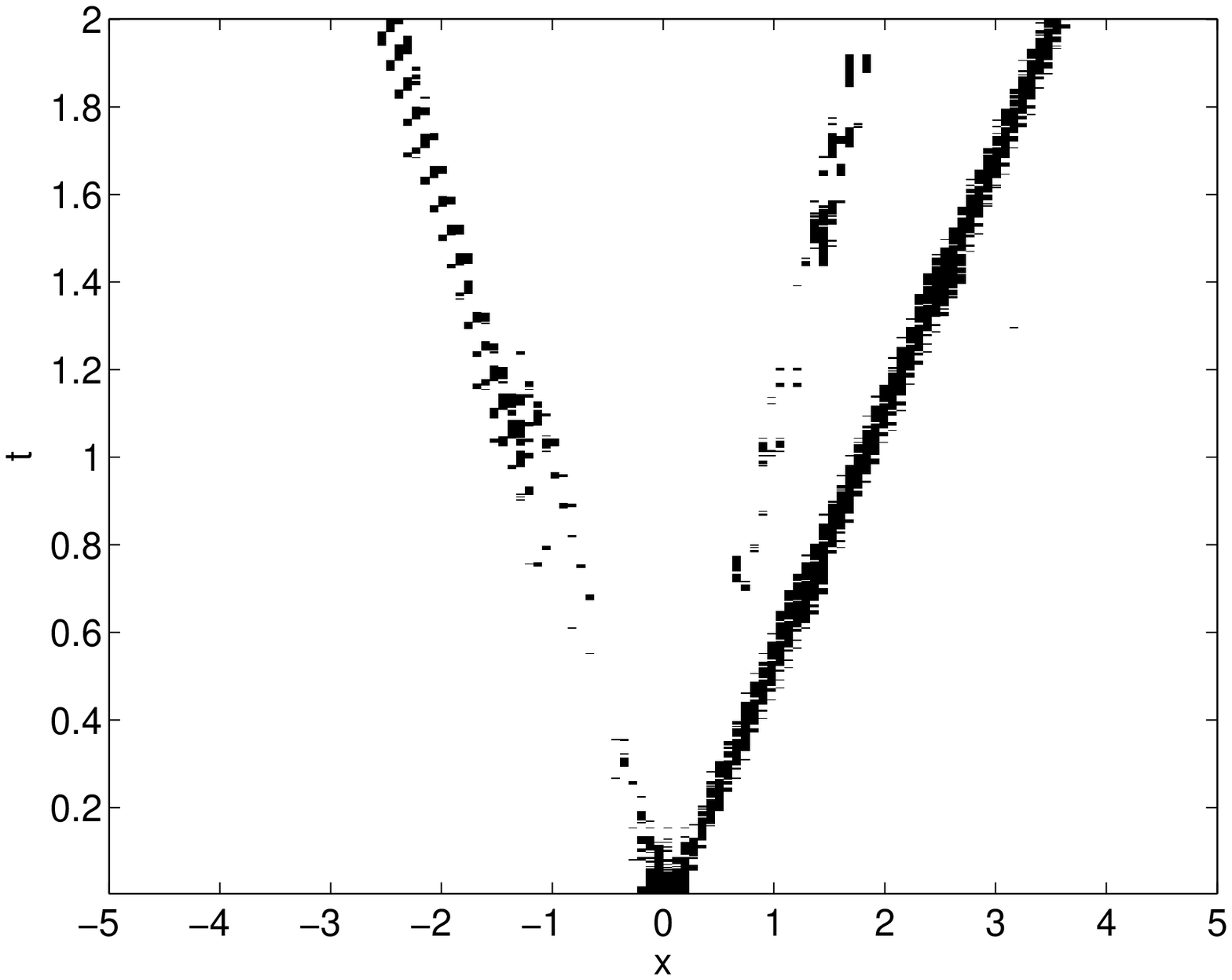}}
\subfigure[Outlier, minmod-based TVB]{\includegraphics[scale = 0.27]{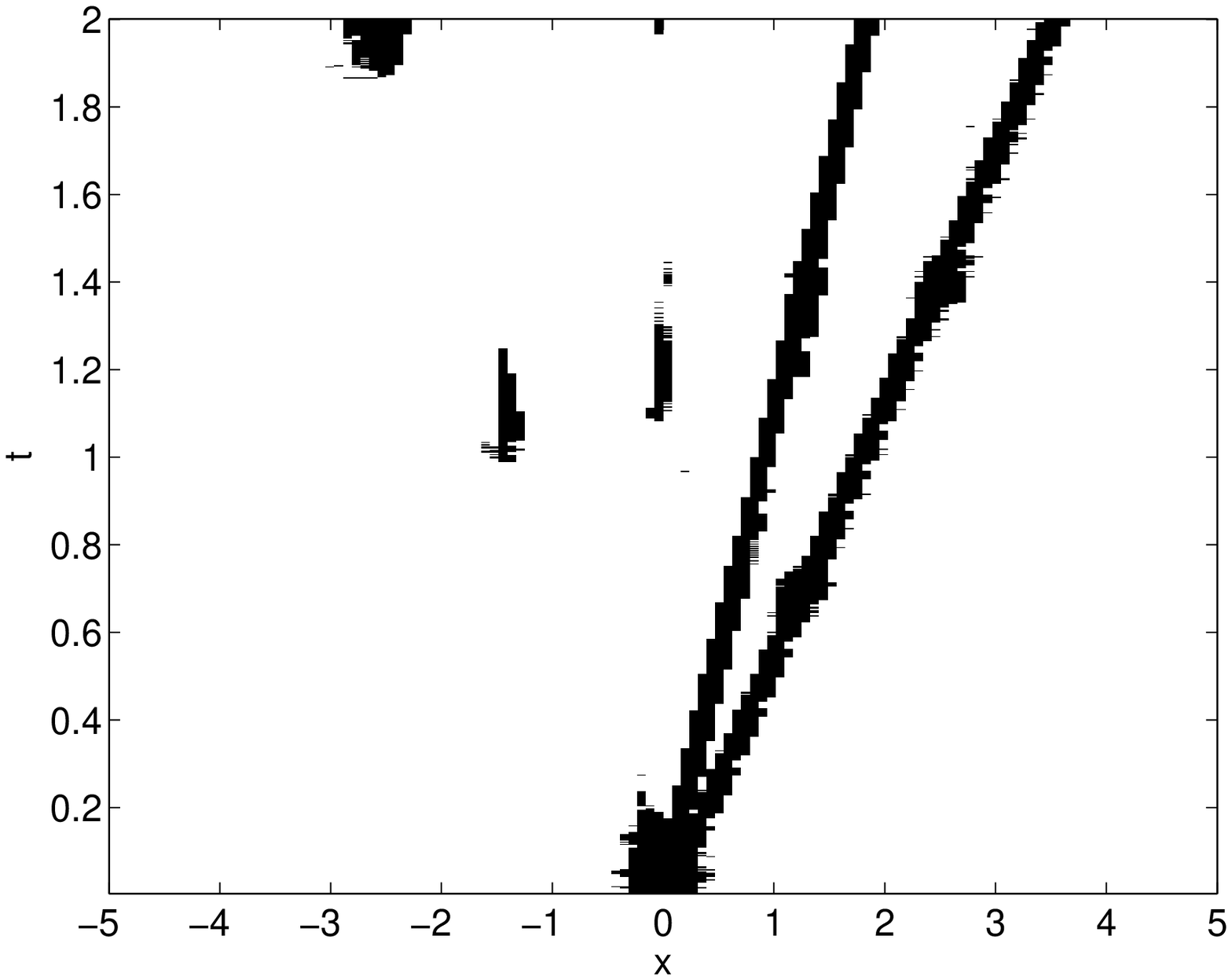}} \\
\vspace{-0.4cm}
\subfigure[Original, $C=0.1$]{\includegraphics[scale = 0.27]{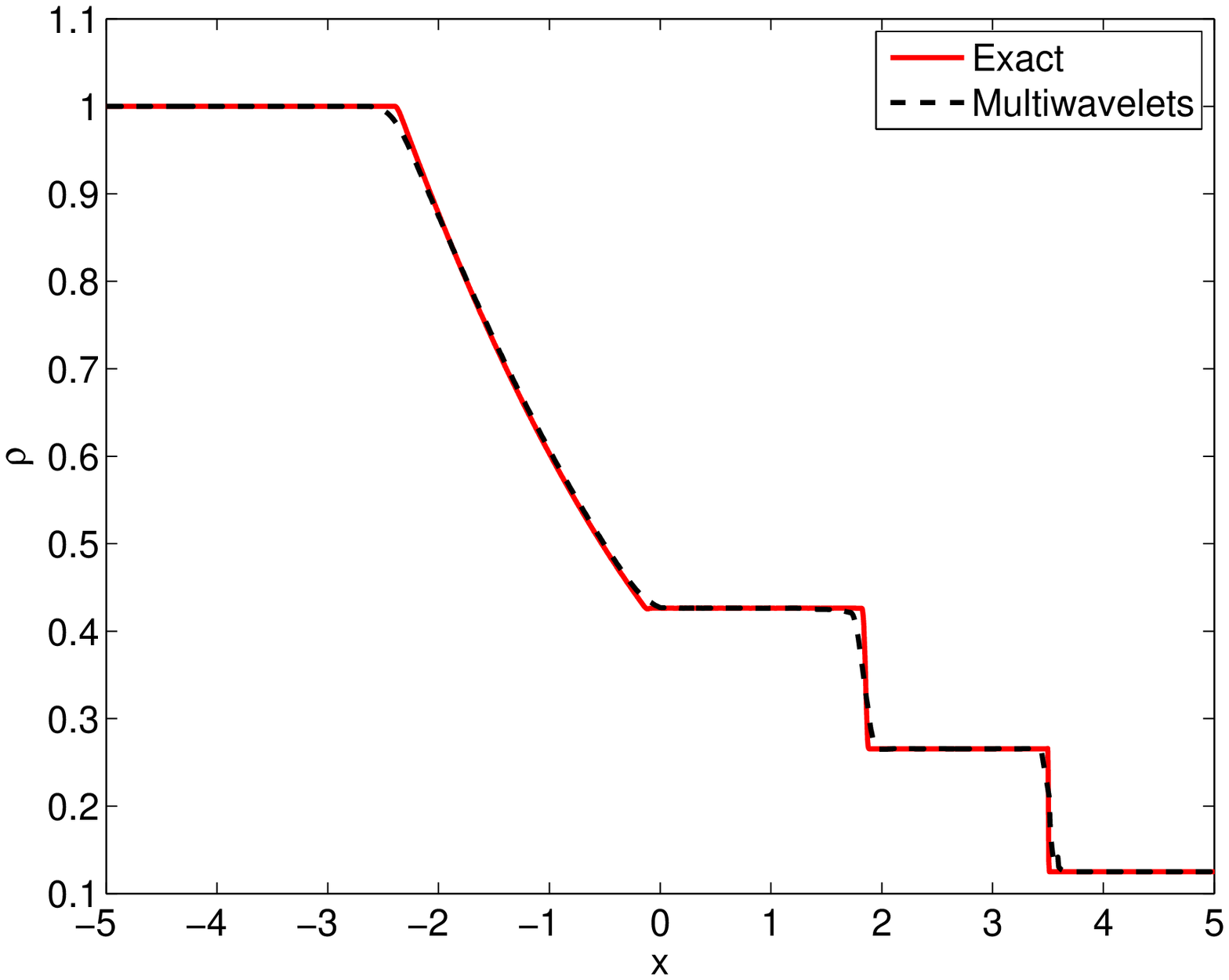}}
\subfigure[Original, KXRCF]{\includegraphics[scale = 0.27]{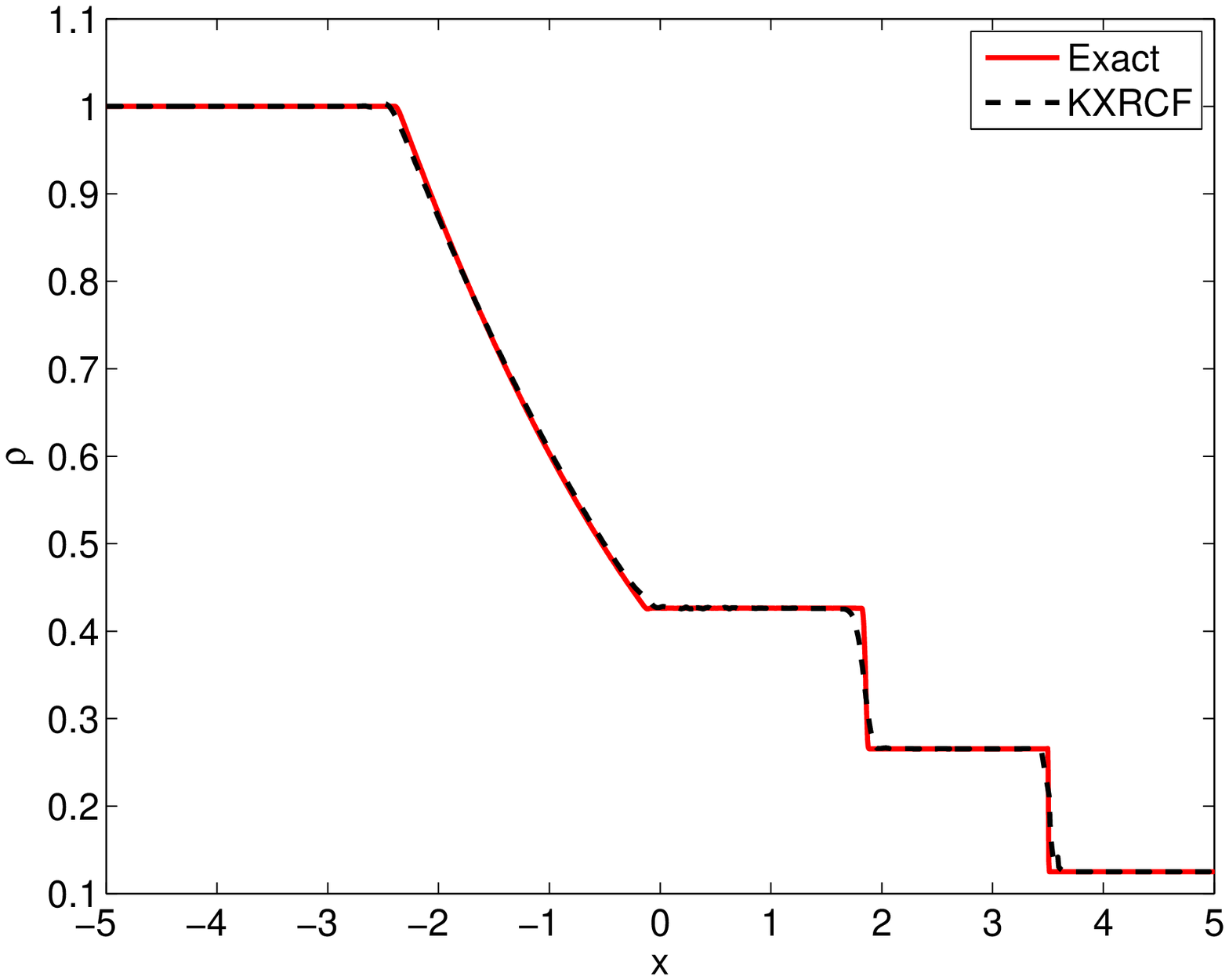}}
\subfigure[Original, $M = 10$]{\includegraphics[scale = 0.27]{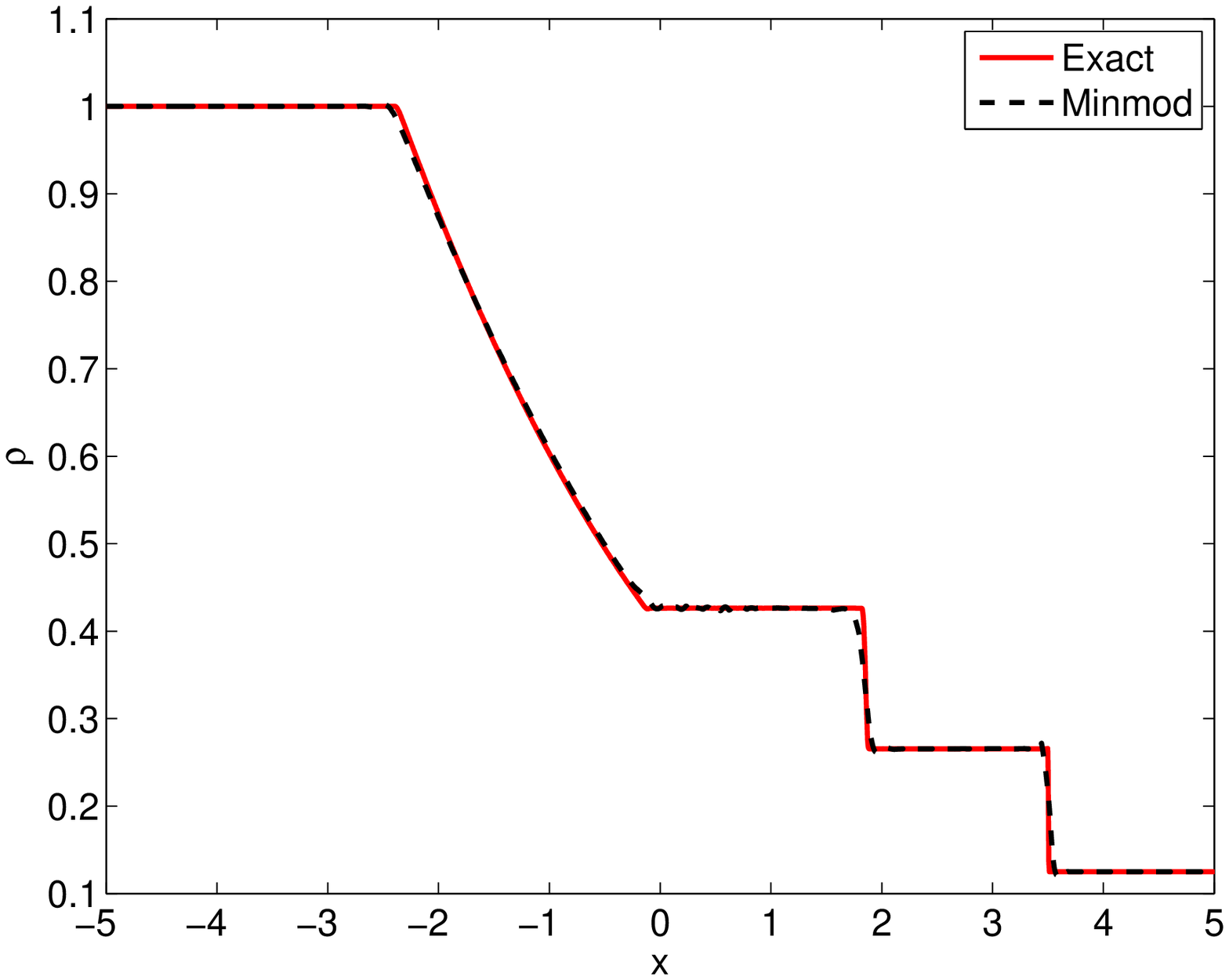}} \\
\vspace{-0.4cm}
\subfigure[Outlier, multiwavelets]{\includegraphics[scale = 0.27]{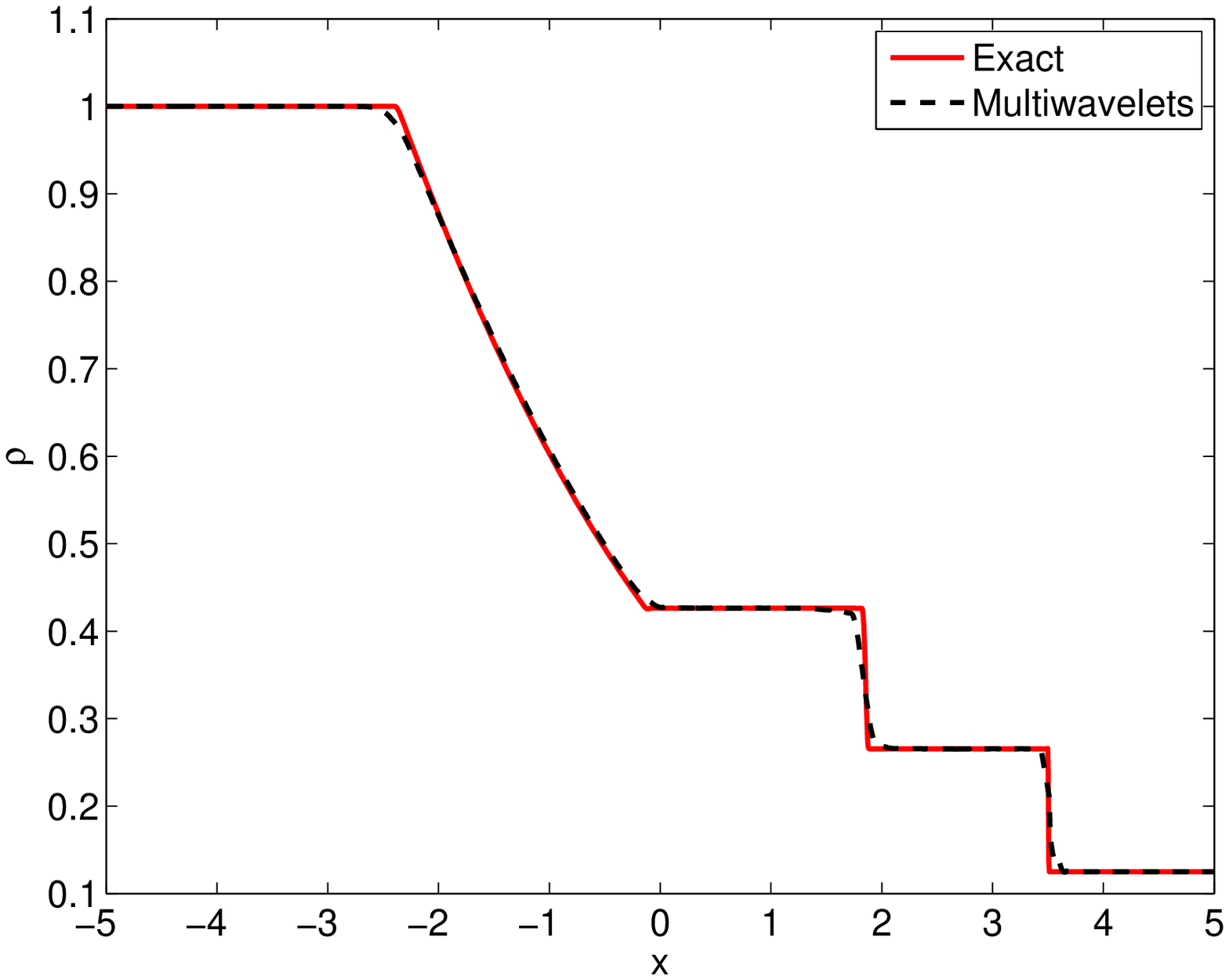}}
\subfigure[Outlier, KXRCF value]{\includegraphics[scale = 0.27]{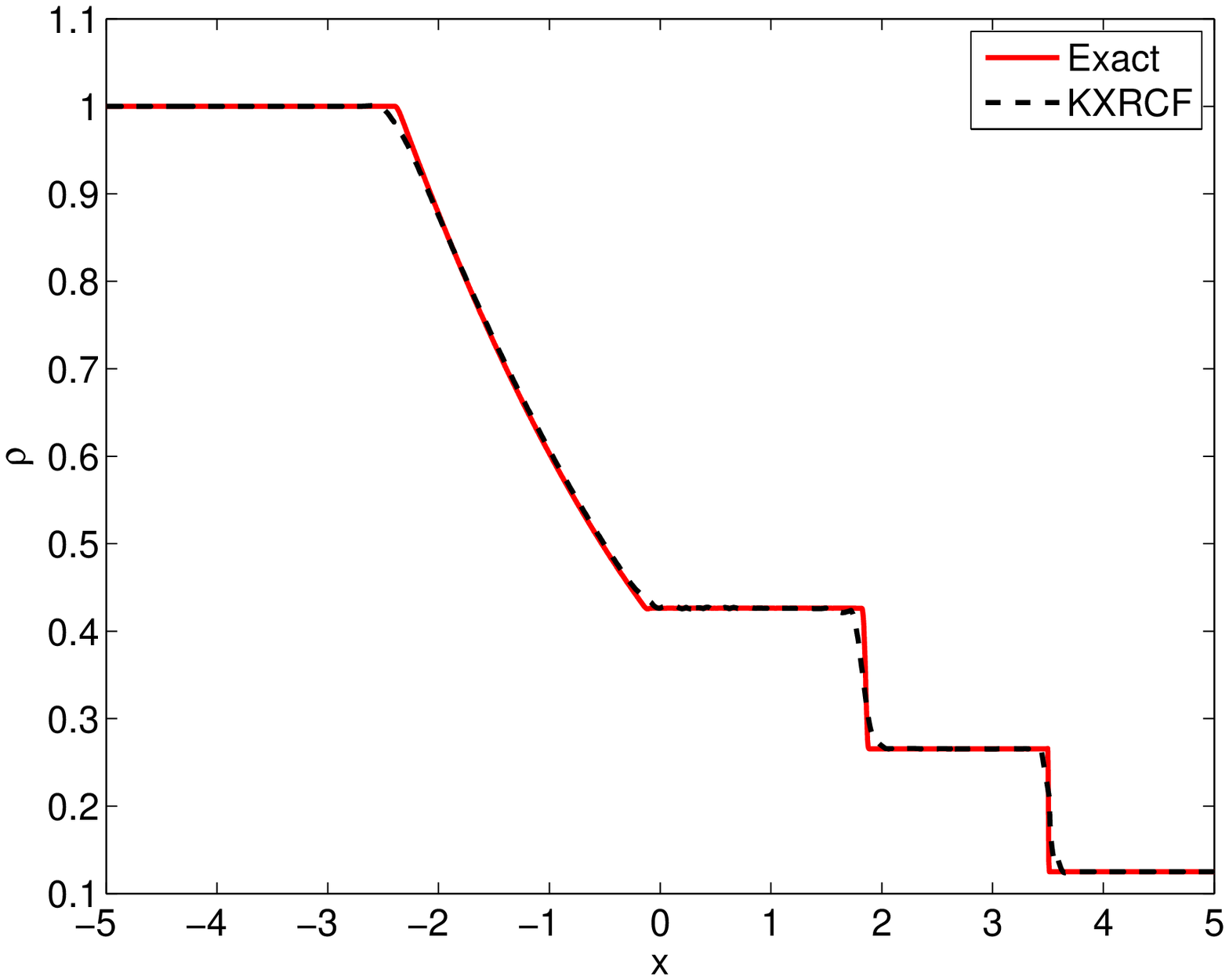}}
\subfigure[Outlier, minmod-based TVB]{\includegraphics[scale = 0.27]{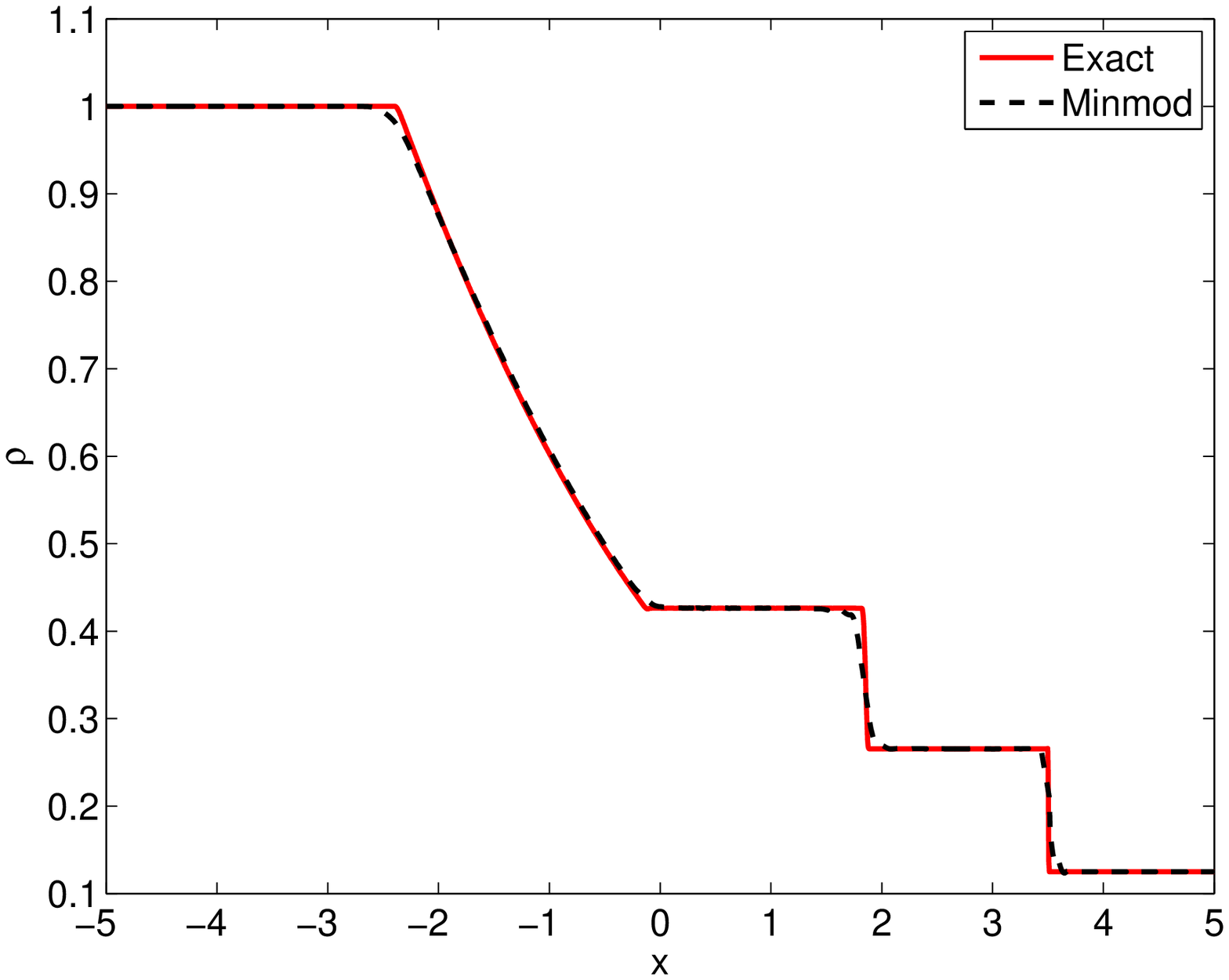}} \\
\vspace{-0.4cm}
\caption{Detected troubled cells (row 1 and 2) and approximation at final time $T=2$ (row 3 and 4), shock tube of Sod, $k=2$, 128 elements.}\label{fig:Sodk2}
\end{figure}

\begin{figure}[h!]
 \centering
\subfigure[Original, $C=0.1$]{\includegraphics[scale = 0.27]{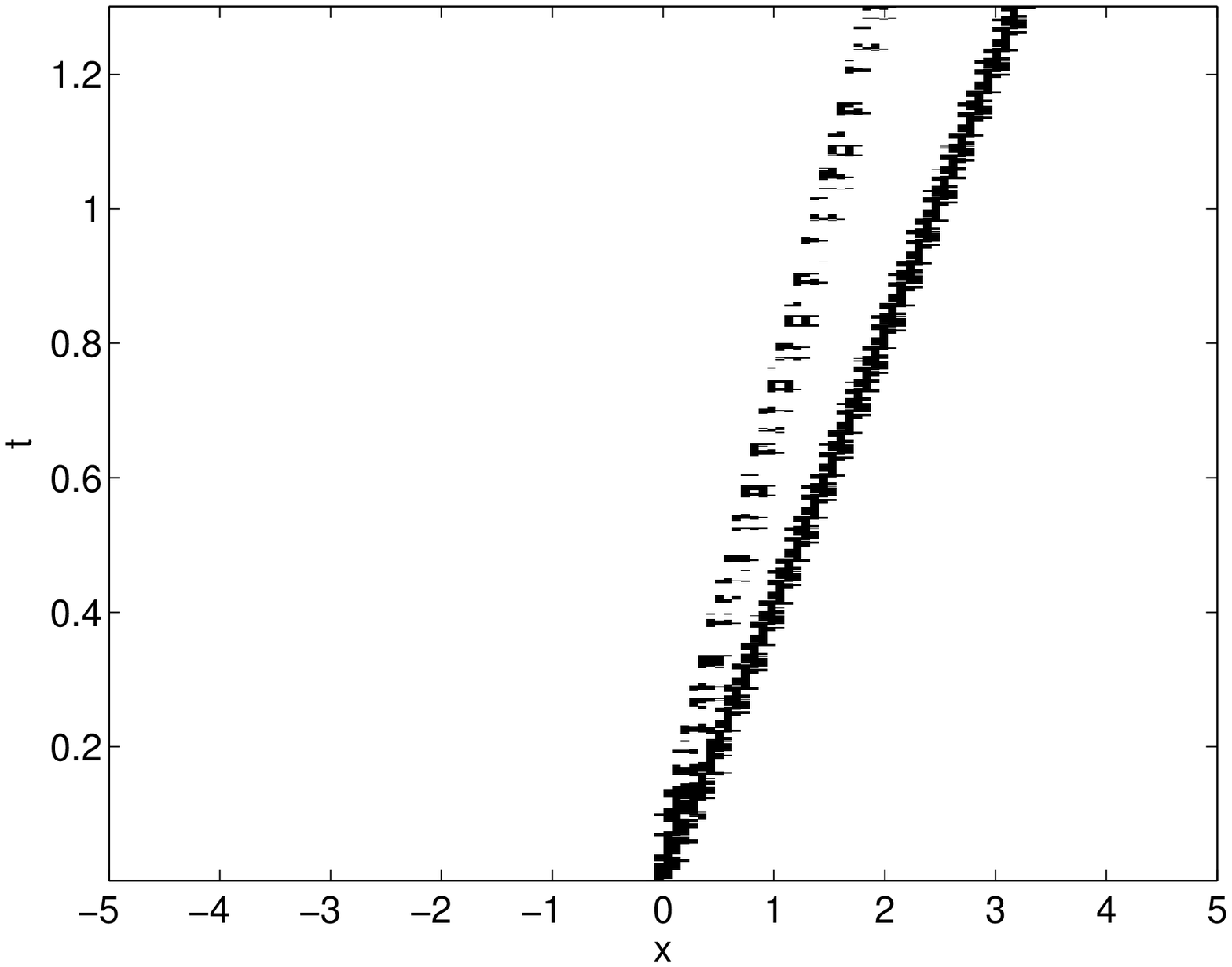}}
\subfigure[Original, KXRCF]{\includegraphics[scale = 0.27]{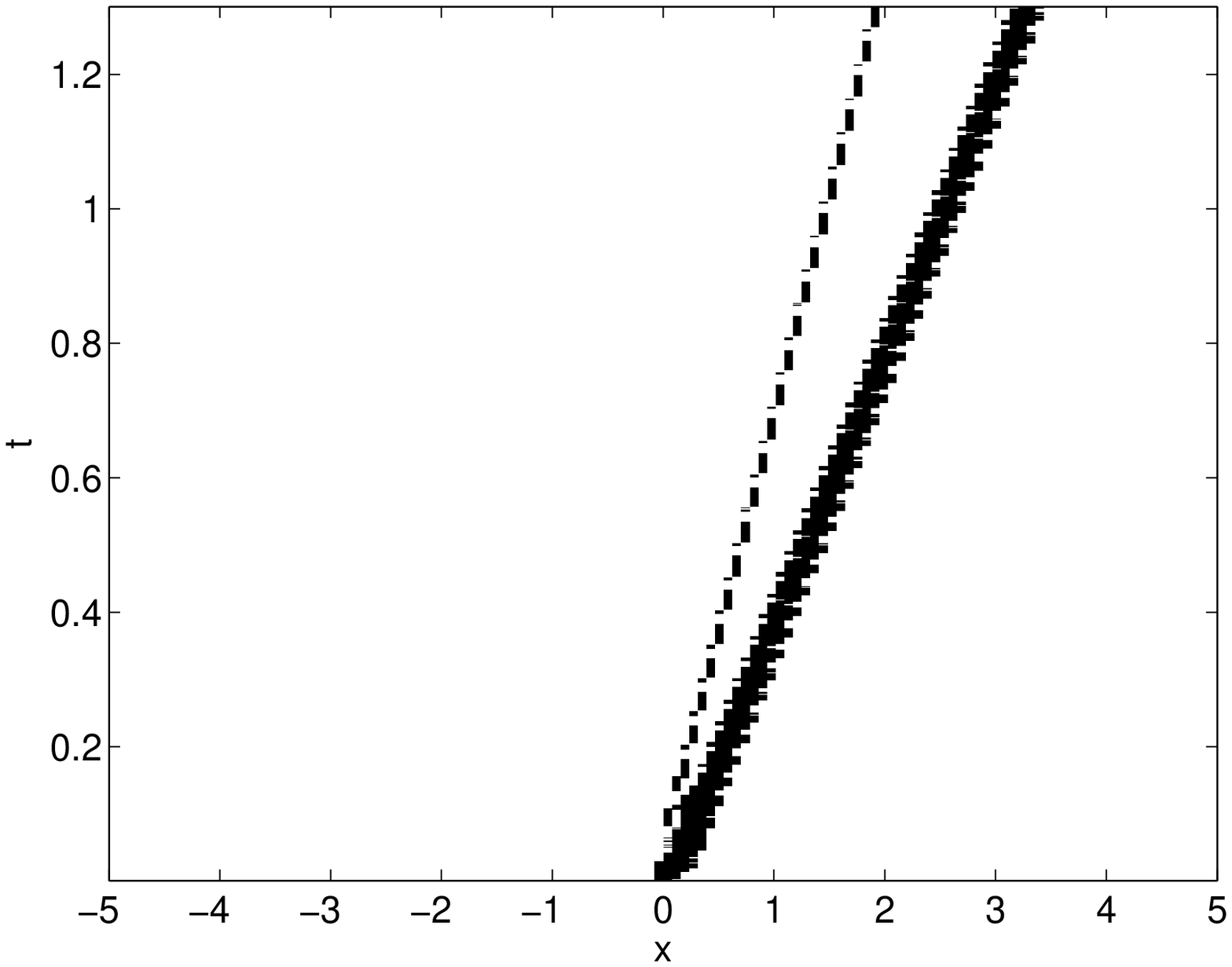}}
\subfigure[Original, $M = 10$]{\includegraphics[scale = 0.27]{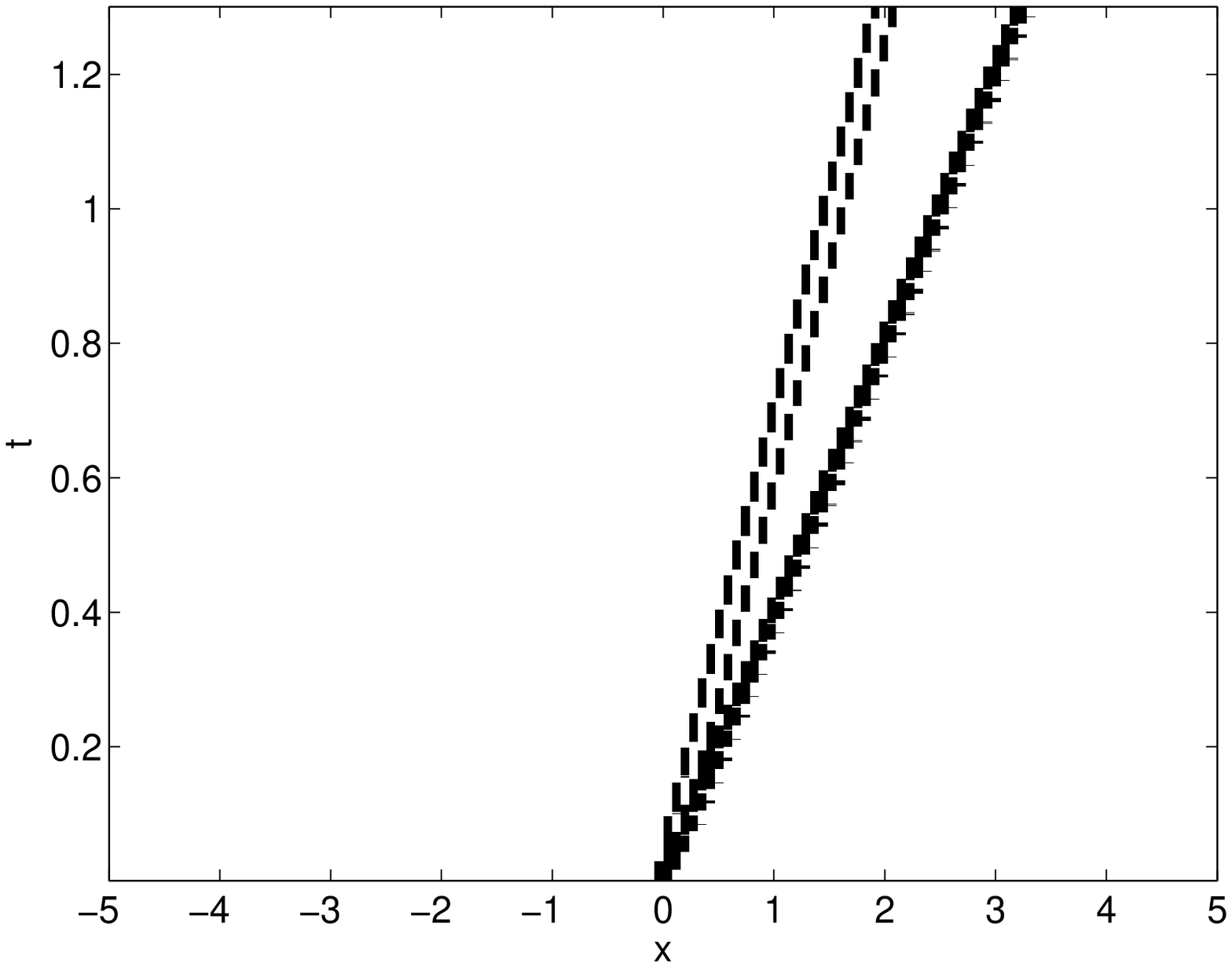}} \\
\vspace{-0.3cm}
\subfigure[Outlier, multiwavelets]{\includegraphics[scale = 0.27]{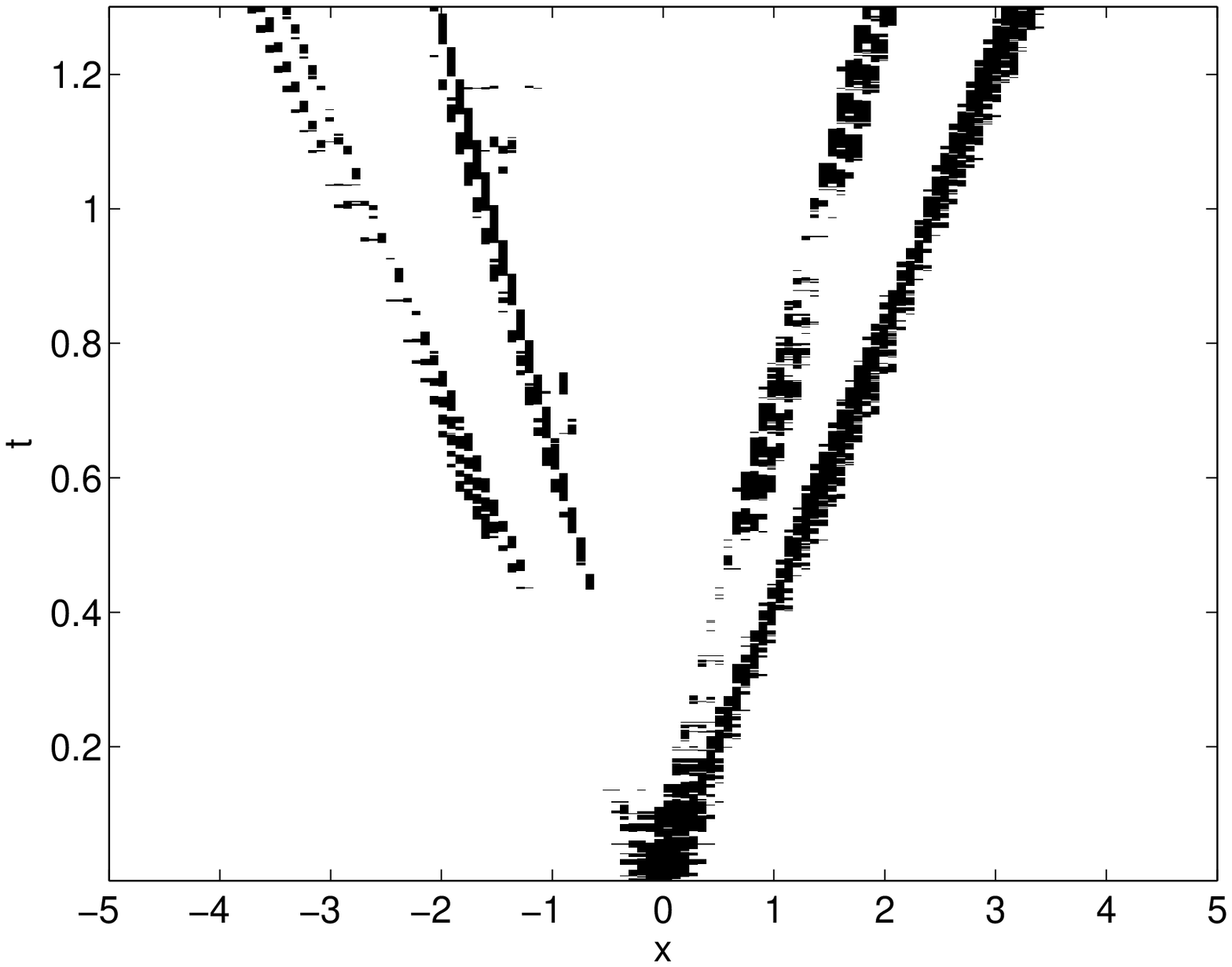}}
\subfigure[Outlier, KXRCF value]{\includegraphics[scale = 0.27]{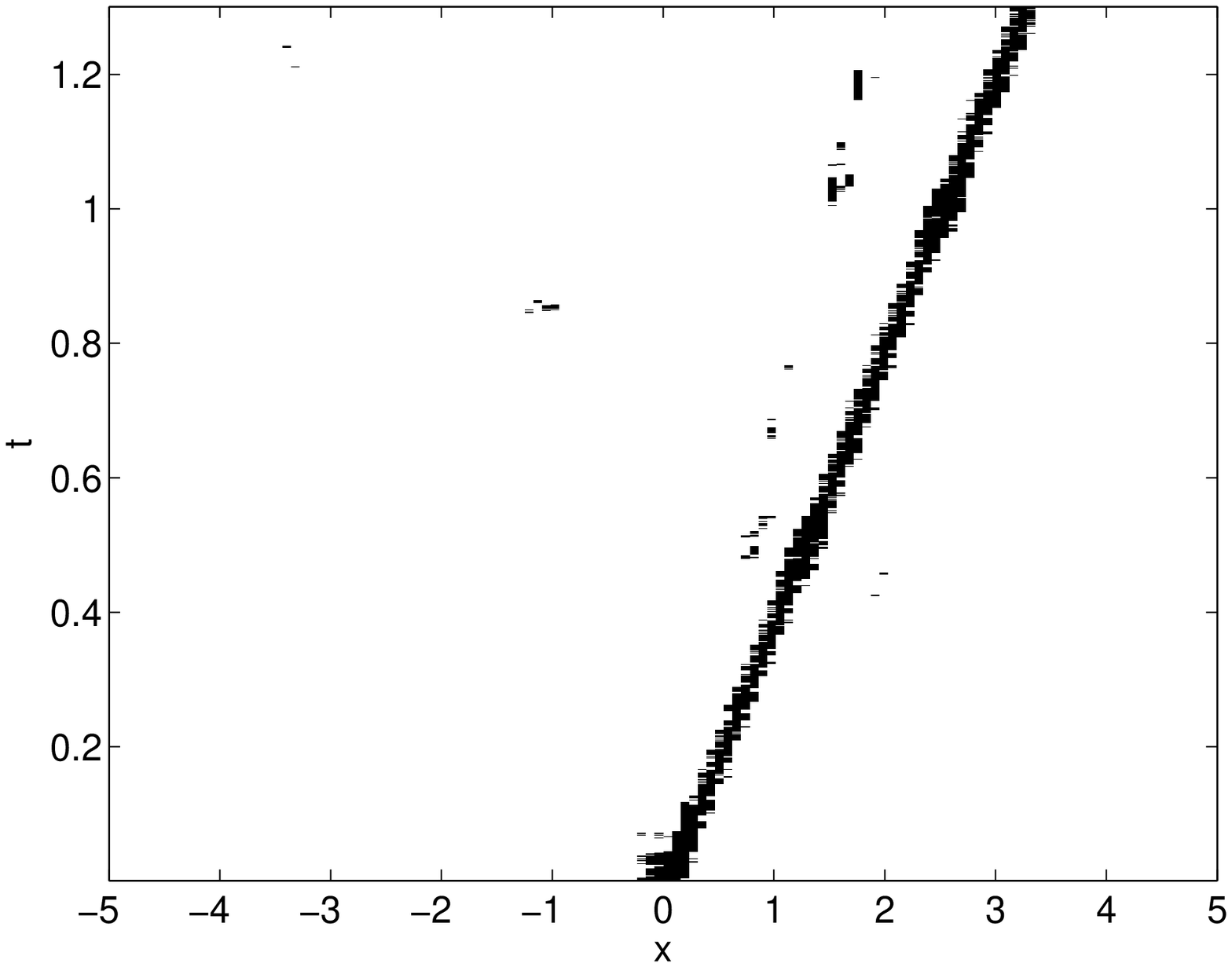}}
\subfigure[Outlier, minmod-based TVB]{\includegraphics[scale = 0.27]{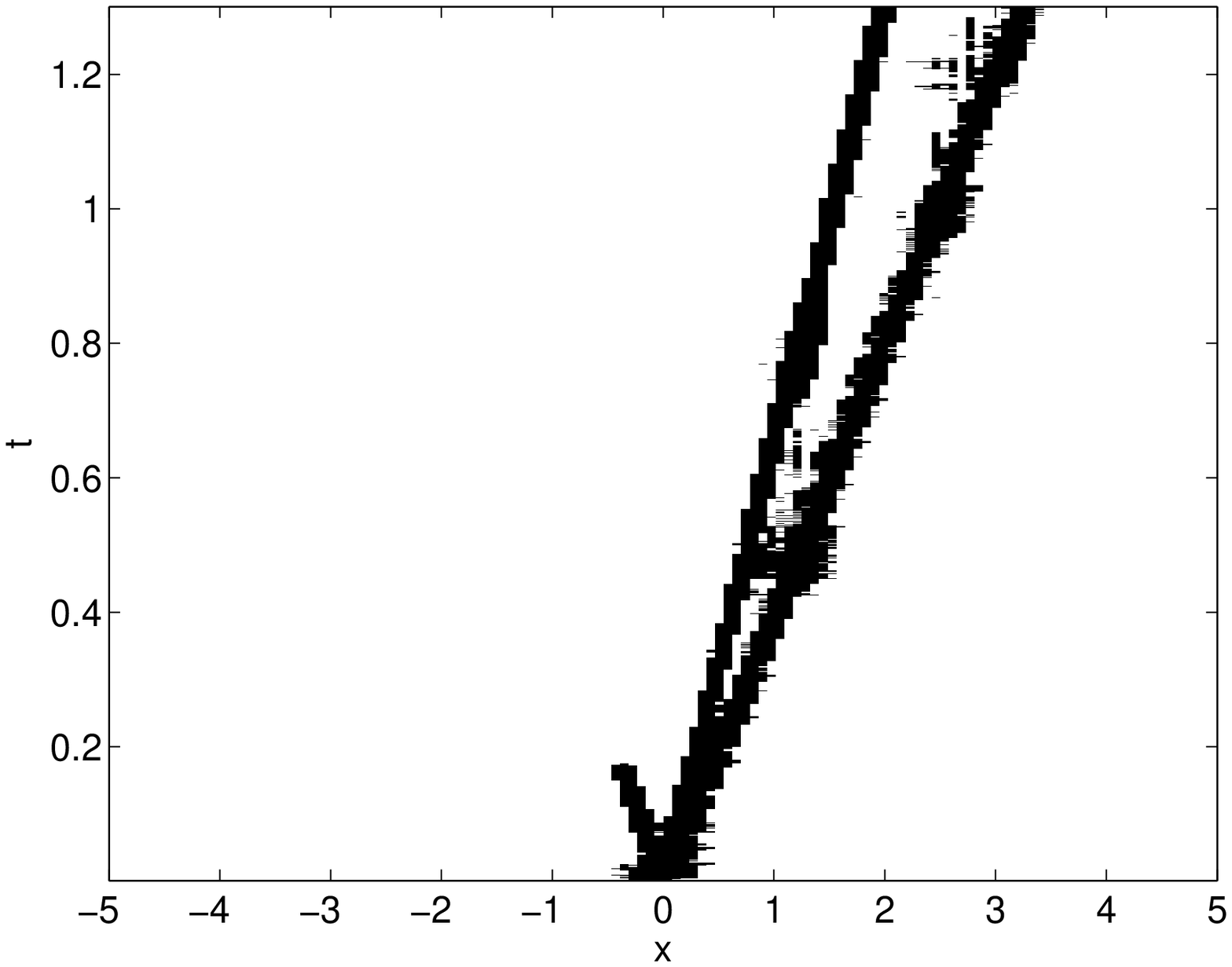}} \\
\vspace{-0.3cm}
\subfigure[Original, $C=0.1$]{\includegraphics[scale = 0.27]{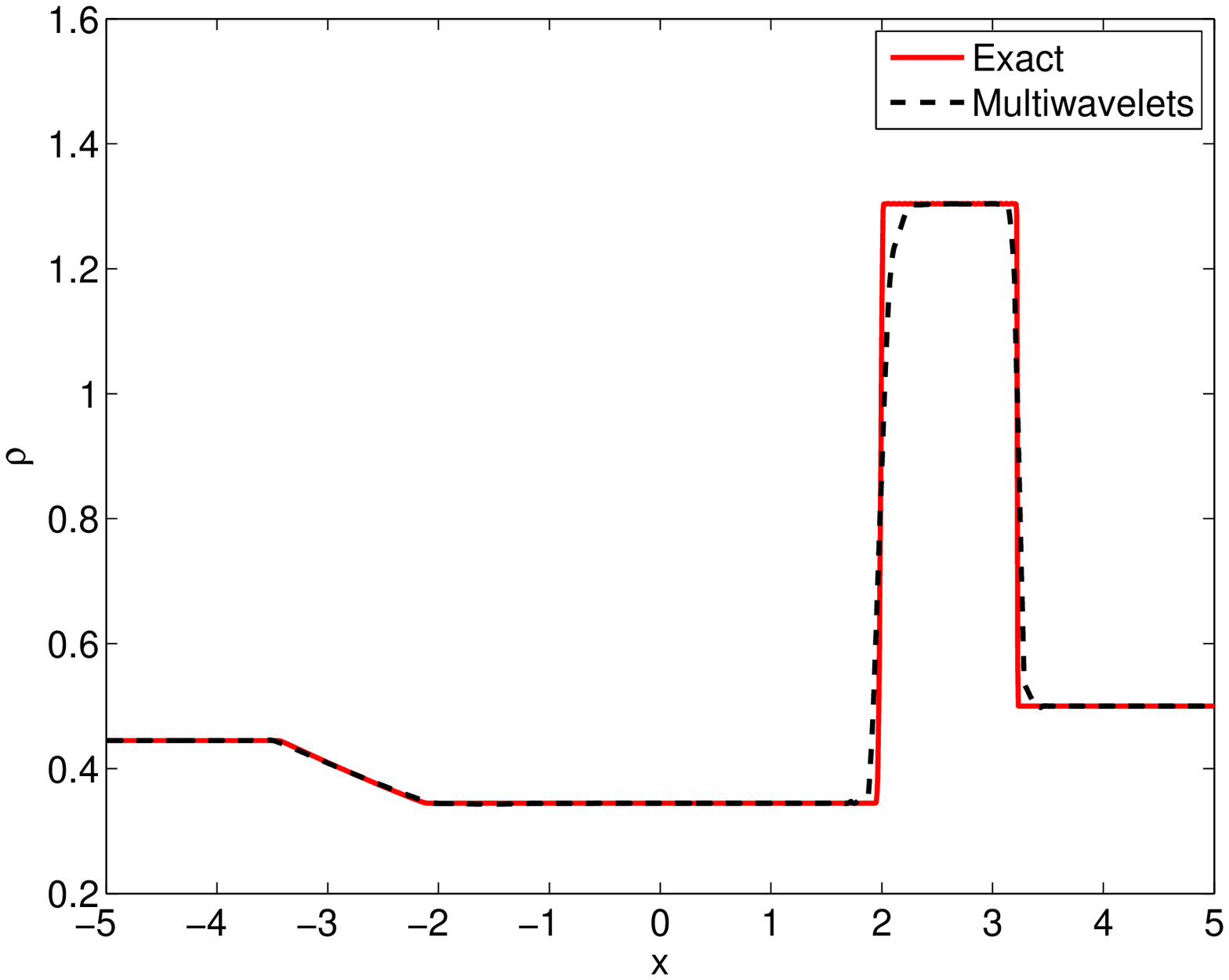}}
\subfigure[Original, KXRCF]{\includegraphics[scale = 0.27]{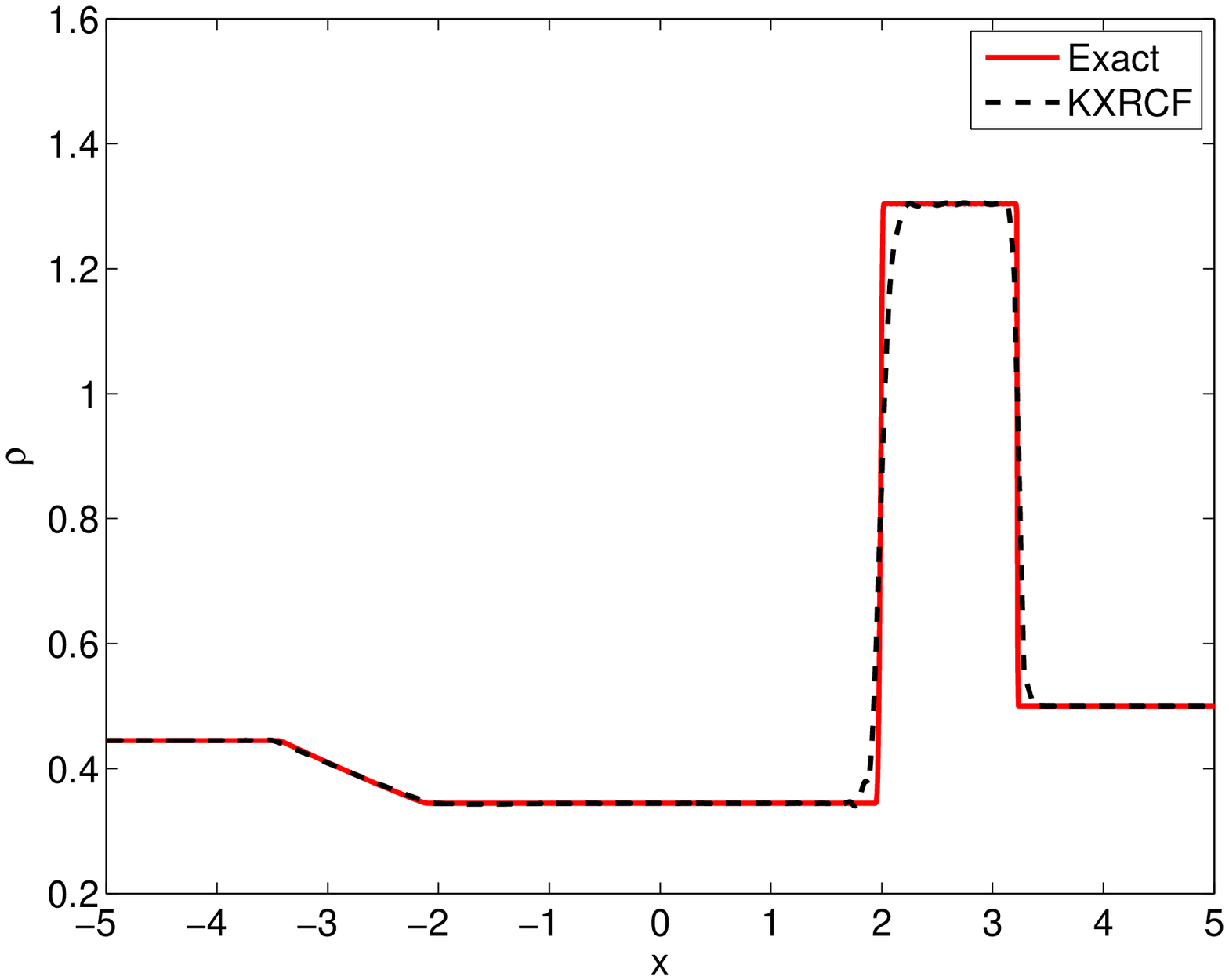}}
\subfigure[Original, $M = 10$]{\includegraphics[scale = 0.27]{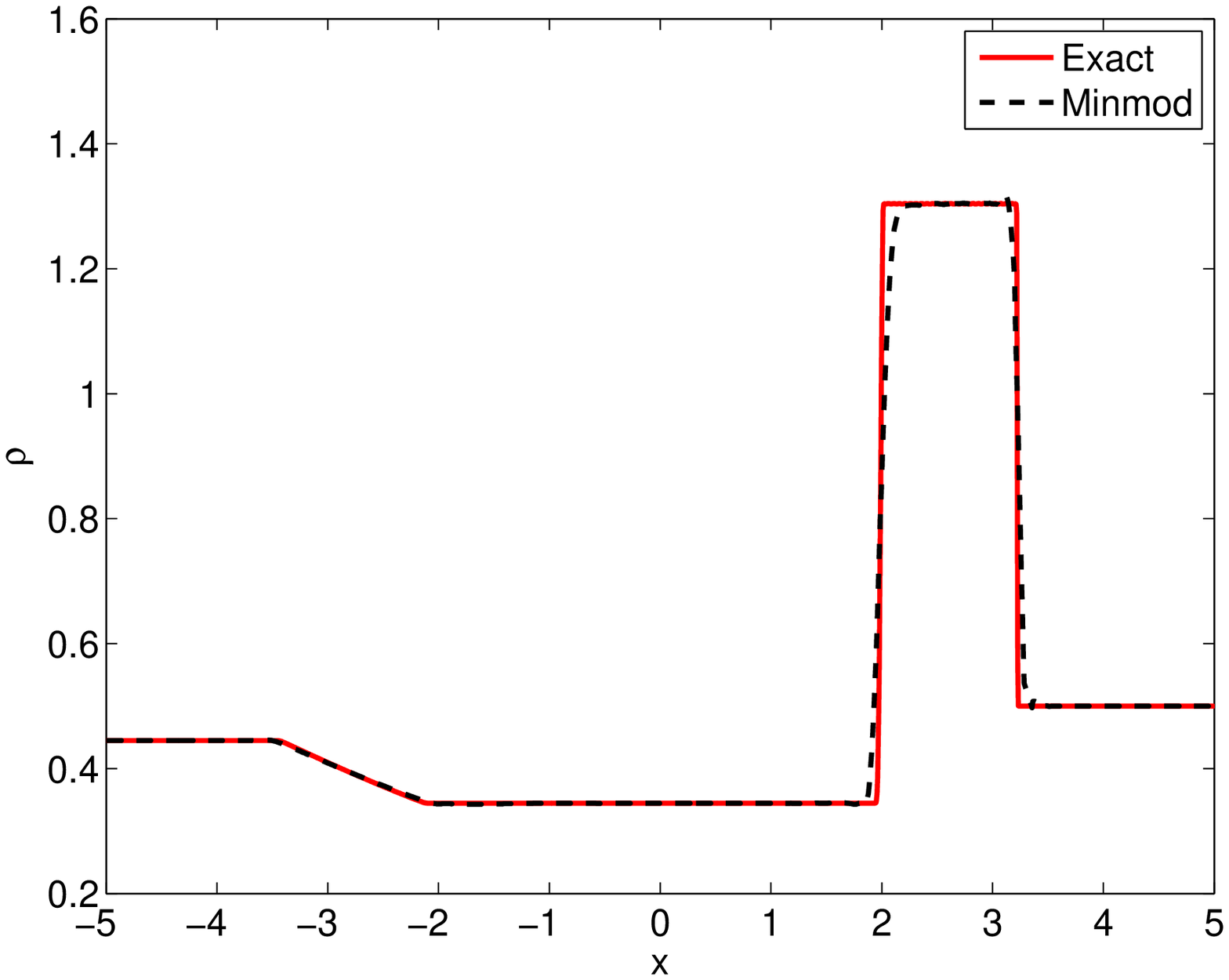}} \\
\vspace{-0.3cm}
\subfigure[Outlier, multiwavelets]{\includegraphics[scale = 0.27]{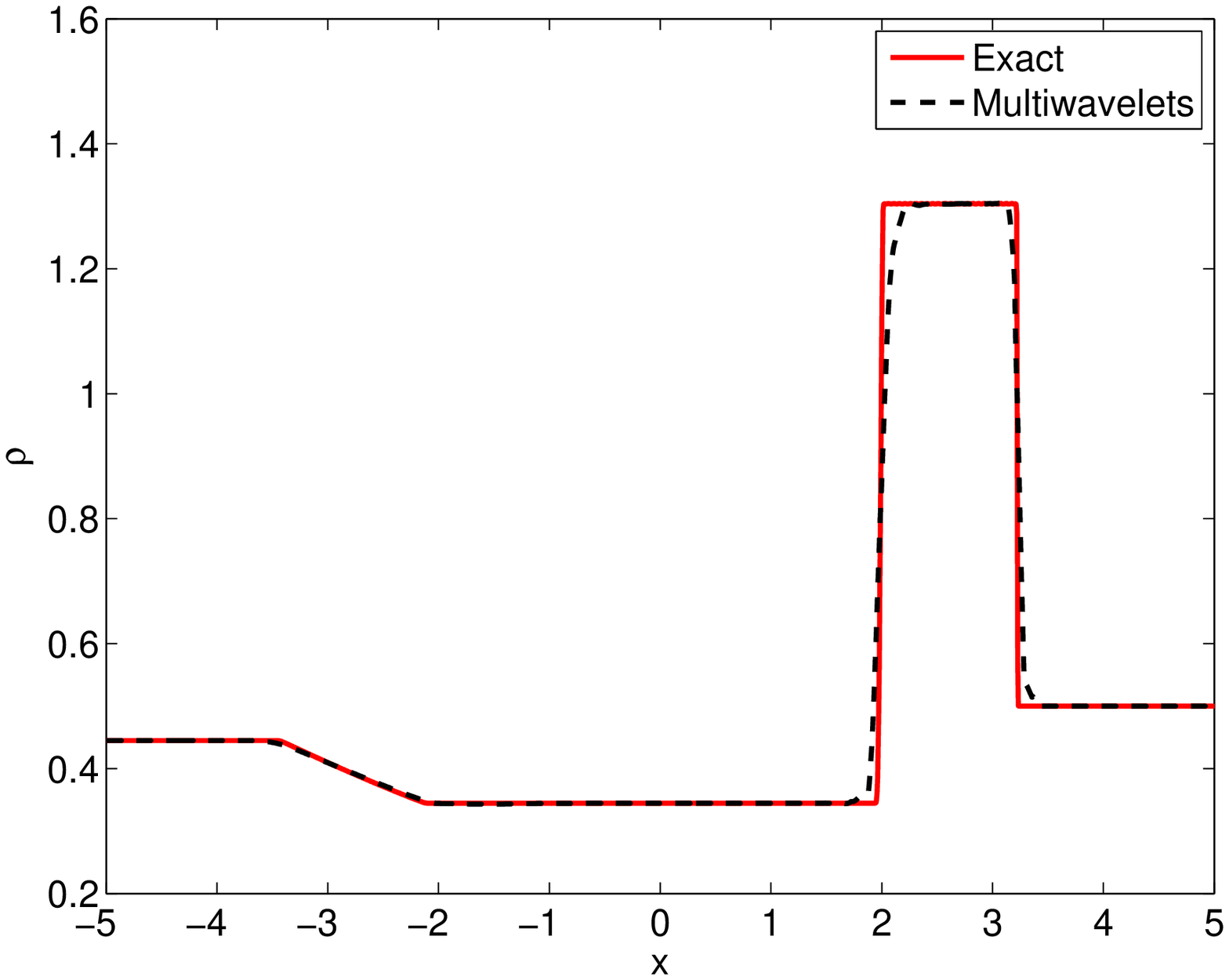}}
\subfigure[Outlier, KXRCF value]{\includegraphics[scale = 0.27]{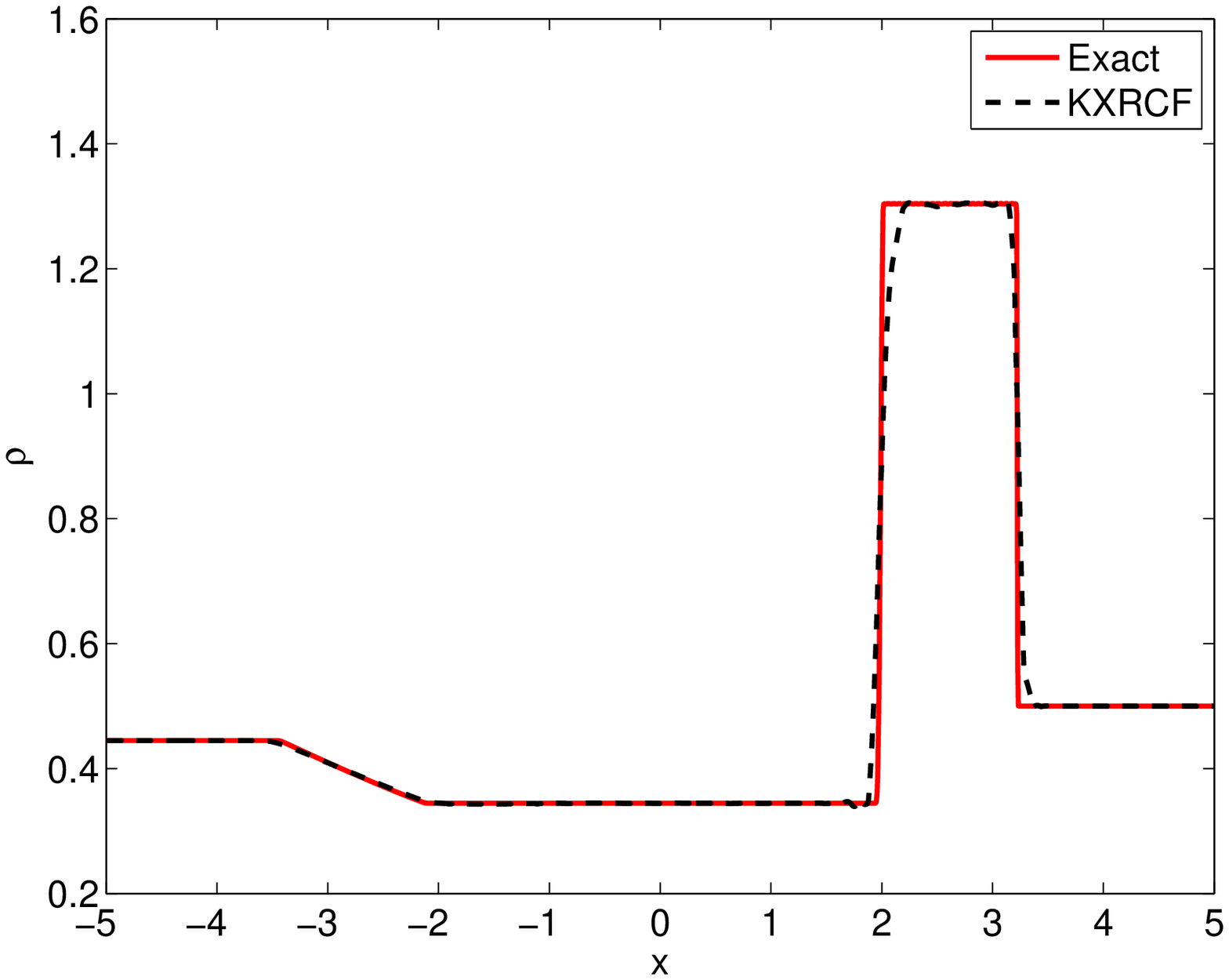}}
\subfigure[Outlier, minmod-based TVB]{\includegraphics[scale = 0.27]{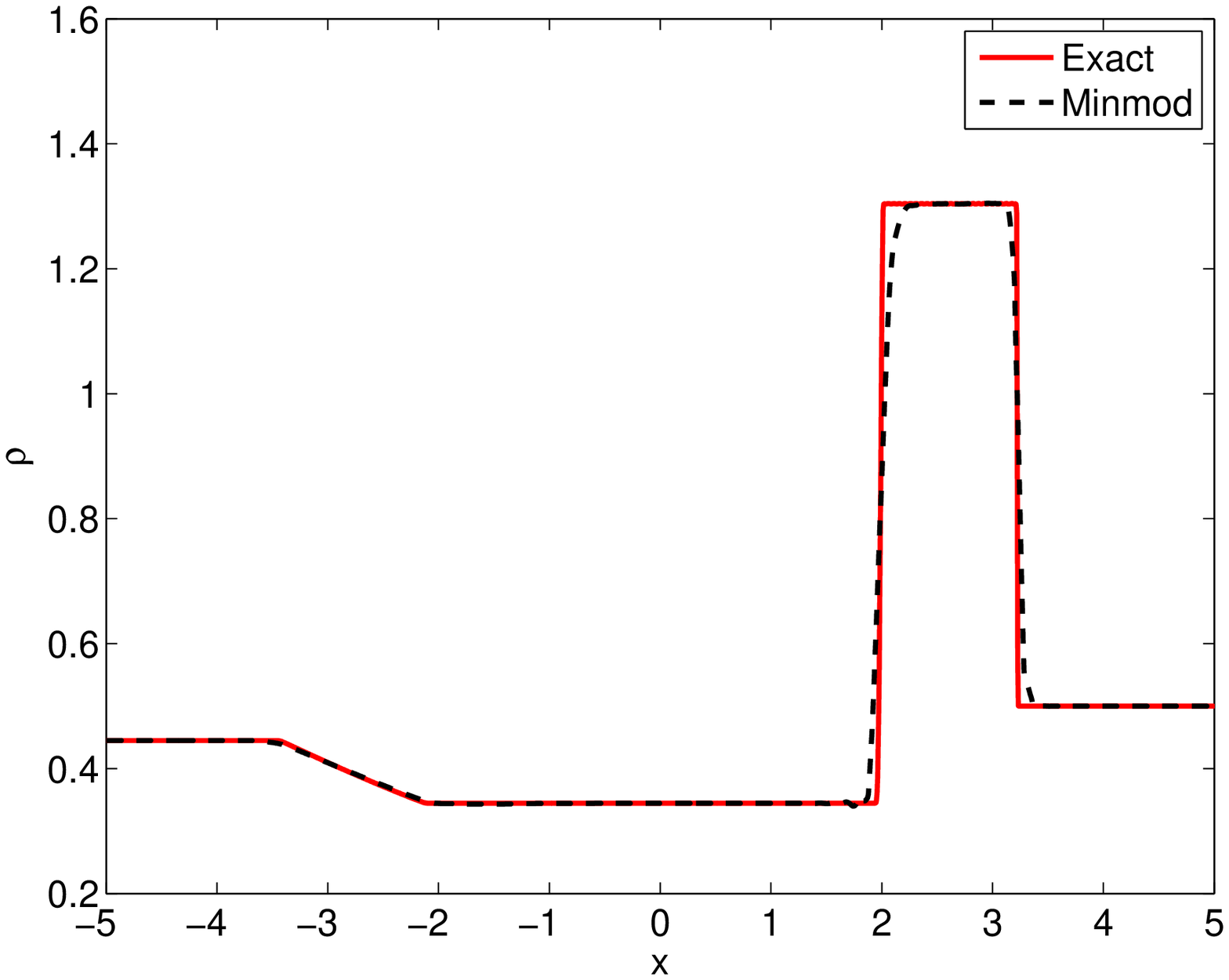}} \\
\vspace{-0.3cm}
\caption{Detected troubled cells (row 1 and 2) and approximation at final time $T=1.3$ (row 3 and 4), shock tube of Lax, $k=2$, 128 elements.}\label{fig:Laxk2}
\end{figure}

\begin{figure}[h!]
 \centering
\subfigure[Original, $C=0.05$]{\includegraphics[scale = 0.26]{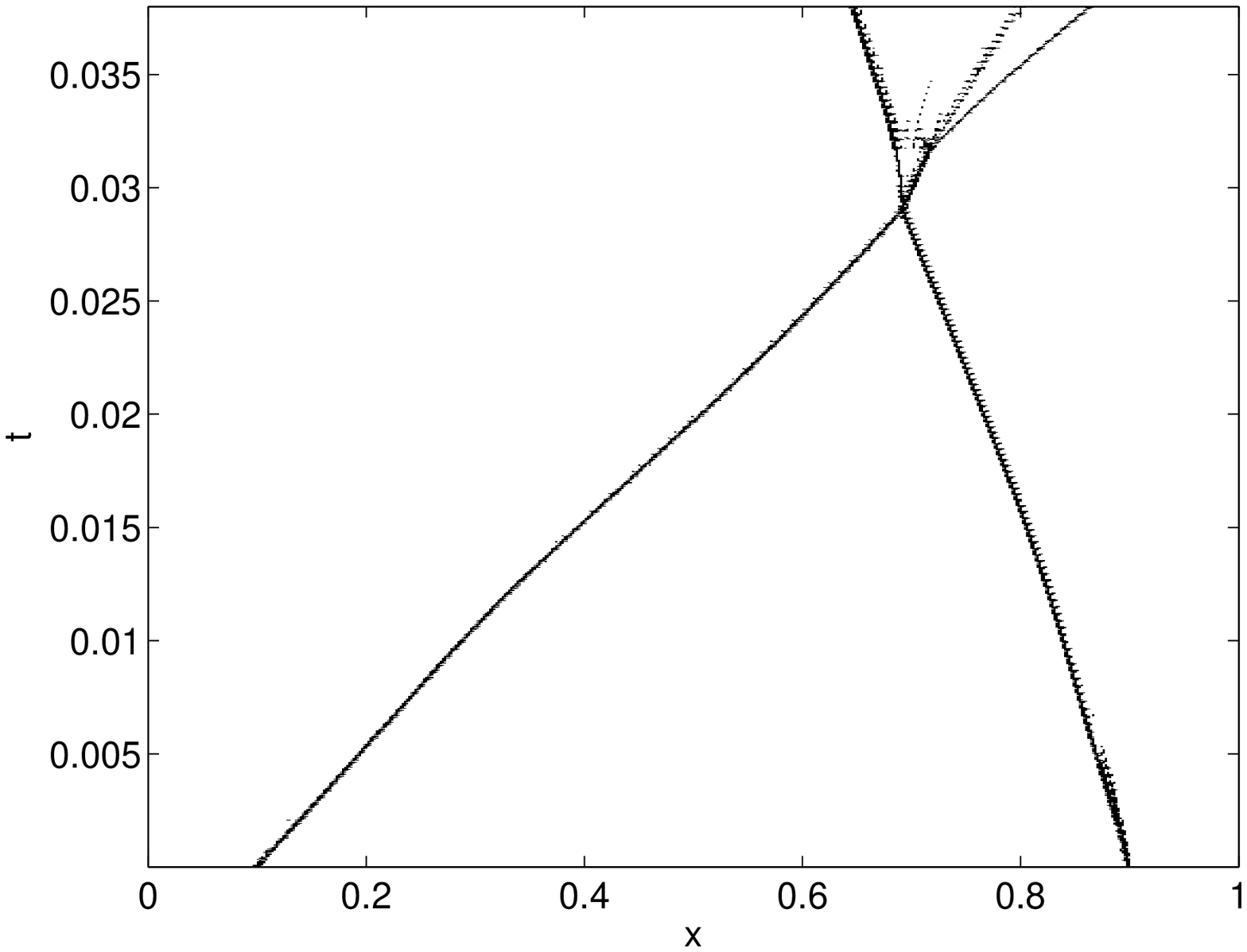}}
\subfigure[Original, KXRCF]{\includegraphics[scale = 0.26]{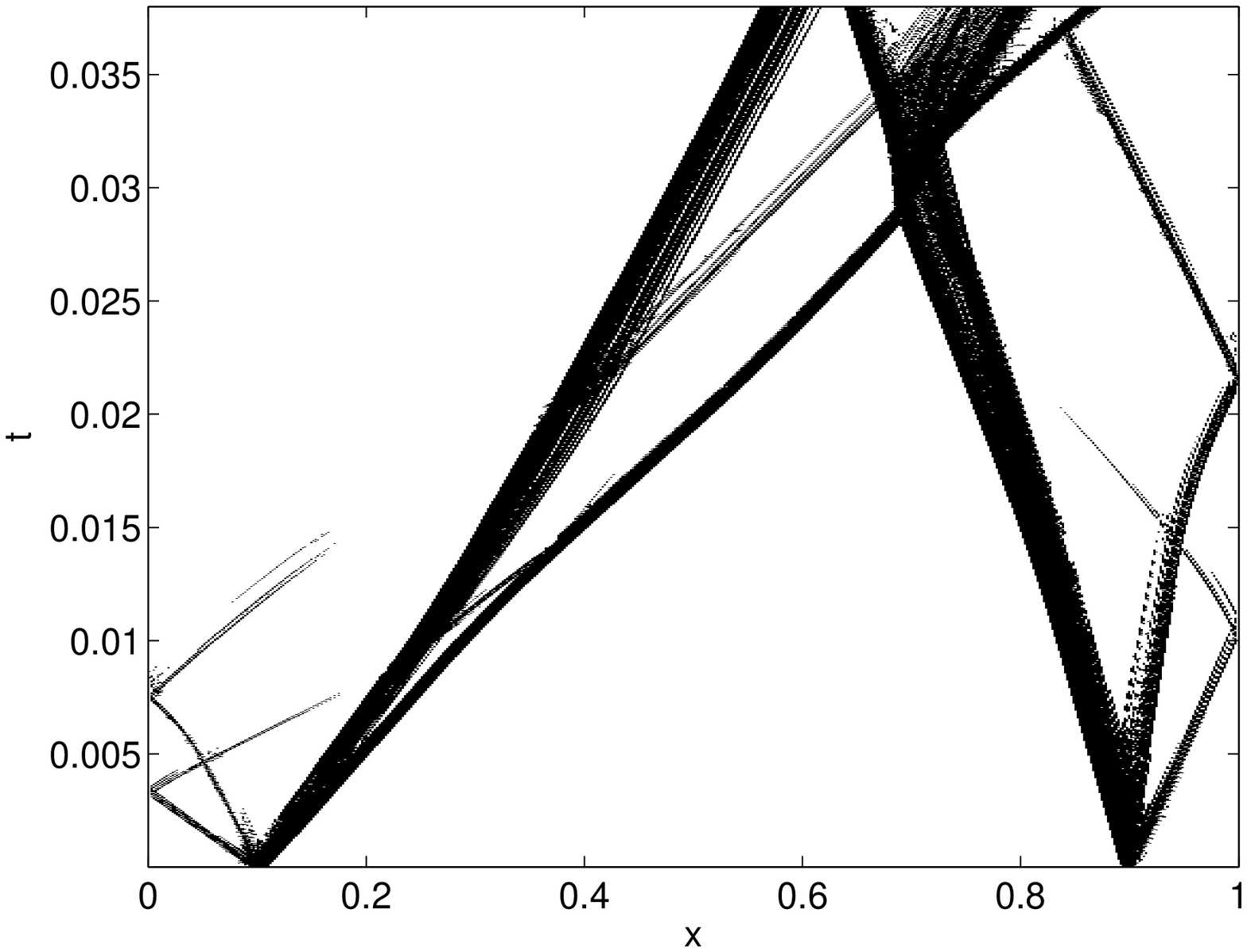}}
\subfigure[Original, $M = 100$]{\includegraphics[scale = 0.26]{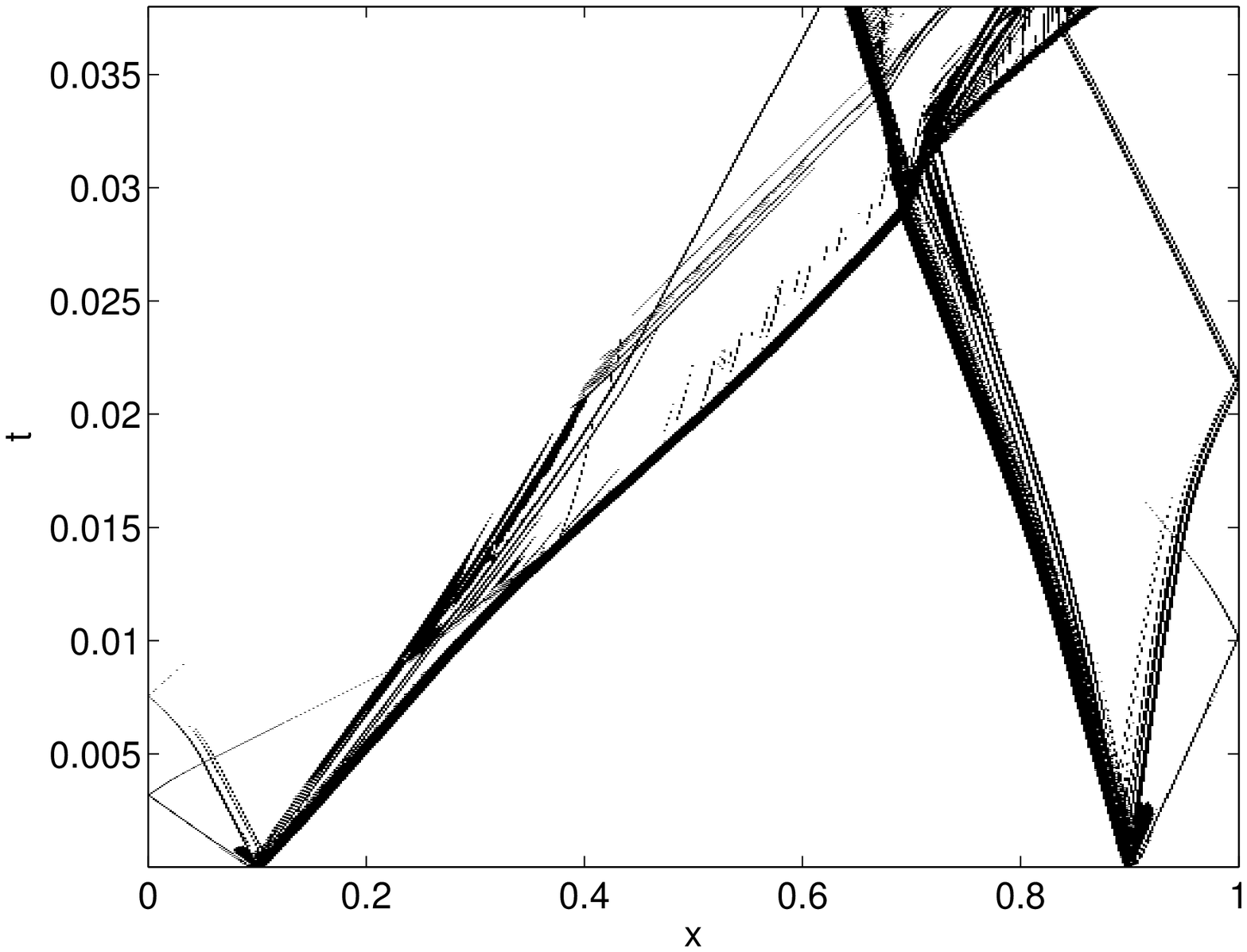}} \\
\vspace{-0.3cm}
\subfigure[Outlier, multiwavelets]{\includegraphics[scale = 0.26]{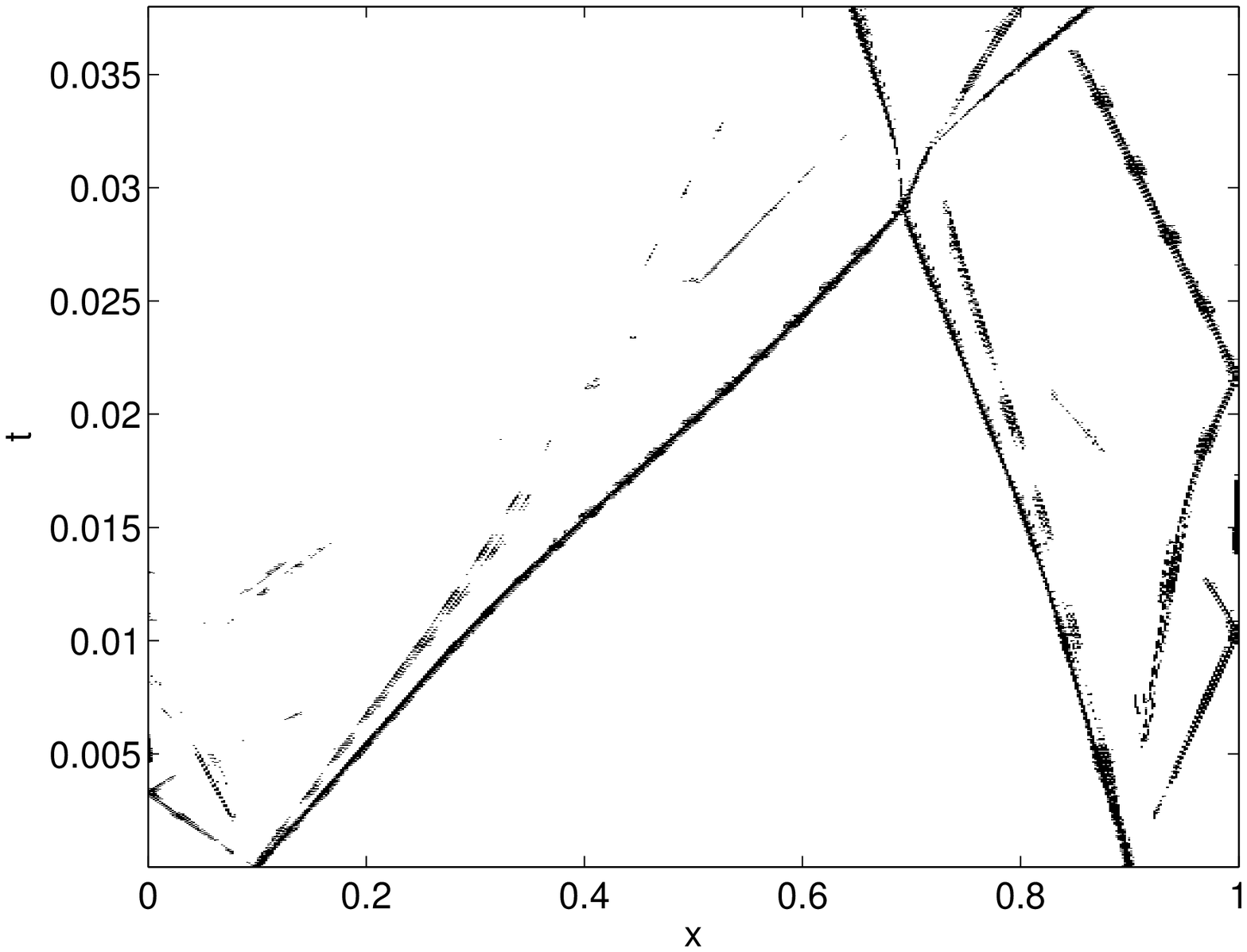}}
\subfigure[Outlier, KXRCF value]{\includegraphics[scale = 0.26]{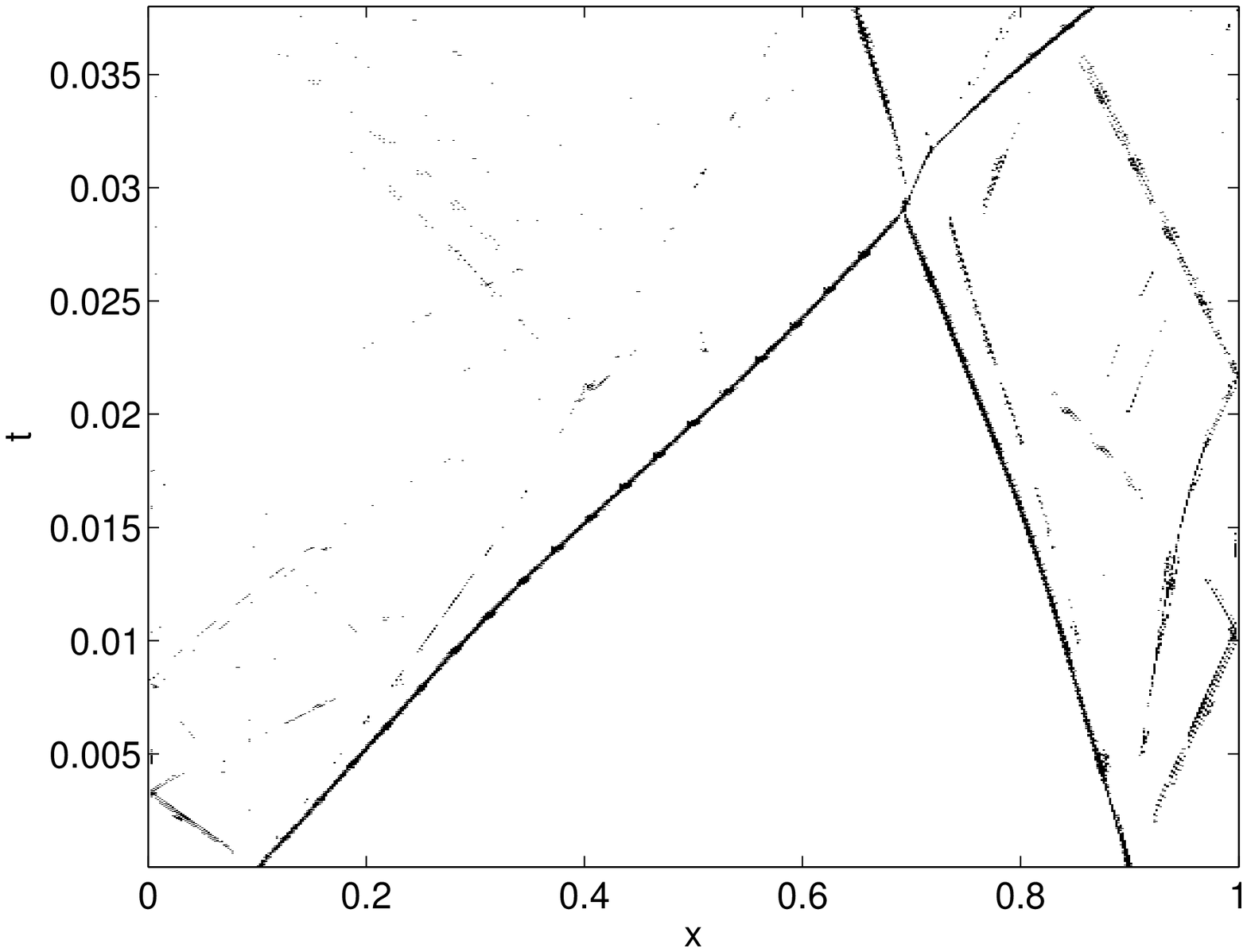}} 
\subfigure[Outlier, minmod-based TVB]{\includegraphics[scale = 0.26]{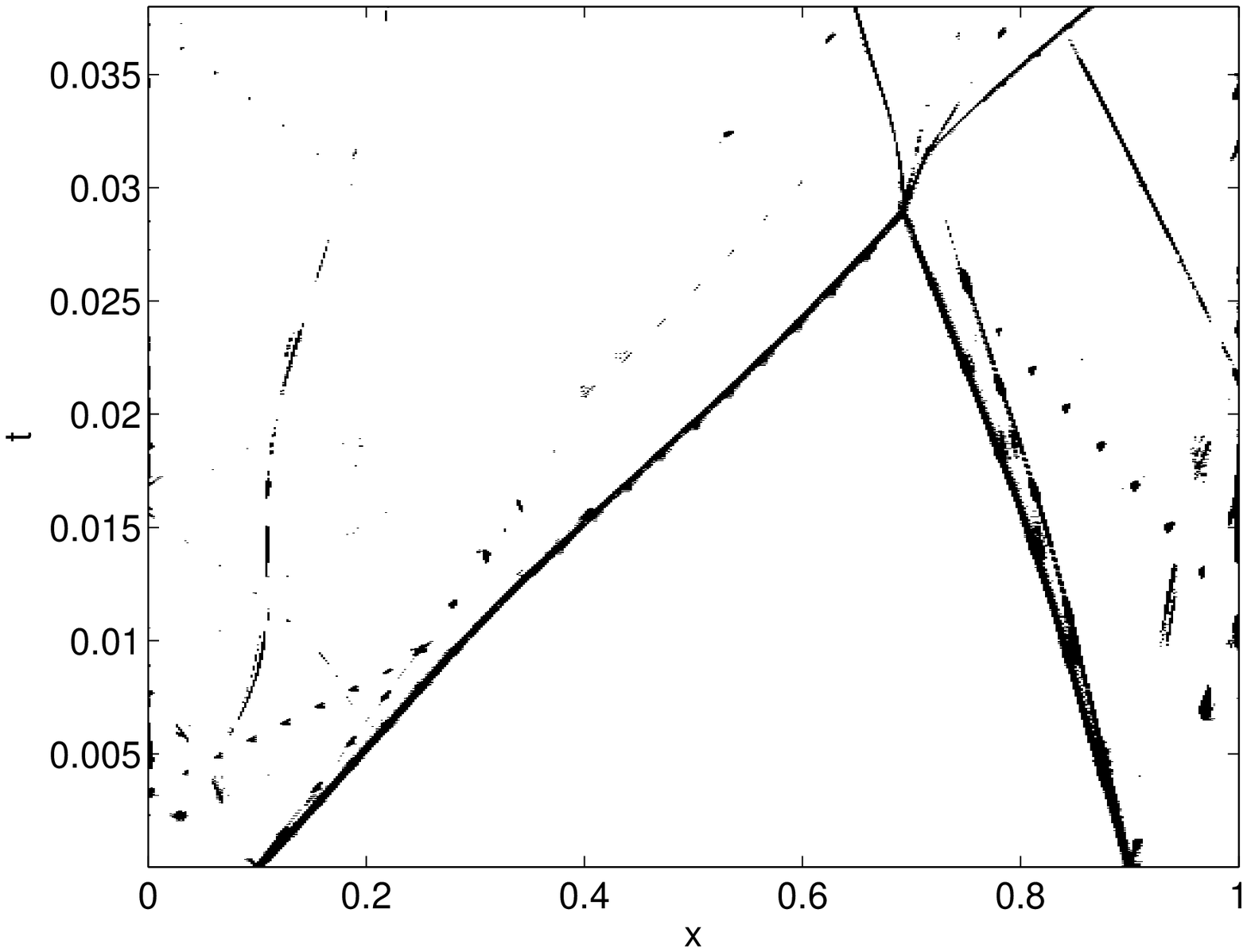}} \\
\vspace{-0.3cm}
\subfigure[Original, $C=0.05$]{\includegraphics[scale = 0.26]{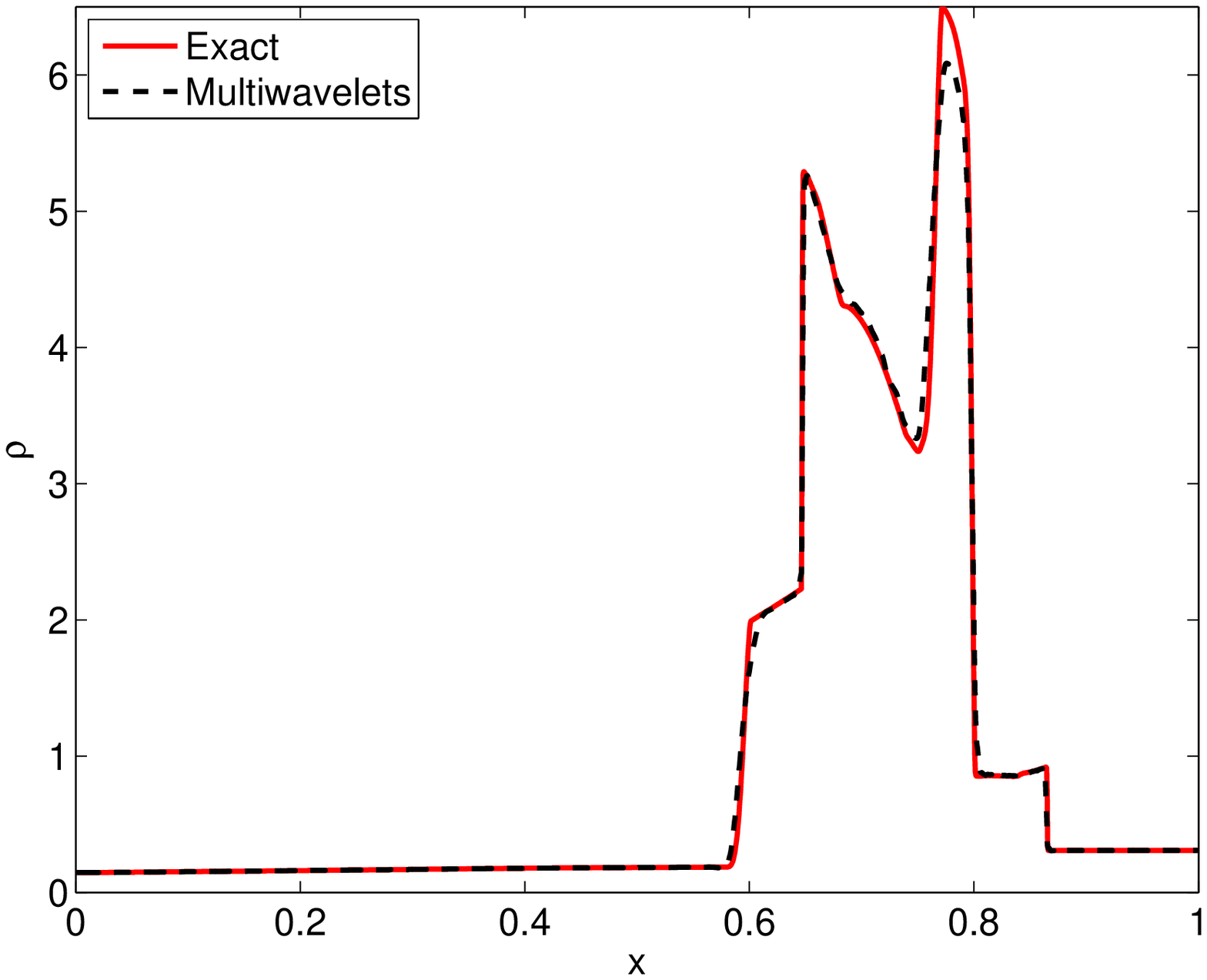}}
\subfigure[Original, KXRCF]{\includegraphics[scale = 0.26]{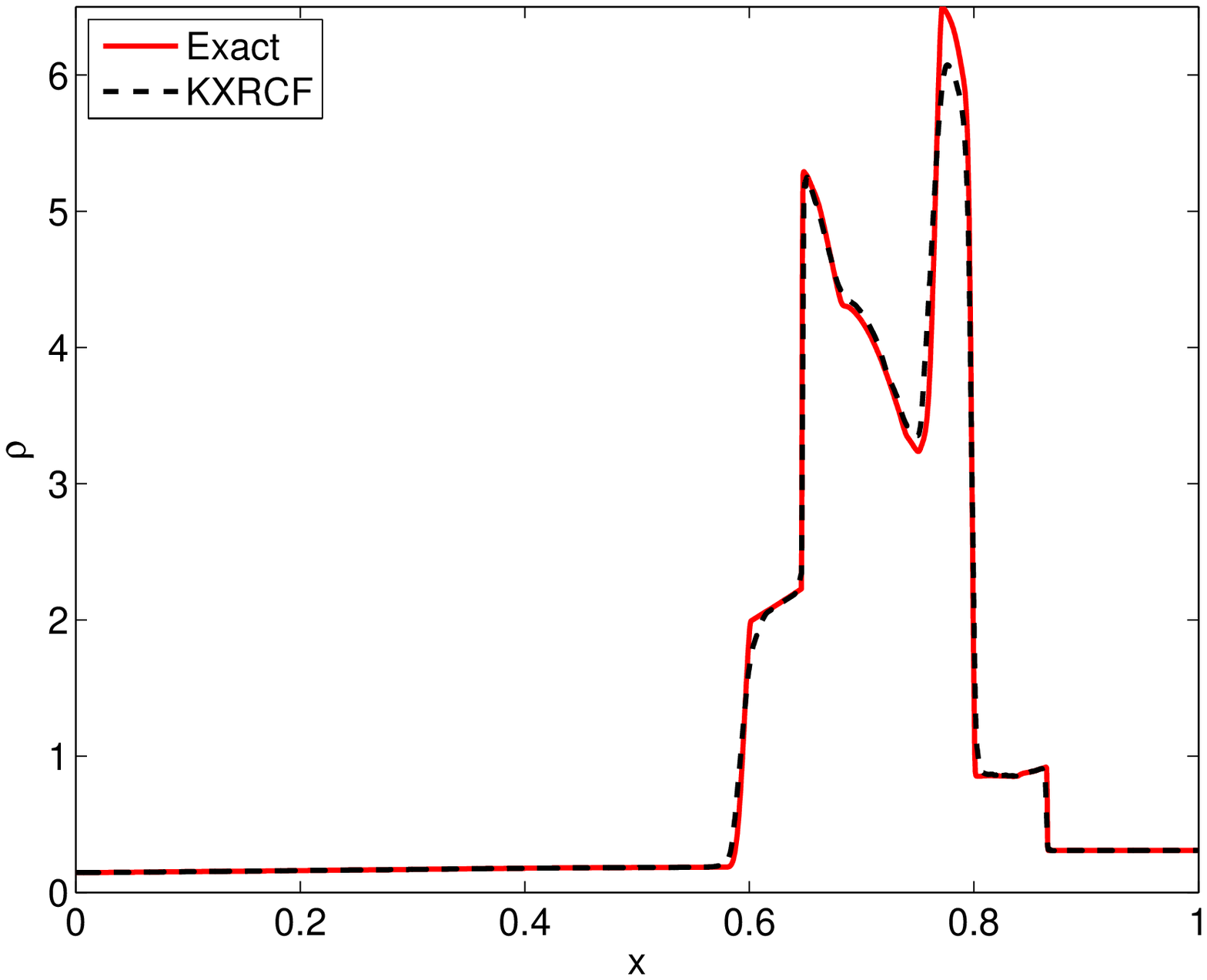}}
\subfigure[Original, $M = 100$]{\includegraphics[scale = 0.26]{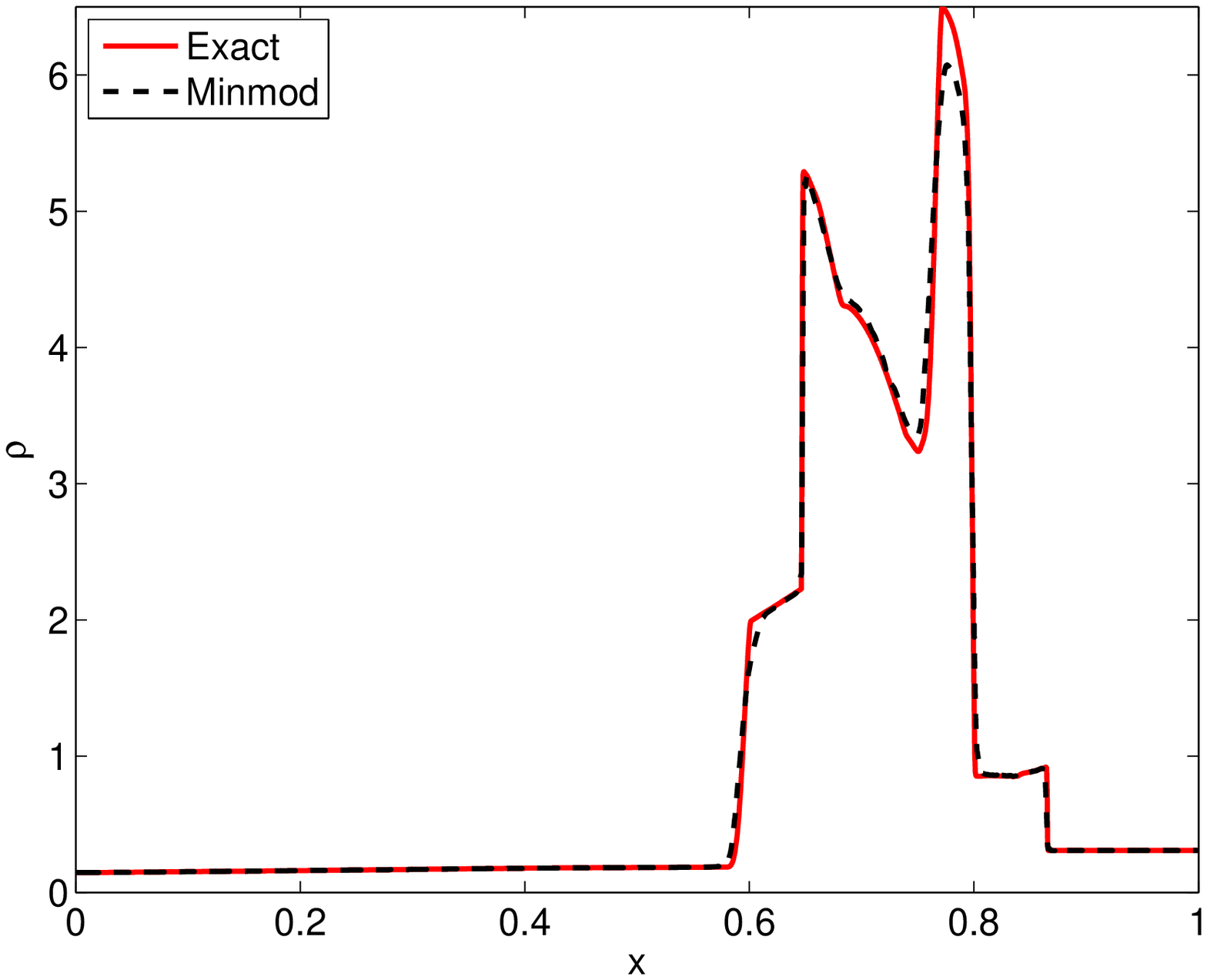}} \\
\vspace{-0.3cm}
\subfigure[Outlier, multiwavelets]{\includegraphics[scale = 0.26]{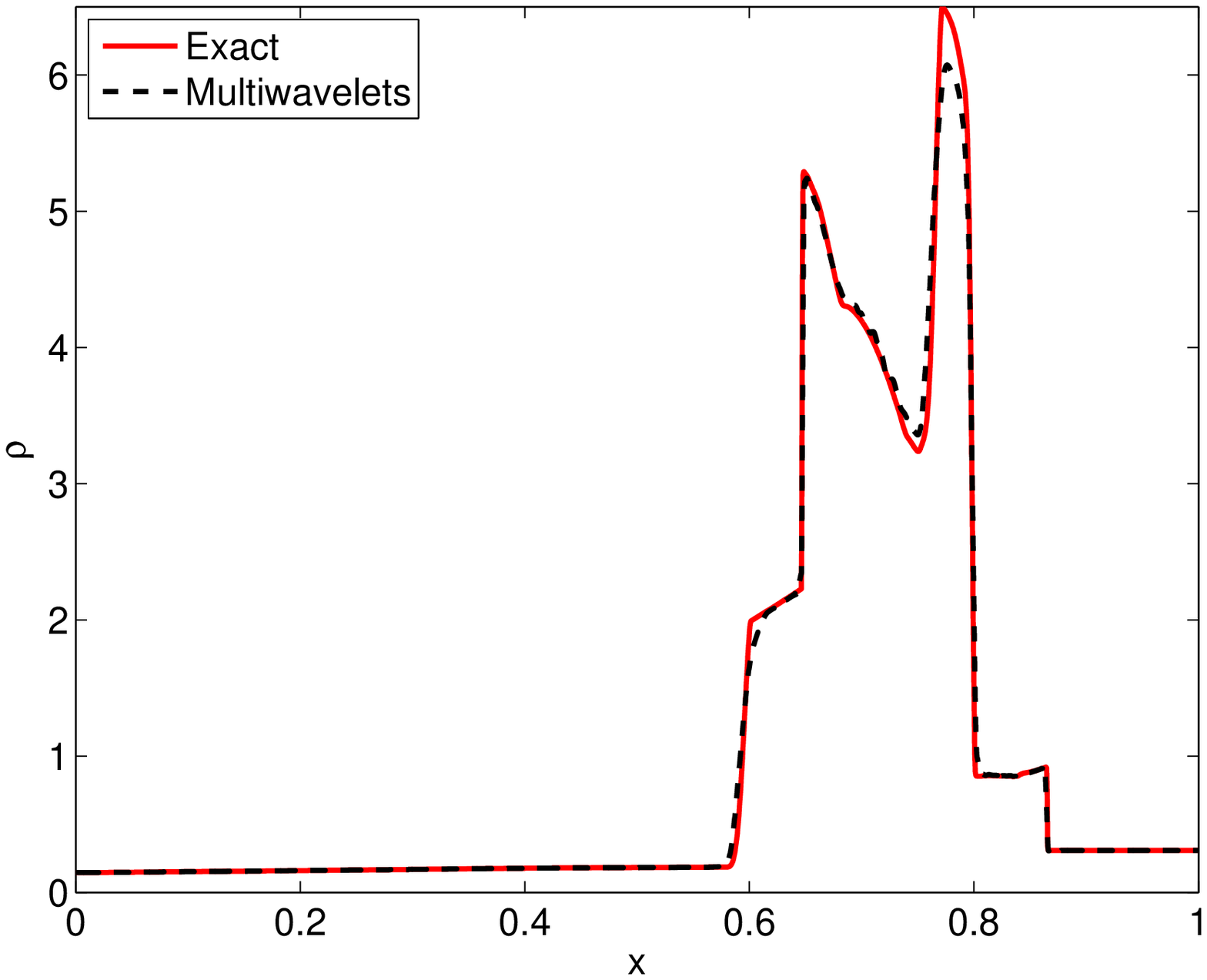}}
\subfigure[Outlier, KXRCF value]{\includegraphics[scale = 0.26]{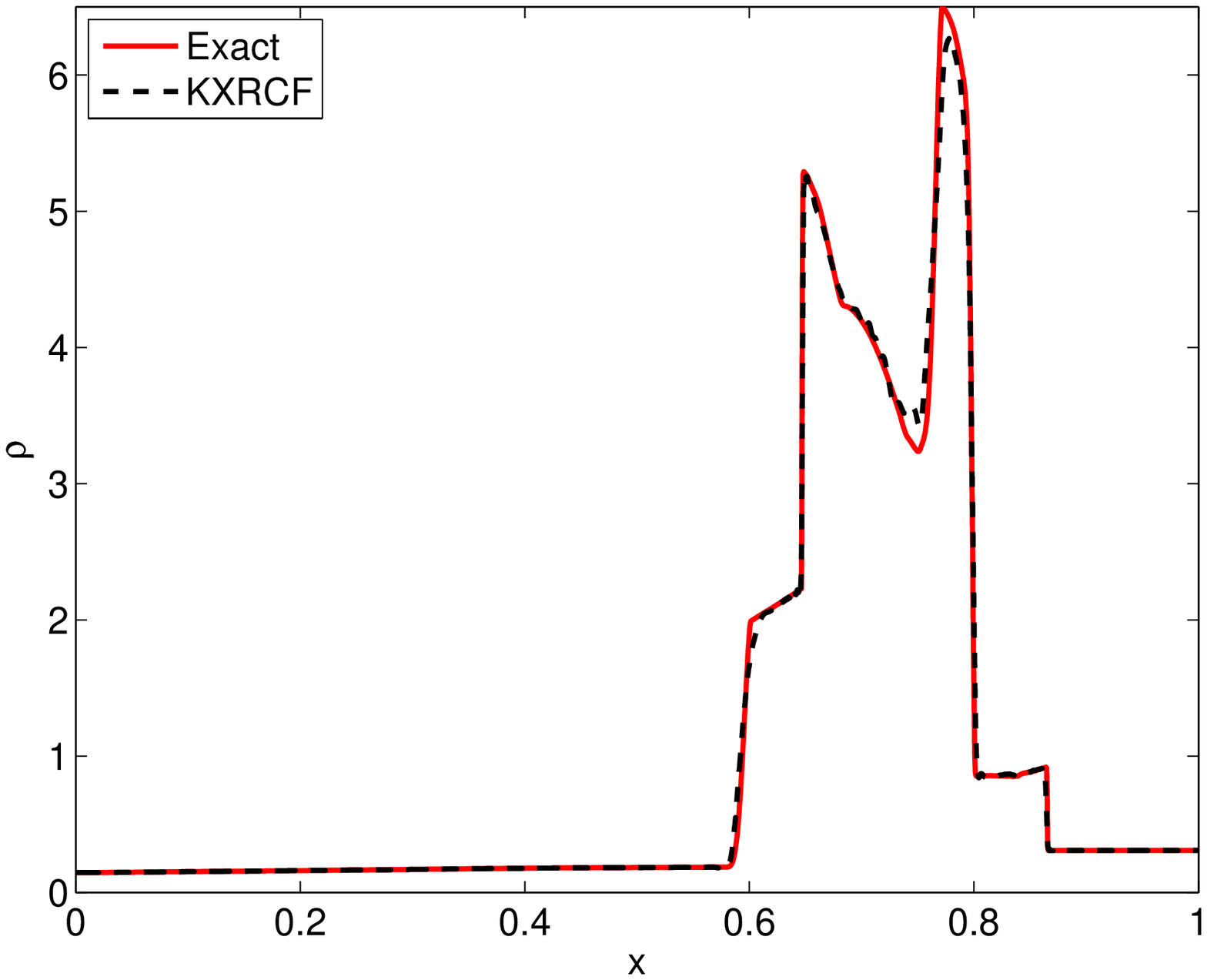}} 
\subfigure[Outlier, minmod-based TVB]{\includegraphics[scale = 0.26]{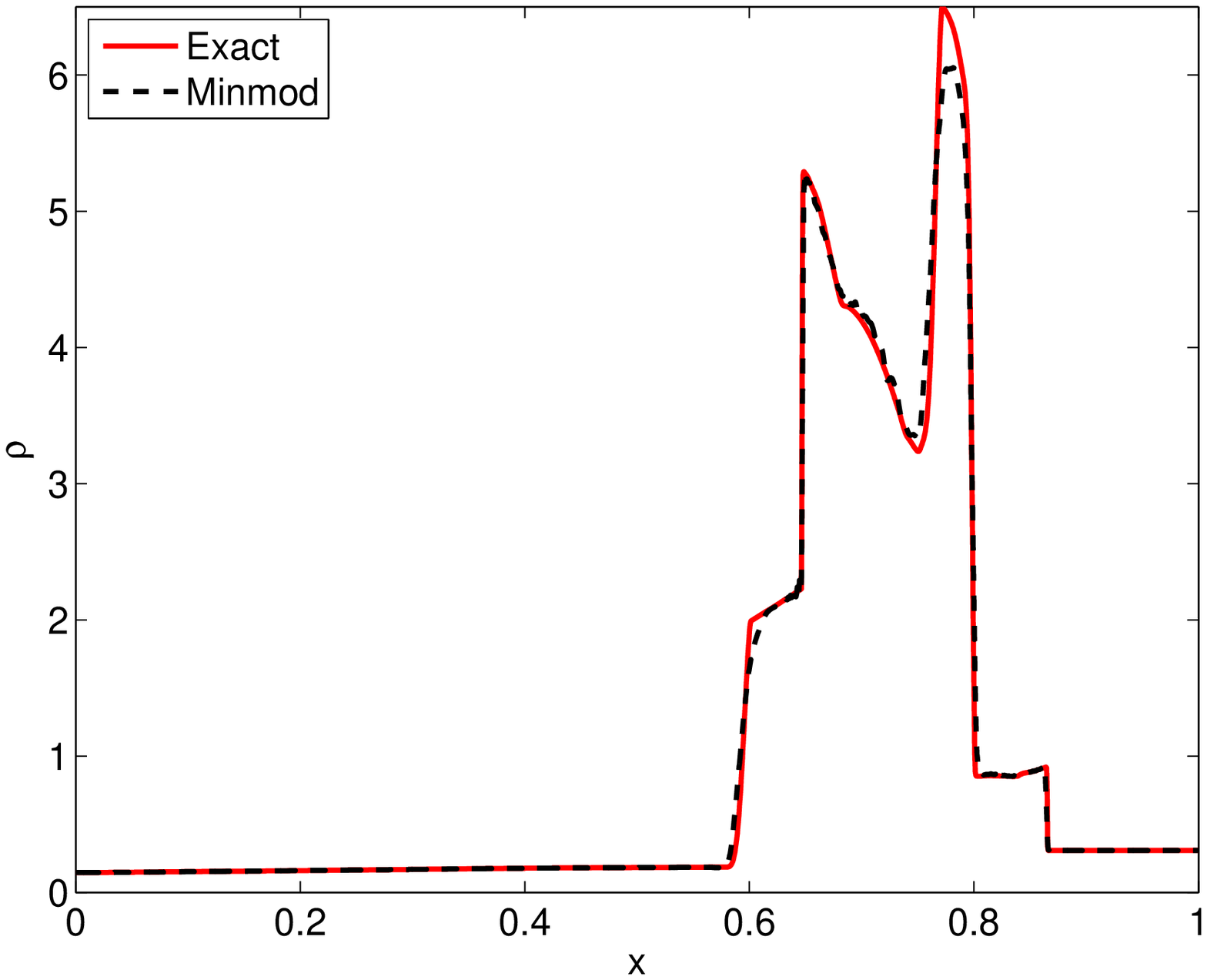}} \\
\vspace{-0.3cm}
\caption{Detected troubled cells (row 1 and 2) and approximation at final time $T=0.038$ (row 3 and 4), blast-wave problem, $k=2$, 512 elements.}\label{fig:Blastk2}
\end{figure}

\begin{figure}[h!]
\subfigure[Original, $C=0.01$]{\includegraphics[scale = 0.27]{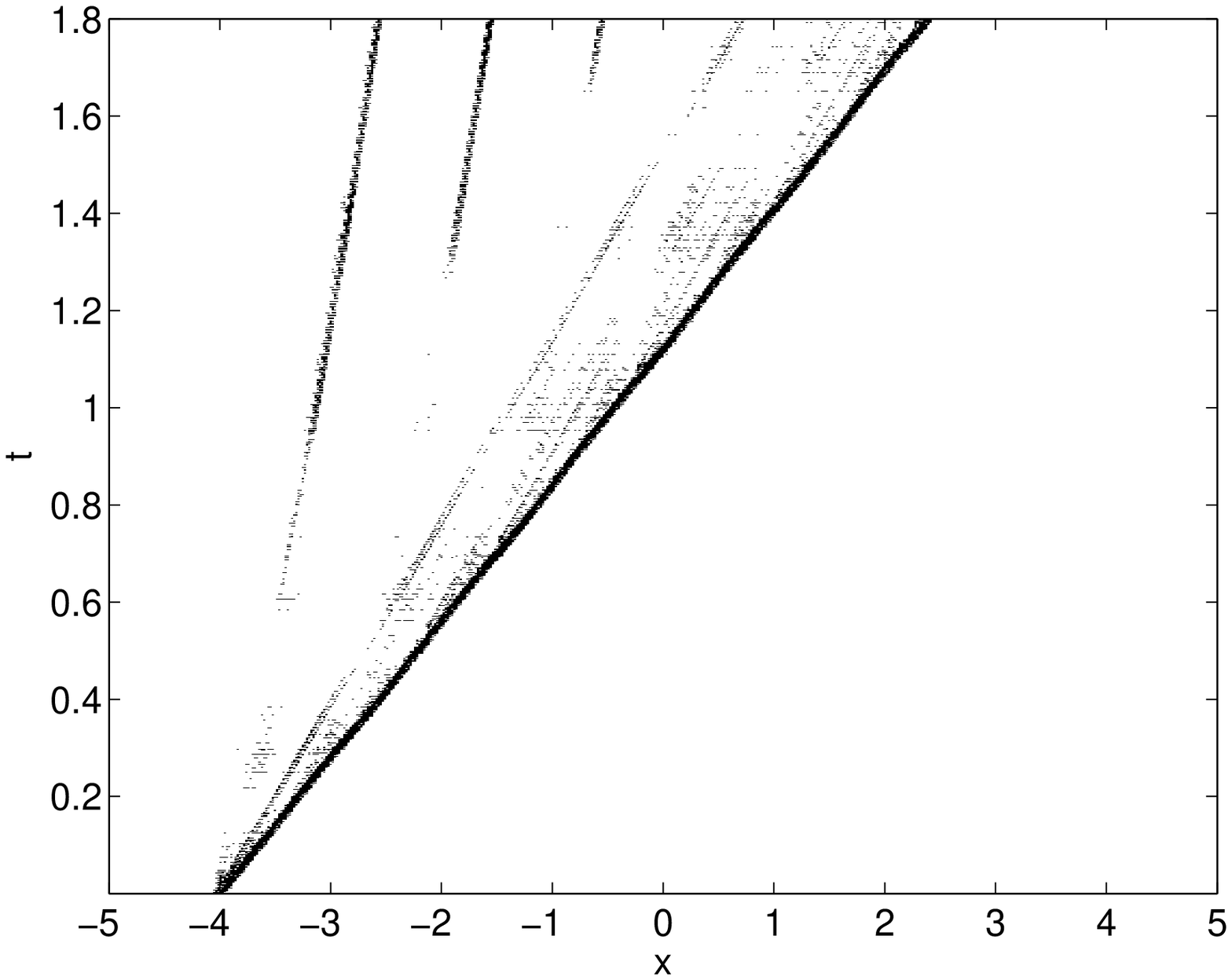}}
\subfigure[Original, KXRCF]{\includegraphics[scale = 0.27]{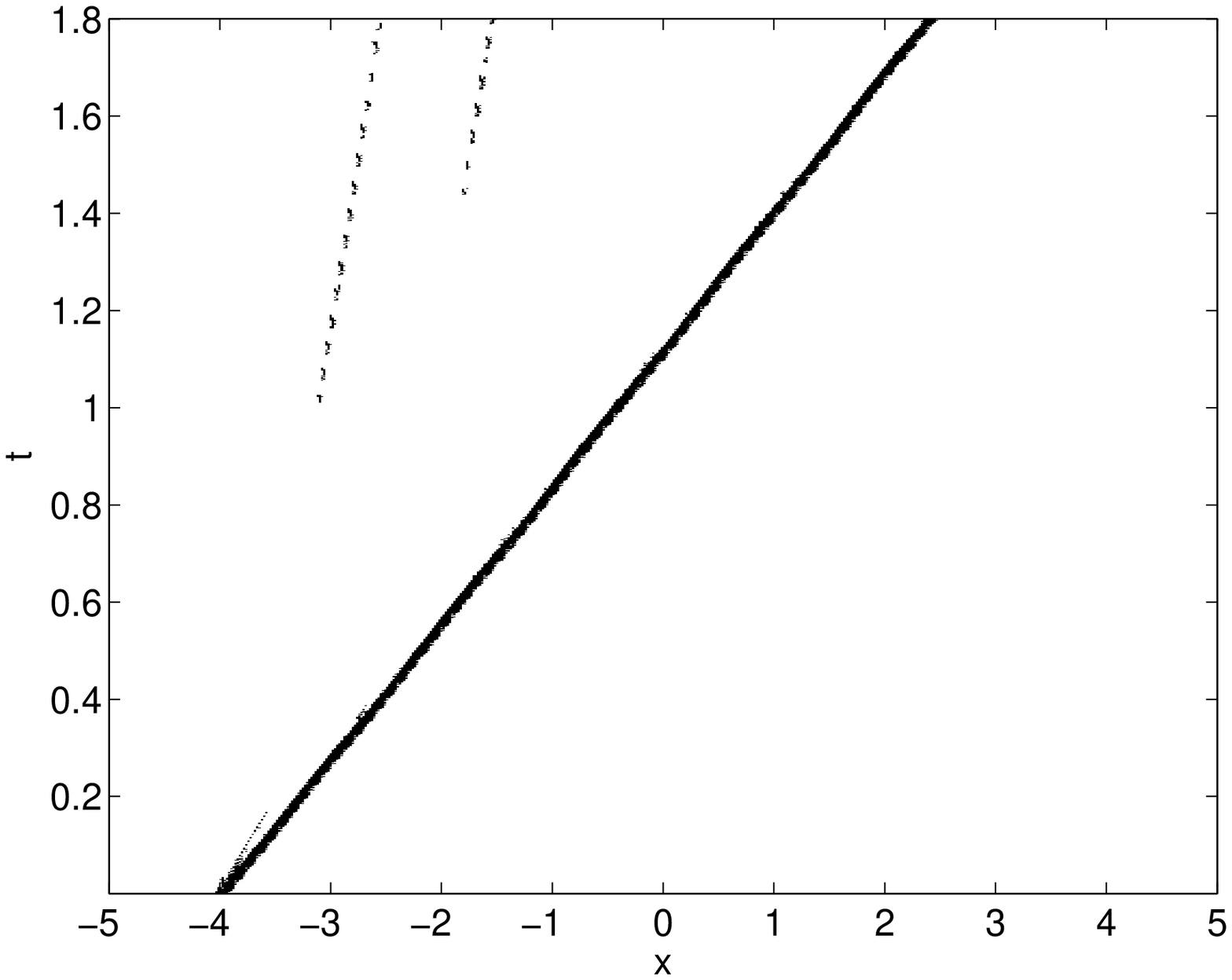}}
\subfigure[Original, $M = 100$]{\includegraphics[scale = 0.27]{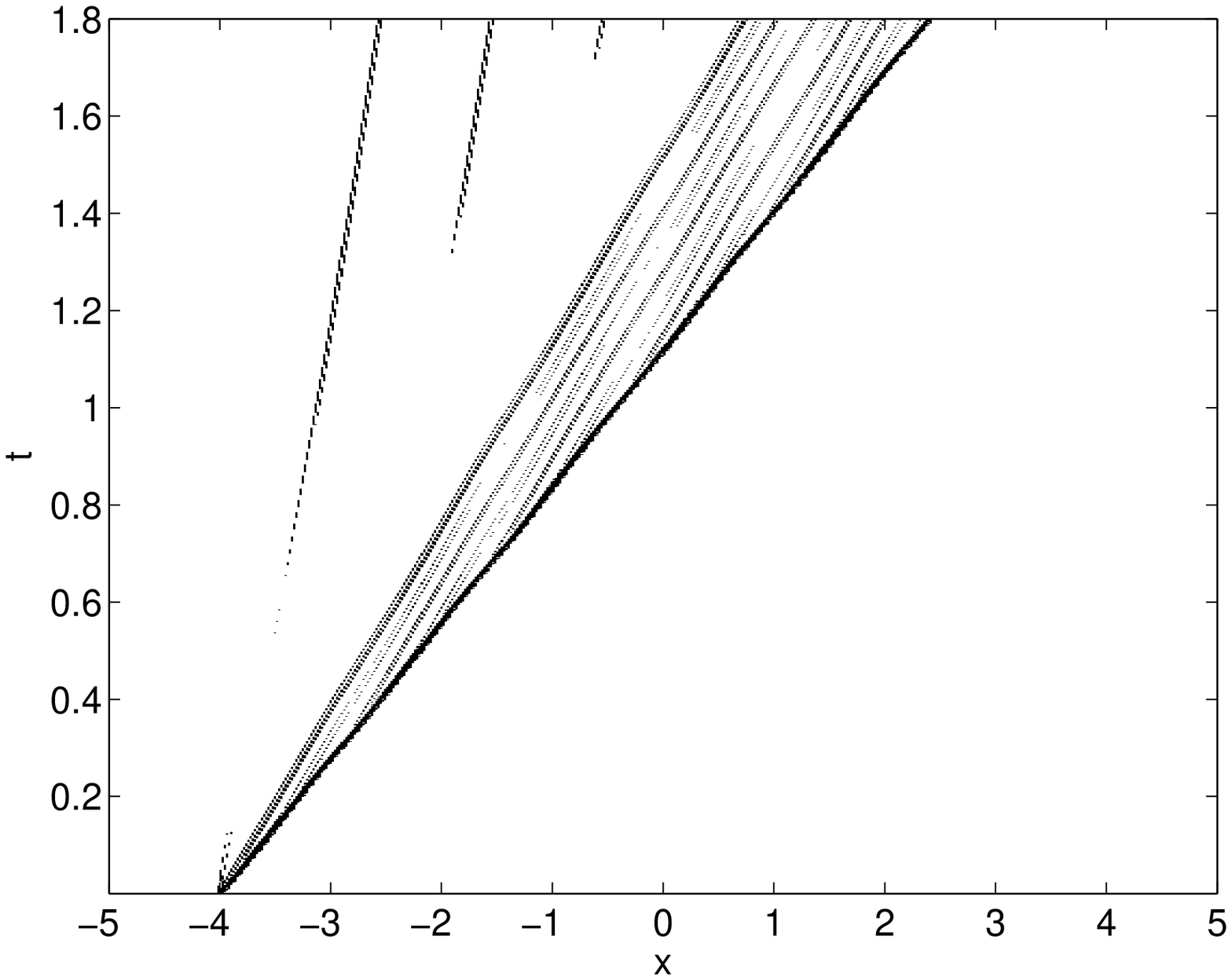}} \\
\vspace{-0.2cm}
\subfigure[Outlier, multiwavelets]{\includegraphics[scale = 0.27]{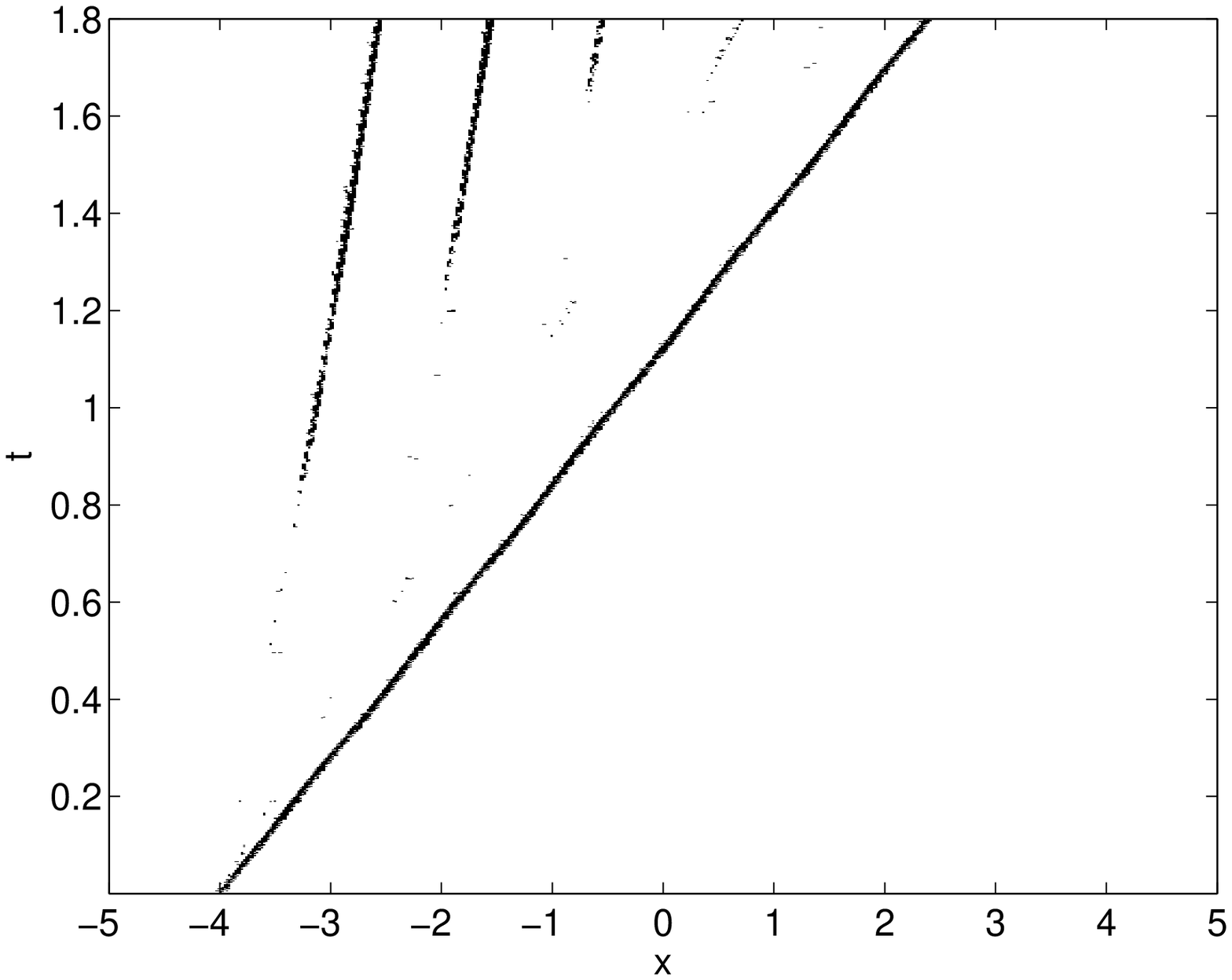}}
\subfigure[Outlier, KXRCF value]{\includegraphics[scale = 0.27]{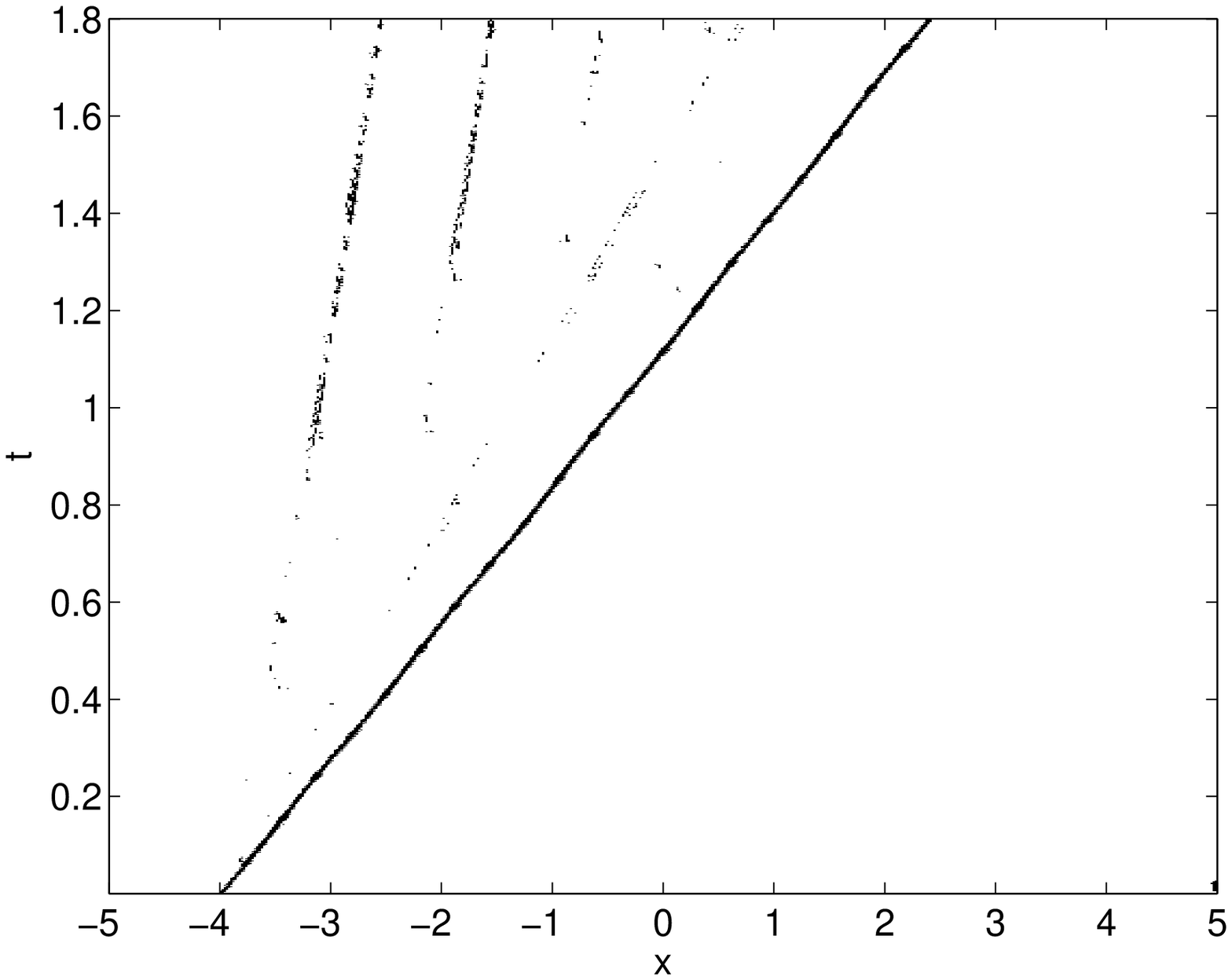}}
\subfigure[Outlier, minmod-TVB]{\includegraphics[scale = 0.27]{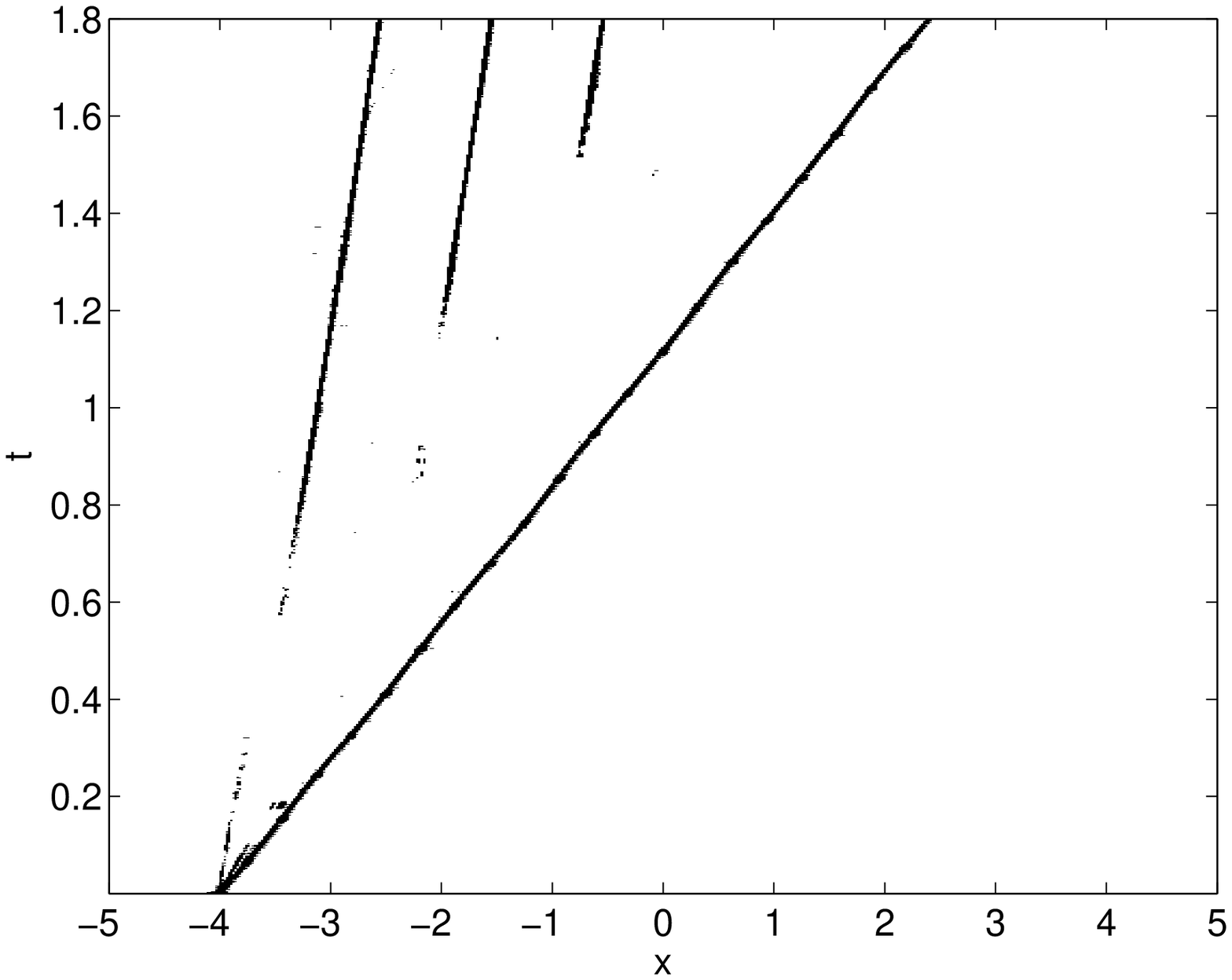}} \\
\vspace{-0.2cm}
\subfigure[Original, $C=0.01$]{\includegraphics[scale = 0.27]{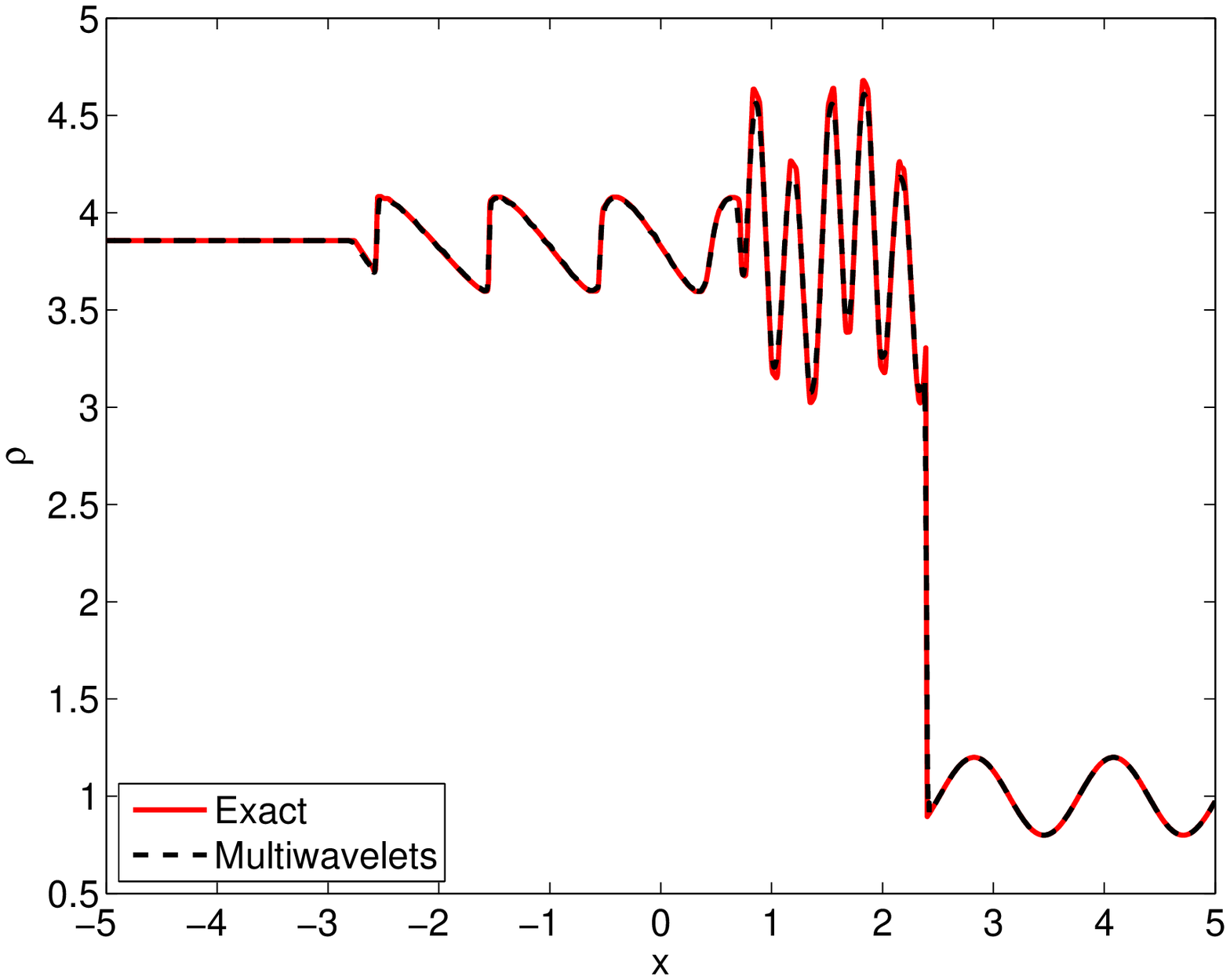}}
\subfigure[Original, KXRCF]{\includegraphics[scale = 0.27]{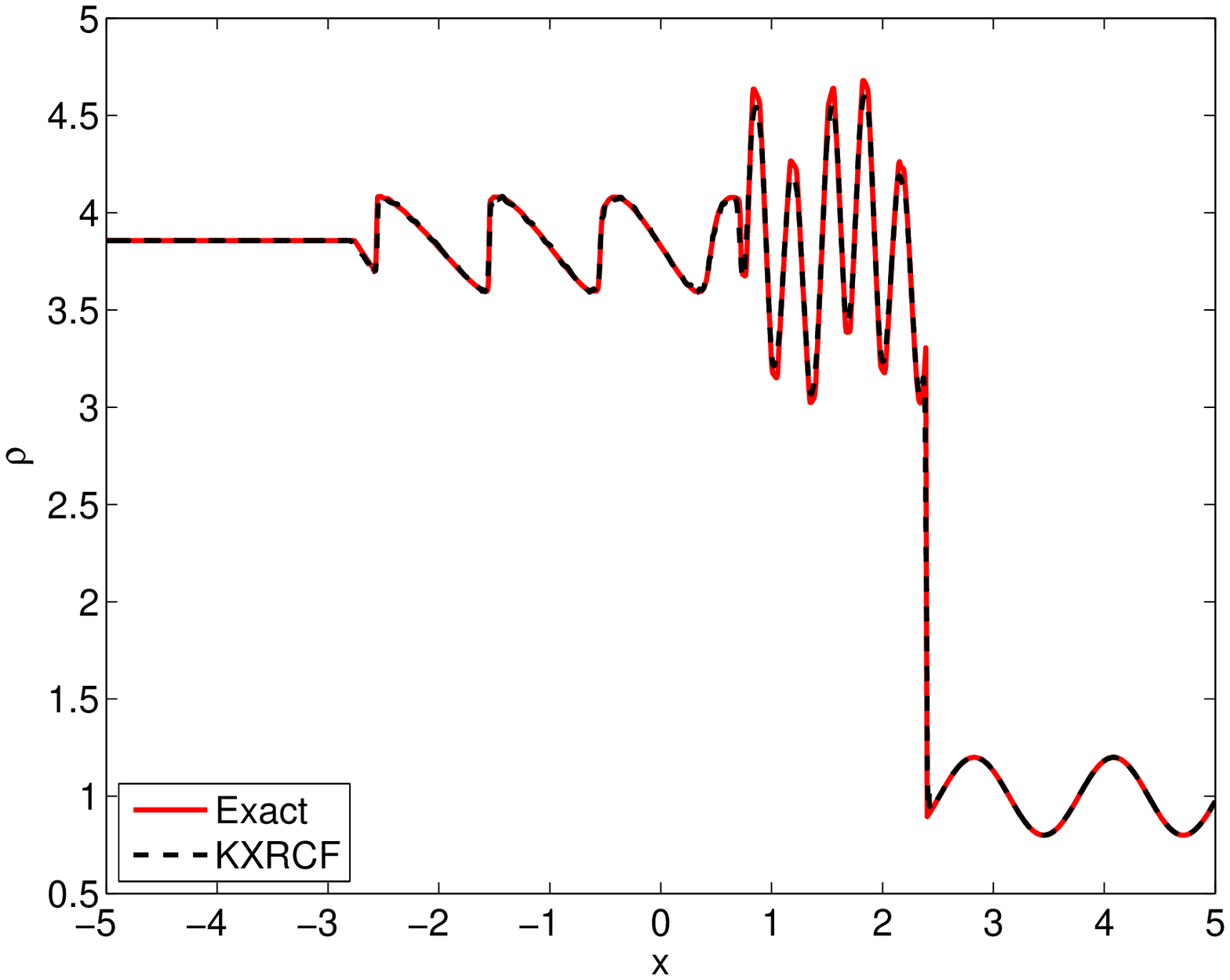}}
\subfigure[Original, $M = 100$]{\includegraphics[scale = 0.27]{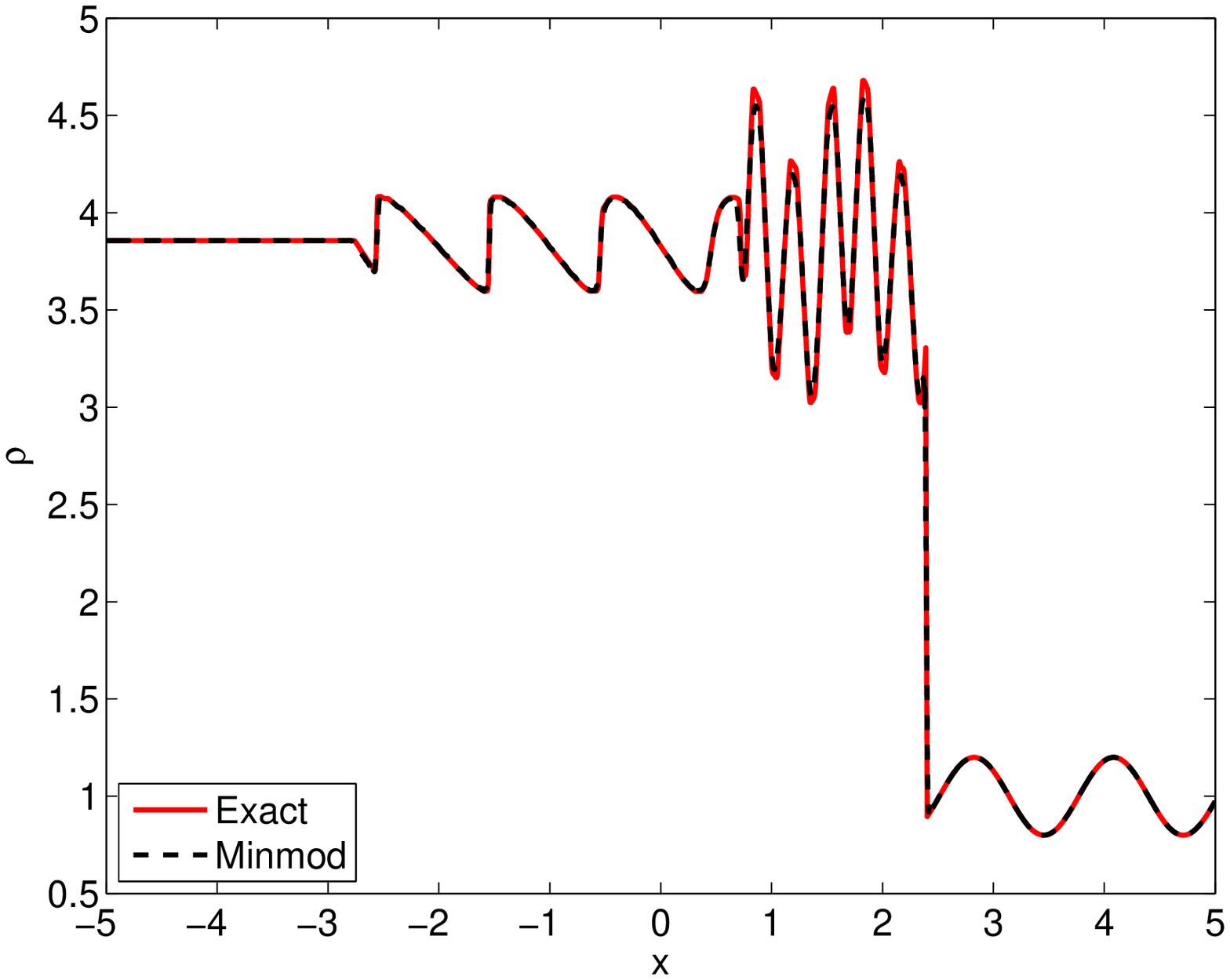}} \\
\vspace{-0.2cm}
\subfigure[Outlier, multiwavelets]{\includegraphics[scale = 0.27]{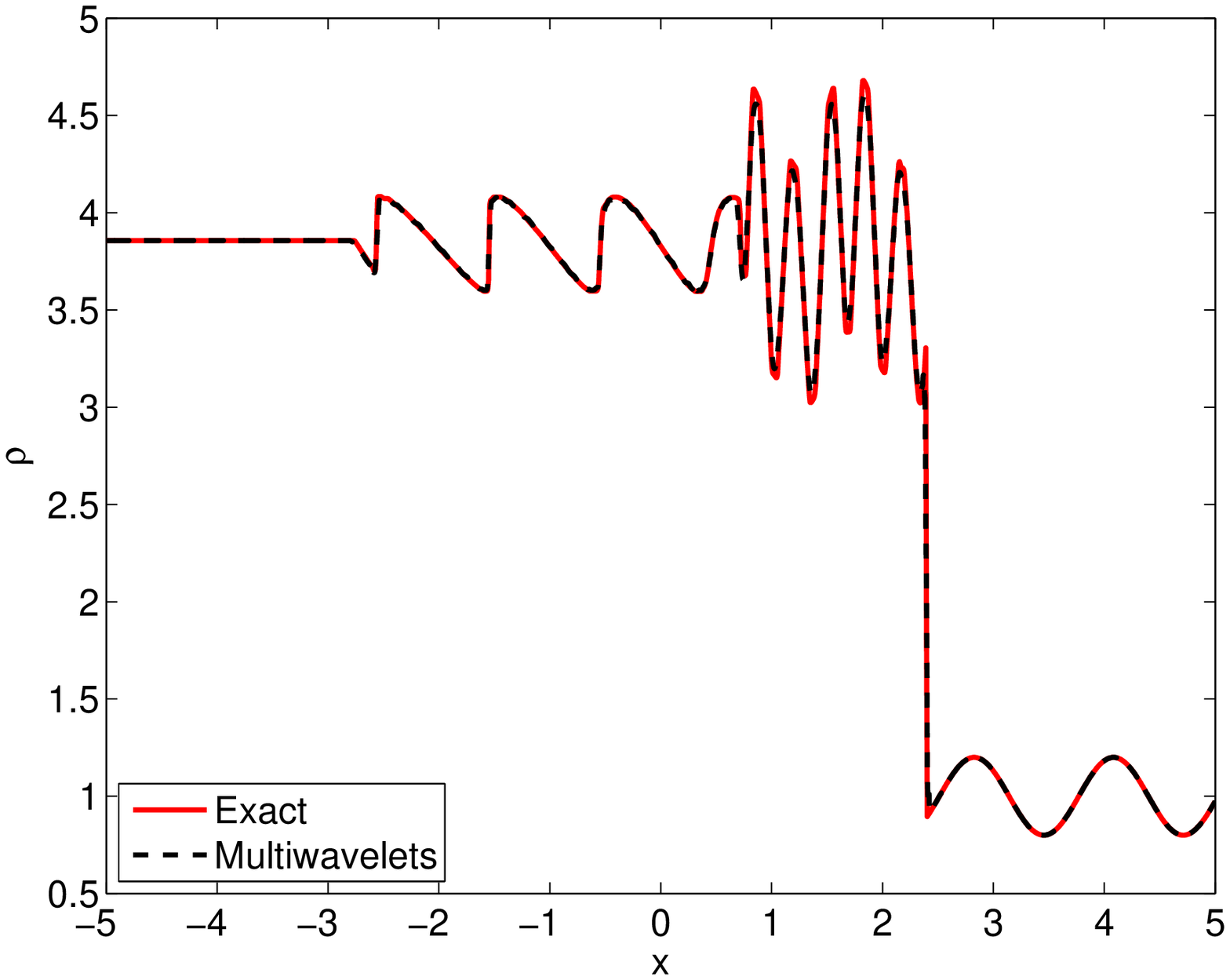}}
\subfigure[Outlier, KXRCF value]{\includegraphics[scale = 0.27]{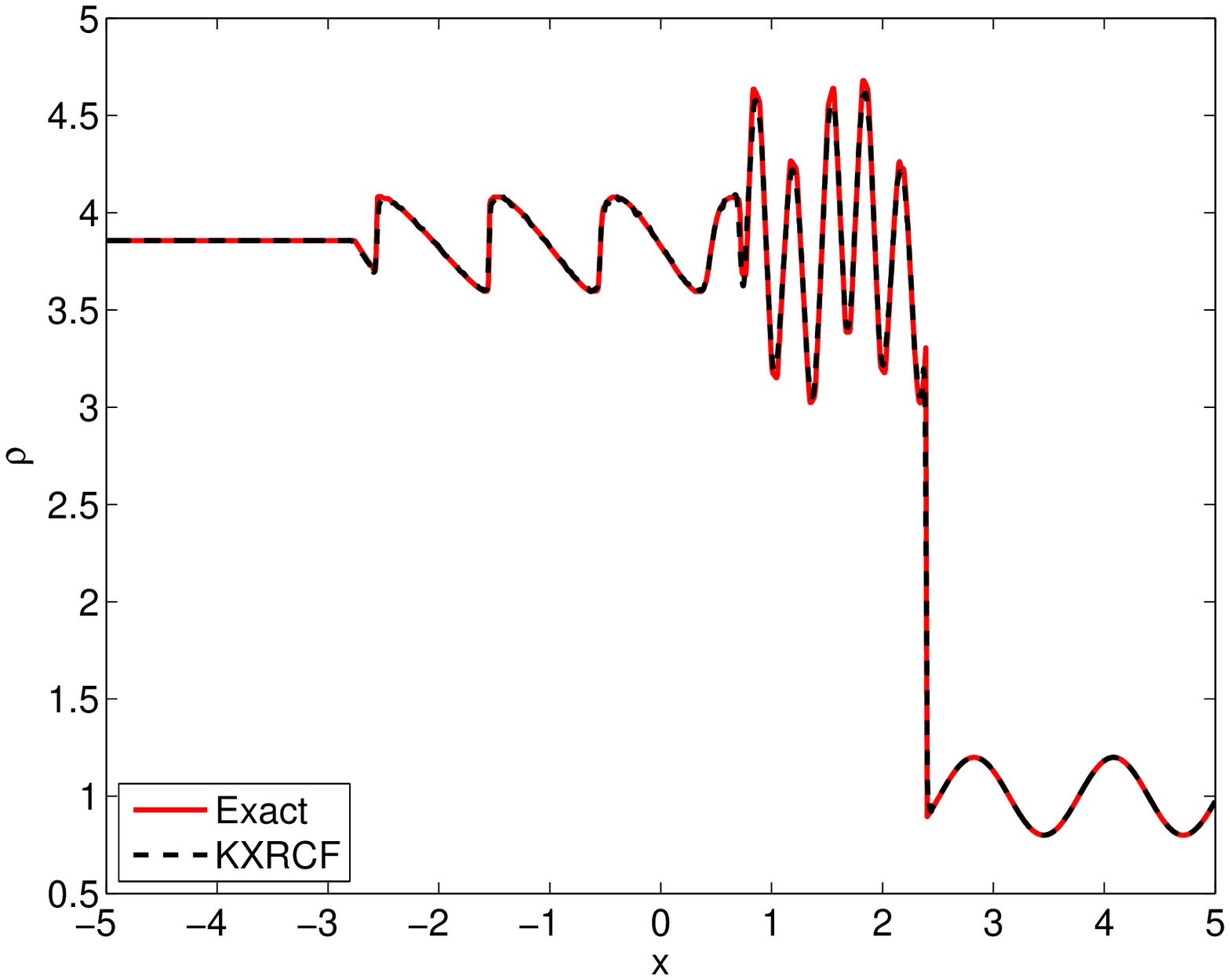}}
\subfigure[Outlier, minmod-TVB]{\includegraphics[scale = 0.27]{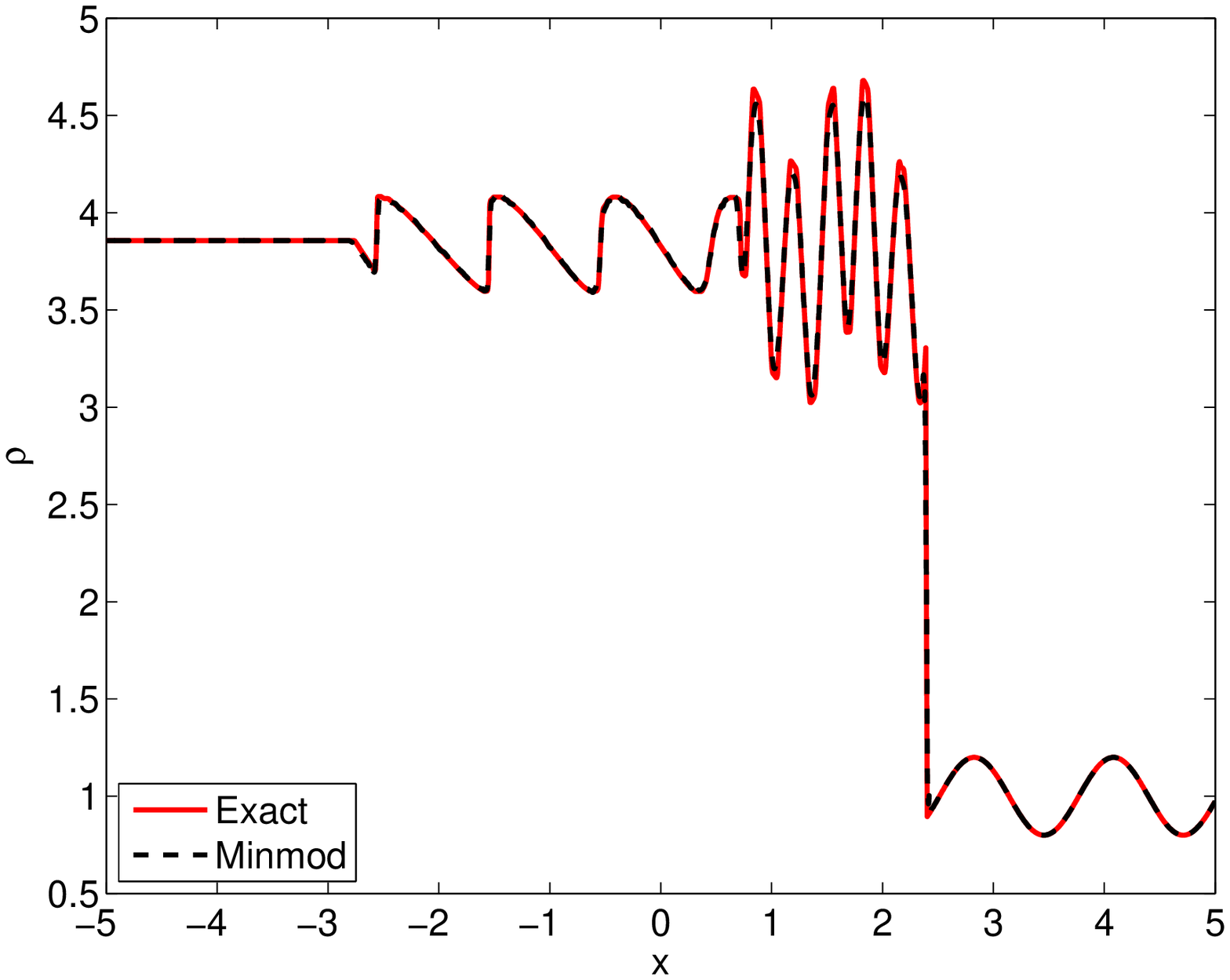}} \\
\vspace{-0.2cm}
\caption{Detected troubled cells (row 1 and 2) and approximation at final time $T=1.8$ (row 3 and 4), Shu-Osher problem, $k=2$, 512 elements.}\label{fig:Sinek2}
\end{figure}

Below the results and comparisons are shown using four different sets of initial conditions: the shock tubes of Sod \cite{Sod78} (Figure \ref{fig:Sodk2}) and Lax \cite{Lax54} (Figure \ref{fig:Laxk2}), the blast-wave problem \cite{Woo84C} (Figure \ref{fig:Blastk2}), and the Shu-Osher problem \cite{Shu89O} (Figure \ref{fig:Sinek2}).  We omit the details of these test problems and refer to \cite{Vui14R} for more information on initial conditions and boundary conditions. We apply the indication technique to density for the modified multiwavelet indicator, density and energy for KXRCF, and the characteristic variables for the minmod-based TVB indicator, as has been done by Qiu et al. \cite{Qiu05S_2}. The first row of each figure consists of time-history plots of detected troubled cells using the original indicators. The second row belongs to the outlier-detected troubled cells. The corresponding approximations at the final times are given in the third and fourth row. Note that these results are computed using the moment limiter in the detected troubled cells. A different choice for the limiter will result in different approximations. In all figures we take $k=2$, and similar results were found for $k=1,3$.

Note that the original troubled-cell indicators are applied using the optimal problem-dependent parameters as found in \cite{Qiu05S_2, Vui14R}.  We stress that the outlier-detected results are computed without problem-dependent parameters, but with a fixed whisker length equal to 3, and with local indication vectors of size 16.

It turns out that the new outlier-detection approach detects the troubled regions very accurately and generally better than the original parameter-using methods for the blast-wave and Shu-Osher problem. For the shock tube problems of Sod and Lax, most discontinuous regions are selected. Note that the outlier-detection indicators sometimes detect jumps in derivatives, as can be seen at the end points of the rarefaction waves. The original indicators, however, do not detect these structures. This difference can be explained by recalling that the original indicators focus on the actual value of the indication variable, whereas the outlier-detection techniques investigate the relative value with respect to the neighboring region. A discontinuity in the derivative usually causes sudden differences, and therefore these regions are detected as troubled. By applying a limiter at these locations, the discontinuity in the derivative is smeared a bit, such that at some time steps these elements are not detected. Note that all approximations are very accurate and close to the exact solution.

The most important improvements are found for the blast-wave and Shu-Osher problem. For the blast waves, the original KXRCF detector and minmod-based TVB indicator detect many elements. However, the new outlier-detection approach combined with these detection variables only selects a few of them, thereby still producing very accurate results.

In the Shu-Osher problem (Figure \ref{fig:Sinek2}) an initial discontinuity is moving to the right, thereby evolving (highly oscillatory) continuous regions and developing new shocks in the left side of the domain. 

The first row of the figure consists of time-history plots of detected troubled cells using the original indicators. Note that both the multiwavelet indicator with $C=0.01$ and the minmod-based TVB indicator with $M=100$ detect the highly-oscillatory region as being discontinuous. In this case, the KXRCF indicator gives more accurate results. For $k=1$ however, the KXRCF indicator only detects the largest discontinuity, and neglects the other three shocks in the left side of the plot, which leads to some spurious oscillations in the approximation.

In the second row of the figure, the time-history plots are shown when the indication vectors are used in the outlier-detection algorithm. All three indication techniques detect the correct regions, and the approximations are as expected (row 3--4 of Figure \ref{fig:Sinek2}). Note that the results are very close to the exact solution: the outlier-detection algorithm is indeed able to replace the problem-dependent parameters in the original indicators.

For $k=1$ and $k=3$, the same behavior is found: the new outlier-detection approach perfectly selects the discontinuous regions in the domain.

\newpage

\newpage

\newpage

\newpage
\subsection{Two-dimensional test}\label{sec:2d}
\begin{figure}[h!]
 \centering
\subfigure[$C=0.05$, $\alpha$]{\includegraphics[scale = 0.41]{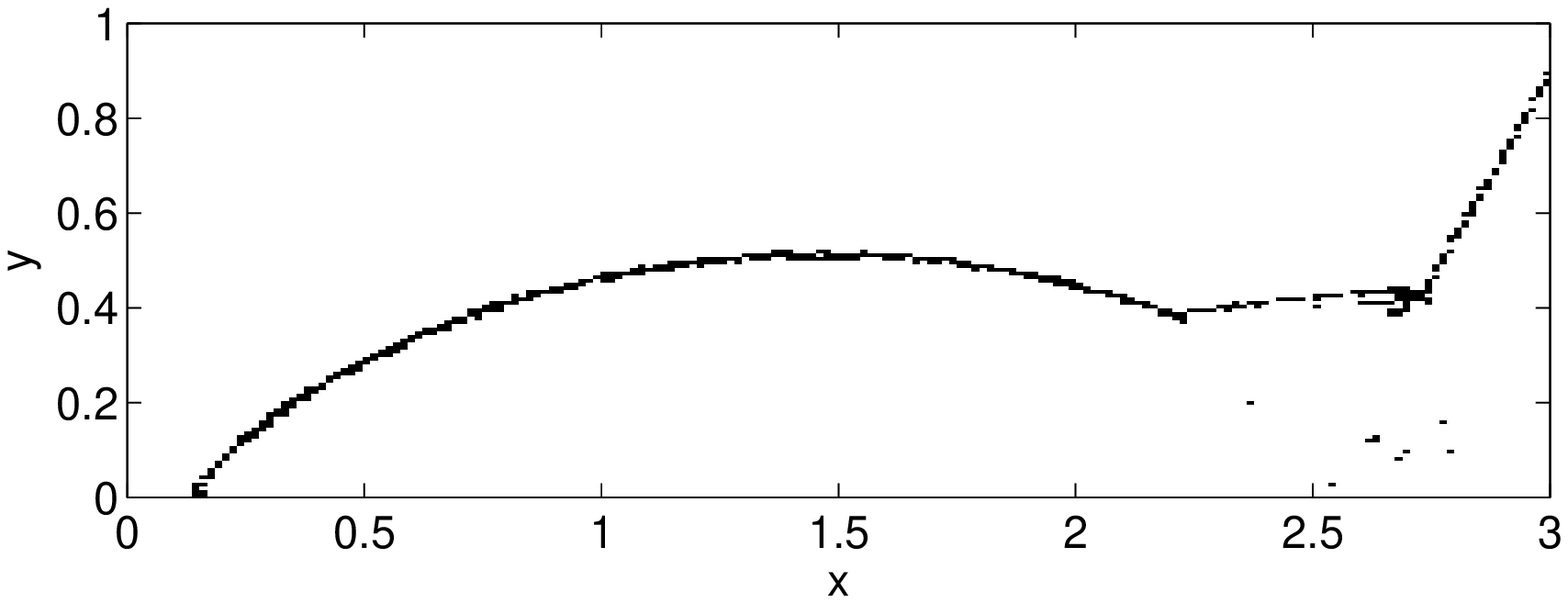}} 
\hfill
\subfigure[Outlier, $\alpha$]{\includegraphics[scale = 0.41]{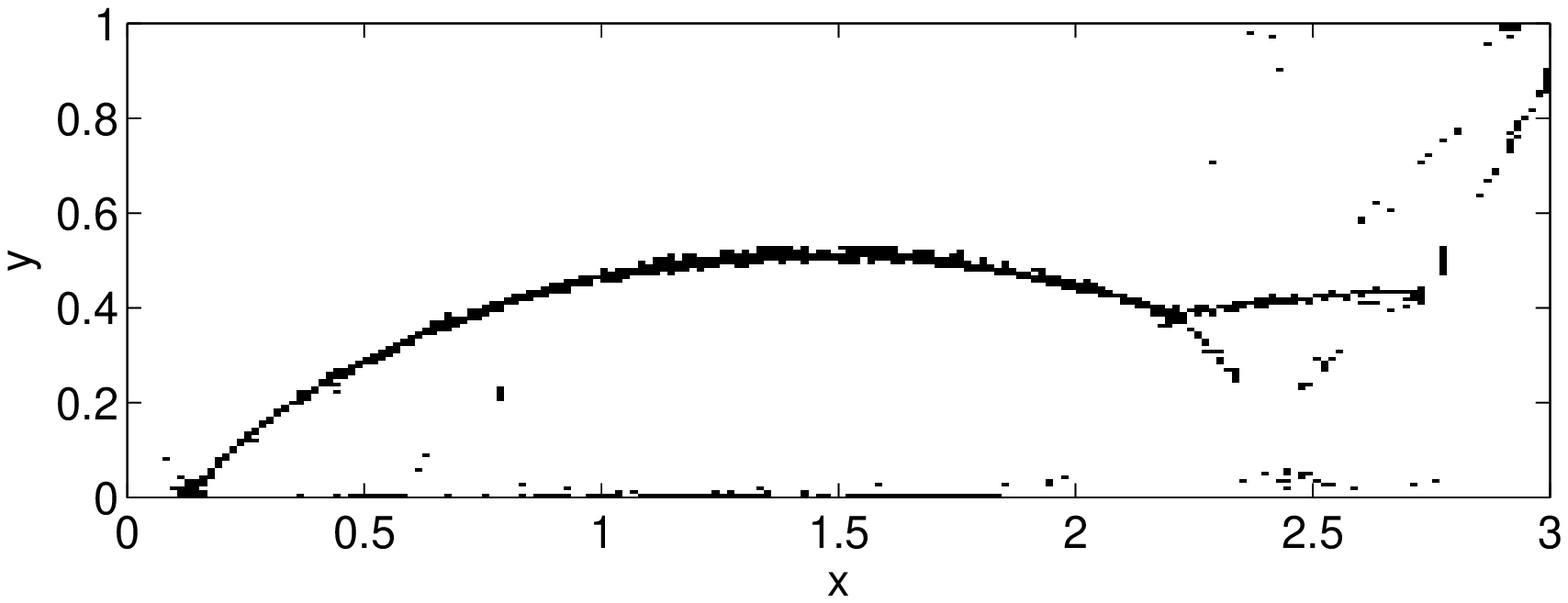}} \\
\subfigure[$C=0.05$, $\beta$]{\includegraphics[scale = 0.41]{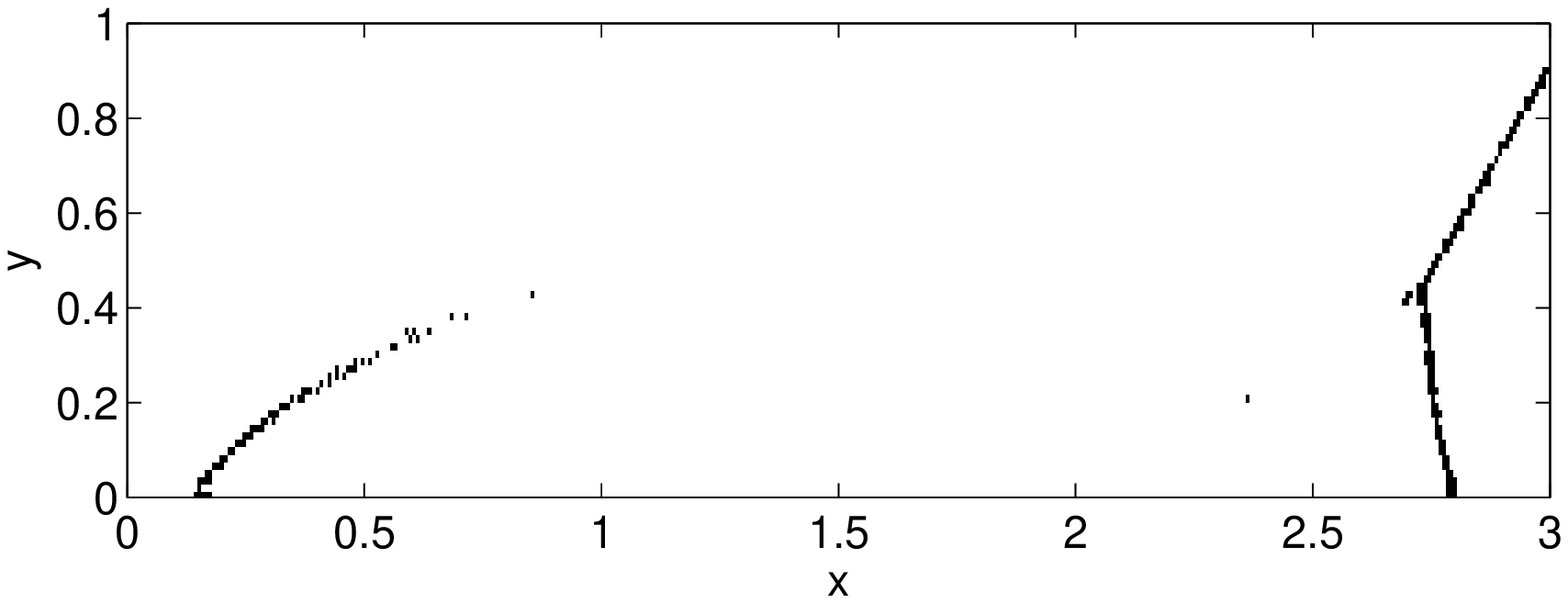}} 
\hfill
\subfigure[Outlier, $\beta$]{\includegraphics[scale = 0.41]{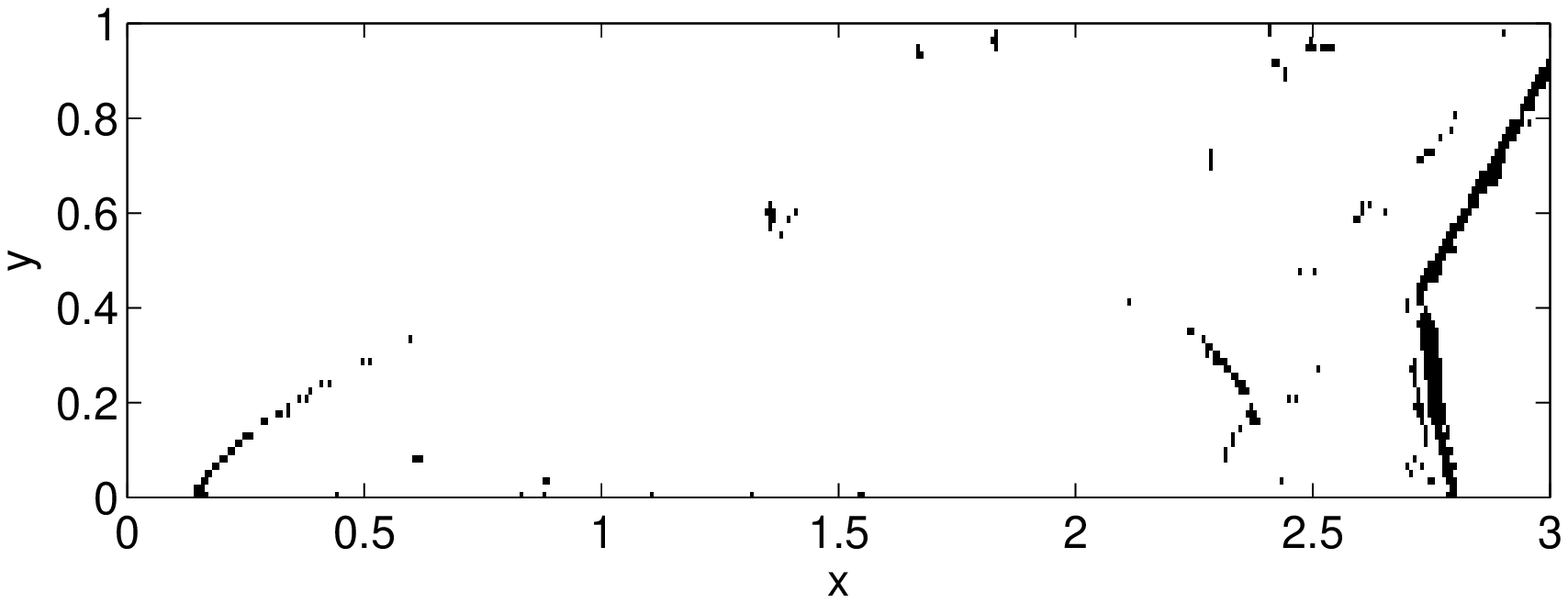}} \\
\flushleft
\subfigure[$C=0.05$, $\gamma$]{\includegraphics[scale = 0.41]{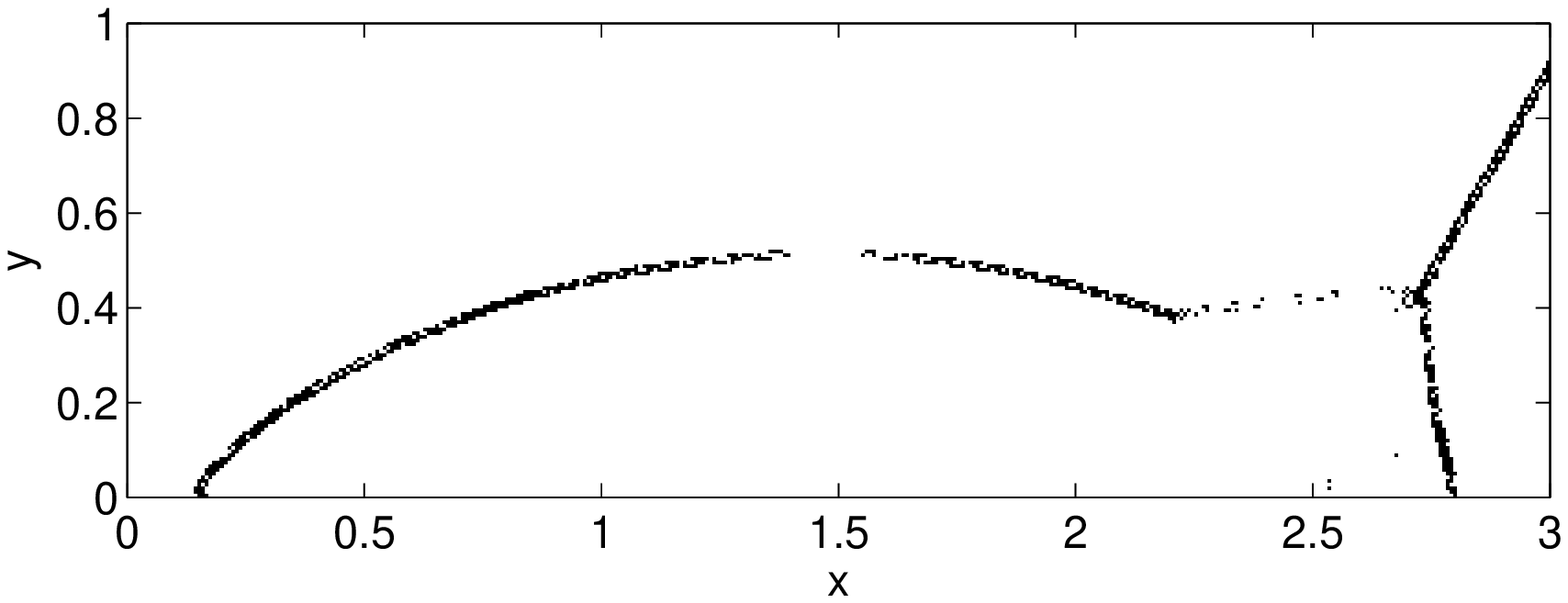}} \\
\centering
\subfigure[$C=0.05$, total detected]{\includegraphics[scale = 0.41]{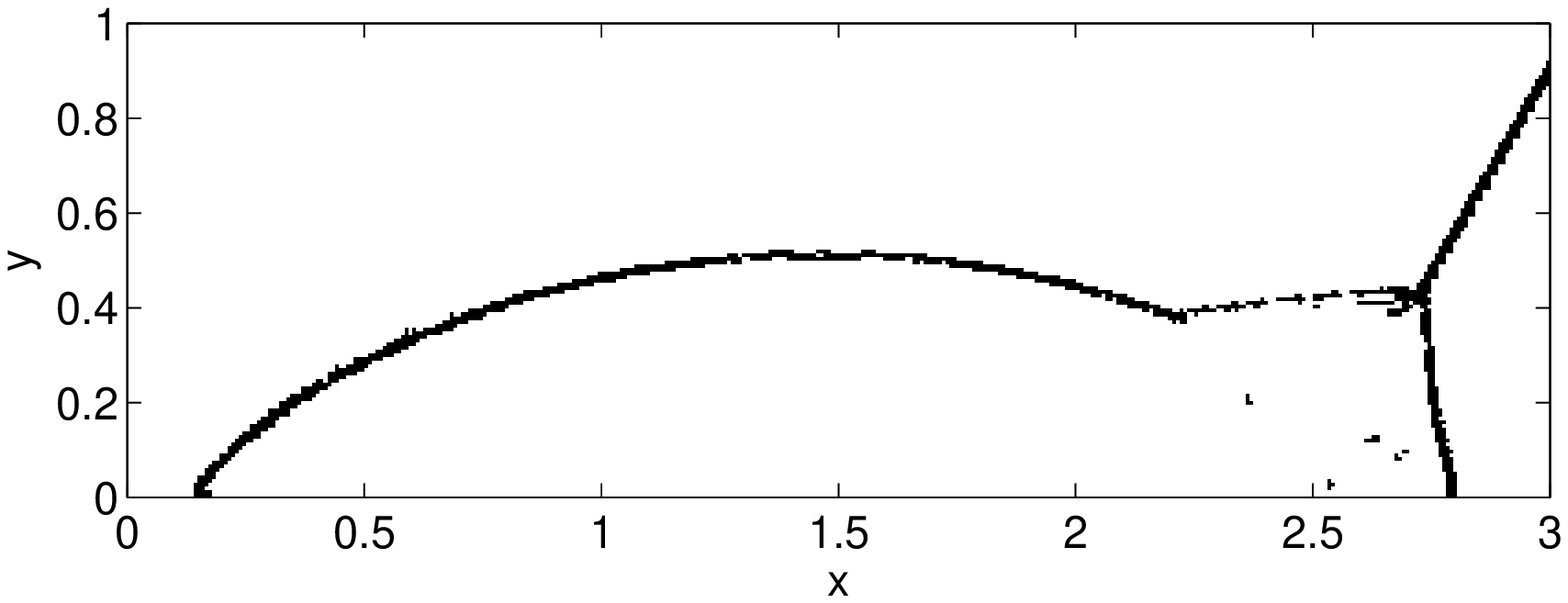}} 
\hfill
\subfigure[Outlier, total detected]{\includegraphics[scale = 0.41]{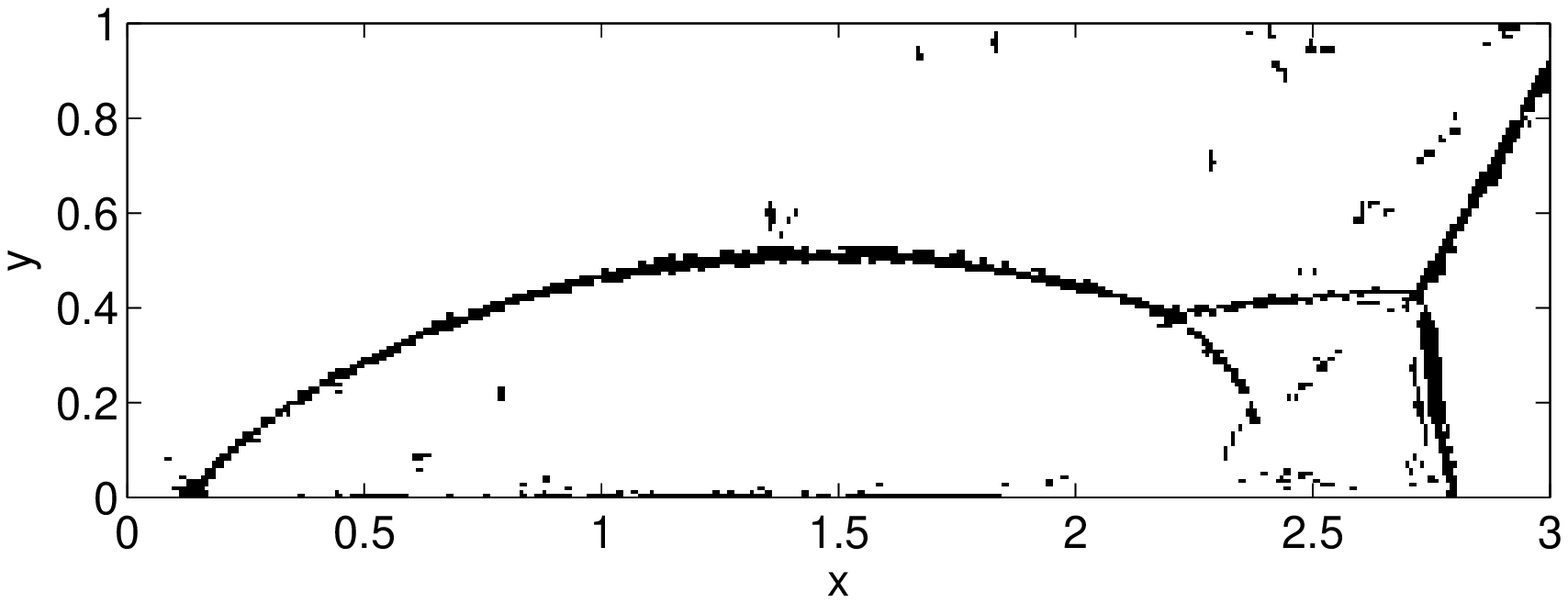}} \\
\vspace{-0.4cm}
\caption{Detected troubled cells at $T=0.2$, double Mach reflection problem, modified multiwavelet troubled-cell indicator, $\Delta x = \Delta y = 1/128$, $k=2$.}\label{fig:dmk2mw}
\end{figure}

\begin{figure}[ht!]
\centering
 \subfigure[Original]{\includegraphics[scale = 0.41]{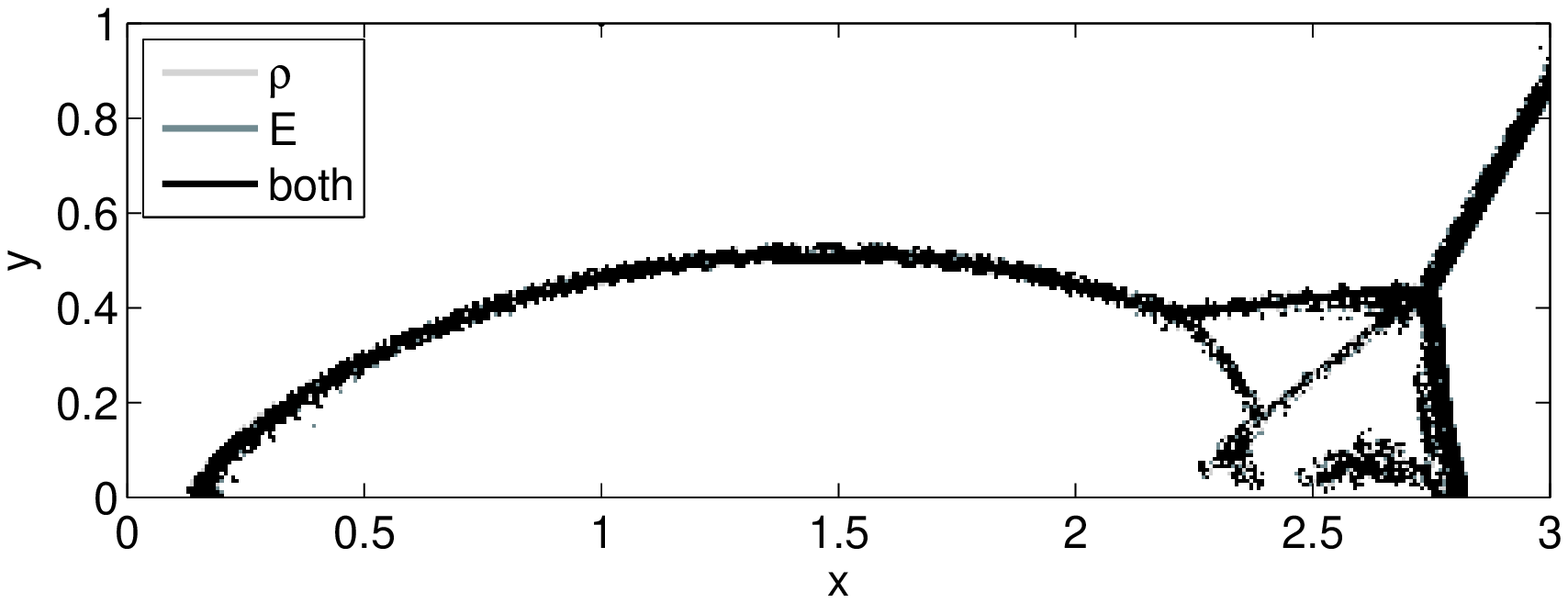}}
 \subfigure[Outlier]{\includegraphics[scale = 0.41]{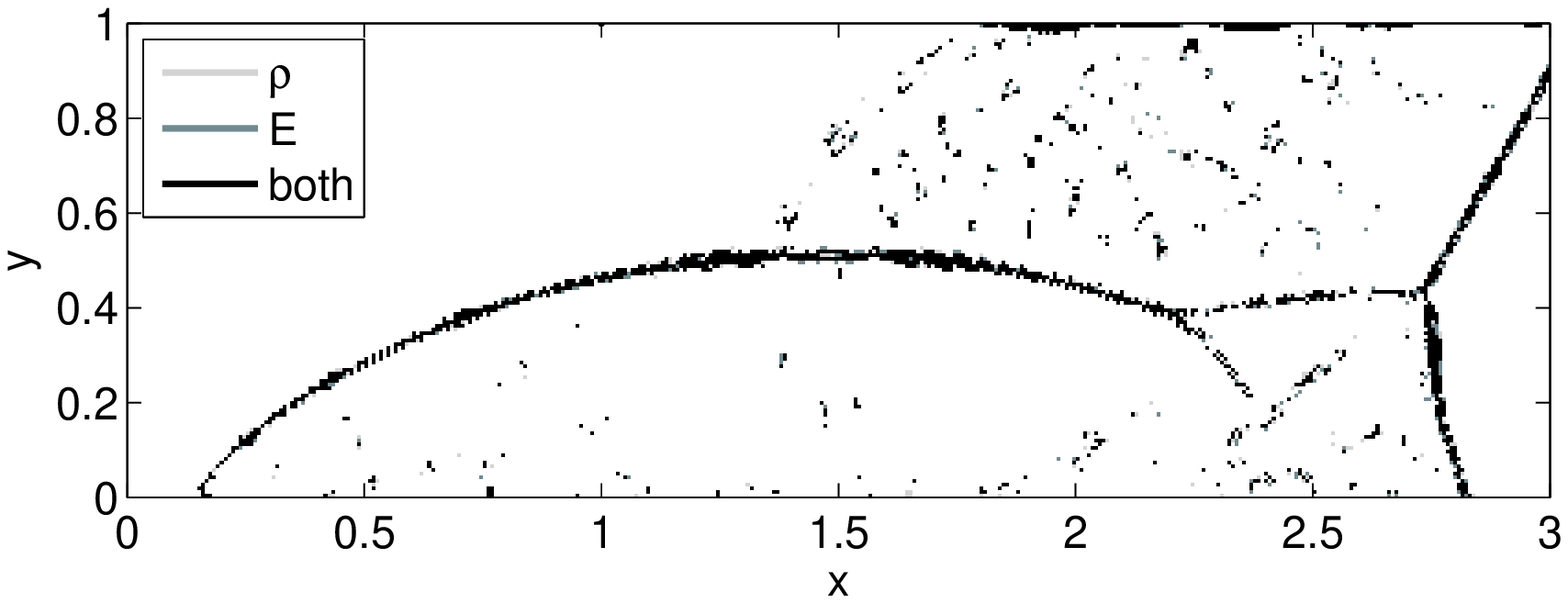}} \\
\caption{Detected troubled cells at $T=0.2$, double Mach reflection problem, KXRCF shock detector, $\Delta x = \Delta y = 1/128$, $k=2$.}\label{fig:dMKXRCFk2}
\end{figure}

\begin{figure}[h!]
\centering
 \subfigure[$M=100$]{\includegraphics[scale = 0.41]{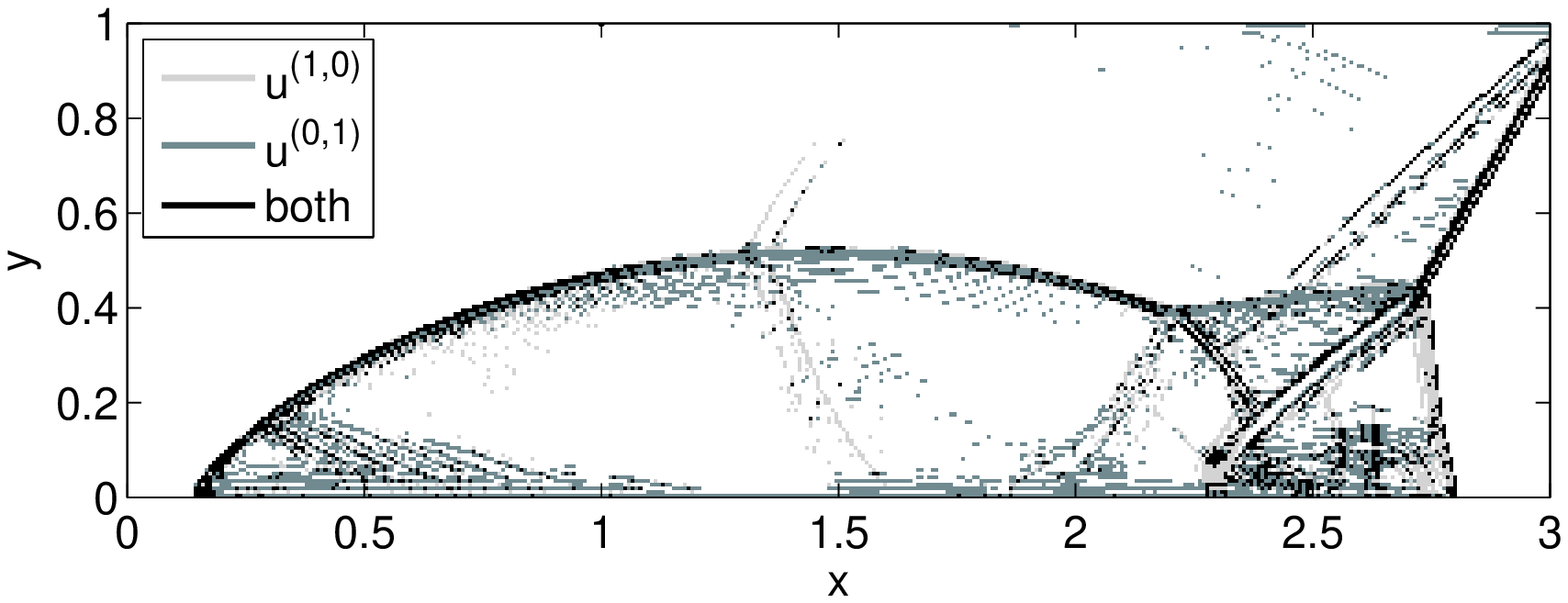}\label{fig:minmodoriginal}}
 \subfigure[Outlier]{\includegraphics[scale = 0.41]{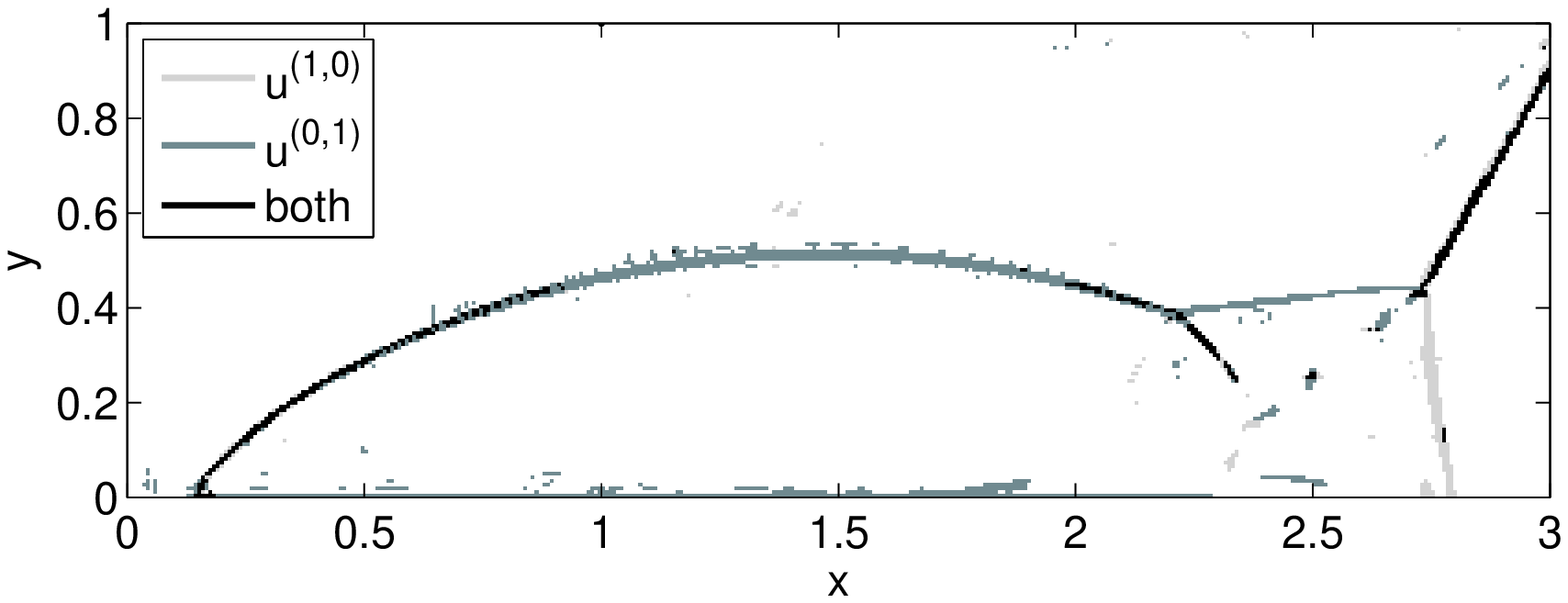}\label{fig:minmodoutlier}} \\
\caption{Detected troubled cells at $T=0.2$, double Mach reflection problem, minmod-based TVB indicator, $\Delta x = \Delta y = 1/128$, $k=1$.}\label{fig:dMminmodk1}
\end{figure}

In two dimensions, we investigate the double Mach reflection of a strong shock \cite{Woo84C}, which satisfies the two-dimensional Euler equations. Again, the original troubled-cell indicators (with optimized parameter) are compared to their outlier-detection approaches. The results for $k=2$ can be compared in Figure \ref{fig:dmk2mw} for the modified multiwavelet troubled-cell indicator, Figure \ref{fig:dMKXRCFk2} for the KXRCF shock detector, and Figure \ref{fig:dMminmodk1} for the minmod-based TVB indicator ($k=1$ only). The spatial domain is split into $2^9 \times 2^7$ rectangular elements, such that $\Delta x = \Delta y = 1/128$. In each figure, the left plots are computed using the original troubled-cell indicators, and the right plots correspond to the outlier-detection approaches.

As mentioned earlier, the multiwavelet technique is able to distinguish between $x$- and $y$-directed discontinuous regions. This is also the case when outlier detection is used. We point out that a sharp detection of the discontinuous region is found. Only a few elements outside the discontinuous region are added, which apparently correspond to discontinuities in derivatives (since jumps in multiwavelet coefficients are found). The approximations at the final time are comparable to the results using the original modified multiwavelet troubled-cell indicator.

The original KXRCF shock detector is compared to the outlier-detection application in Figure \ref{fig:dMKXRCFk2} for $k=2$. The detected troubled cells at the final time using either the original or the outlier-detection approach are similar for $k=1$. For $k=2$ and especially for $k=3$ fewer elements are detected by the outlier-detection scheme. However, more elements are detected in the top region of the domain. This is due to the fact that in this region neighboring jumps across the inflow edges of the element differ substantially from each other.

The results using the minmod-based TVB indicator improve considerably when using outlier detection. In Figure \ref{fig:minmodoriginal}, the detected troubled cells at the final time are shown for the original minmod-based TVB indicator. Note that too many elements are detected: also continuous regions are selected. However, the outlier-detection technique applied to the DG coefficients only selects the correct discontinuity profile (Figure \ref{fig:minmodoutlier}). It should be noticed that this approach detects discontinuities in the $x$- and $y$-direction, since DG coefficients $u_{ij}^{(1,0)}$ are related to the first derivative in the $x$-direction, and $u_{ij}^{(0,1)}$ to the first derivative in the $y$-direction. Fewer elements are detected in this case, and the approximation at time $T=0.2$ is still accurate.

\section{Computational costs}\label{sec:costs}
This section contains a discussion about the computational costs of the outlier-detection algorithm. First, we sort $2^{n-4}$ vectors of length 16 each. We use the 'Selection sort' sorting algorithm, which finds the minimum value of the vector, swaps it with the value in the first position, and repeats these steps for the remainder of the list. The method is of order $\mathcal{O}(N^2)$ time complexity, but it is possible to use a more efficient sorting algorithm (for example of order $\mathcal{O}(N)$) \cite{SorA}. Once the vectors are sorted, we compute the quartiles and outer fences. Outliers are determined by comparing the smallest vector values with $Q_1 - 3(Q_3-Q_1)$ and the biggest components with $Q_3+3(Q_3-Q_1)$. For the smallest values we start with testing whether $d_0^s < Q_1 - 3(Q_3-Q_1)$.  If $d_0^s$ is not an outlier, then there are no other outliers, since $d_1 \geq d_0 \geq Q_1 -3(Q_3-Q_1)$. If $d_0^s$ is an outlier, then we test $d_1^s$ in the same way. Note that (by construction) at maximum $d_0^s, d_1^s$ and $d_2^s$ should be tested (as they are the only possible low outliers). Similarly we test $d_{15}^s$ and (possibly) $d_{14}^s$ and $d_{13}^s$ against $Q_3+3(Q_3-Q_1)$ (depending on the outcome). Finally, the detected outliers in the left half of the local regions are compared with the bounds of the left-neighboring region, and the outliers in the right half are compared with the right-neighboring region.

It should be noticed that this novel method works well on a CPU. The local vectors can also be considered using parallel architectures. However, in that case the costs for communication will be higher, since local information should be distributed along the devices. On the other hand, it also typically results in fewer places where a limiter must be applied.

\bigskip
In Table \ref{tab:computationtime1d} the computational times are shown for the test problems of \S \ref{sec:results}, using either the original or the outlier-detection indication technique. Notice that the computational times using outlier detection are slightly longer than the original times, except for the KXRCF indicator. In that case, the number of detected elements for the original algorithm is much larger than when outlier detection is applied, such that the moment limiter is applied more often. For the rest of the examples, the increase in computational time is on average $2.9\%$, which is reasonable. It should be emphasized that the new method also reduces the number of tests by not having to find a problem-dependent parameter.

\begin{table}[ht!]
 \centering
\begin{tabular}{|l|r|r|r|r|r|r|}
 \hline
& \multicolumn{2}{c|}{Multiwavelets} & \multicolumn{2}{c|}{KXRCF} & \multicolumn{2}{c|}{Minmod} \\
\hline
& Original & Outlier & Original & Outlier & Original & Outlier \\
\hline
Sod        &  0.187 & 0.208  &  0.208 &  0.212 &  0.231 &  0.256 \\
Lax        &  0.263 & 0.280  &  0.299 &  0.290 &  0.329 &  0.366 \\
blast wave & 10.539 & 11.045 & 13.505 & 12.313 & 14.776 & 14.855 \\
Shu-Osher  &  5.683 &  5.845 &  6.520 &  6.512 &  7.669 &  7.973 \\
\hline
\end{tabular}
\caption{Total computation time in seconds for the one-dimensional problems of \S \ref{sec:results}.}\label{tab:computationtime1d}
\end{table}

The total computation times for the double Mach reflection problem are presented in Table \ref{tab:computationtime2d}. Note that the case $k=1$ (minmod-based TVB indicator) is much faster than $k=2$. Also here, the computation time increases, on average by $2.6\%$. Since no tests for parameter finding are needed, the new method will still provide the results much faster.

\begin{table}[ht!]
 \centering
\begin{tabular}{|c|c|c|c|c|c|}
 \hline
\multicolumn{2}{|c|}{Multiwavelets} & \multicolumn{2}{c|}{KXRCF} & \multicolumn{2}{c|}{Minmod} \\
\hline
Original & Outlier & Original & Outlier & Original & Outlier \\
\hline
 312 & 316  & 313   &  324 &  93  & 97  \\
\hline
\end{tabular}
\caption{Total computation time in minutes for the double Mach reflection problem ($k=2$ for multiwavelet and KXRCF indicator, $k=1$ for minmod-based TVB indicator).}\label{tab:computationtime2d}
\end{table}

\vspace{-0.5cm}
\section{Conclusion}\label{sec:conclude}
In this paper, we have introduced a new outlier-detection technique which can be applied to existing troubled-cell indication variables. In this way, problem-dependent parameters are no longer required. We showed the performance of this method for various test problems in one and two dimensions, using the modified multiwavelet troubled-cell indicator, the KXRCF shock detector, and the minmod-based TVB indicator. The results were generally better than the original troubled-cell indicators using an optimized parameter: both the weak and the strong shock regions were detected, whereas smooth regions were not selected. Future work will be to include the local spatial information in the statistical approach, and to extend this to unstructured meshes.

\bigskip
{\bf Acknowledgments:}  The authors gratefully wish to acknowledge the useful comments provided by Dennis den Ouden, Jianxian Qiu, Chi-Wang Shu, Mahsa Mirzagar and Arnold Heemink that helped to shape this work. In particular the authors are grateful to the reviewers for their remarks which helped us to greatly improve this paper.

\newpage

\newpage

\newpage
\bibliographystyle{siam}
\bibliography{References}

\end{document}